\documentclass[aps,rmp]{revtex4}
\pdfoutput=1

\usepackage{amsmath,amssymb}
\usepackage{breakurl}
\usepackage{psfrag}
\usepackage{graphicx}
\usepackage{epsf}
\usepackage{dcolumn}
\usepackage{bm}

\usepackage{color}
\usepackage{environ}
\usepackage{scrextend}
\usepackage[outline]{contour}

\usepackage{hyperref}
\newcommand{\HD}[2]{\href{#1}{{\color{blue} #2}}}

\hypersetup{
     colorlinks   = true,
     citecolor    = blue,
     linkcolor=blue
}

\contourlength{0.2pt} 
\contournumber{4} 

\newcommand{\be}{\begin{equation*}}
\newcommand{\ee}{\end{equation*}}

\definecolor{darkgreen}{rgb}{0,0.7,0}

\newcommand{\fett}[1]{\boldsymbol{#1}}
\newcommand{\ws}{\hspace{0.036cm}} 
\newcommand{\deltavarrho}{\delta\hspace{-0.03cm}\varrho}

\makeatletter
\newcommand{\leqnomode}{\tagsleft@true}
\newcommand{\reqnomode}{\tagsleft@false}
\makeatother

\NewEnviron{myequation}{%
    \begin{equation}
    \scalebox{1.3}{$\BODY$}
    \end{equation}
    }

\NewEnviron{mynewequation}{%
    \begin{equation}
    \scalebox{1.04}{$\BODY$}
    \end{equation}
    }

\NewEnviron{tinyequation}{%
    \begin{equation}
    \scalebox{0.88}{$\BODY$}
    \end{equation}
    }

\NewEnviron{smallequation}{%
    \begin{equation}
    \scalebox{0.93}{$\BODY$}
    \end{equation}
    }

\begin{document}
\title{
Hermann Hankel's \\ ``On the general theory of motion of fluids'' \\
\large   An essay including an English translation of the complete  \contour{black}{\it Preisschrift} from 1861\footnote{Originally published in German as ``Zur allgemeinen Theorie der Bewegung der Fl\"ussigkeiten. 
Eine von der philosophischen Fakult\"at der
Georgia Augusta am 4.\ Juni 1861 gekr\"onte Preisschrift''.}
}

\author{Barbara Villone}
\email{barbara.villone@inaf.it}
\affiliation{INAF, Osservatorio Astrofisico, via Osservatorio, 30,  I-10025 Pino Torinese, Italy}

\author{Cornelius Rampf}
\email{rampf@thphys.uni-heidelberg.de}
\affiliation{Institut f\"ur Theoretische Physik, University of Heidelberg, Philosophenweg 16, D--69120 Heidelberg,
Germany}
\affiliation{Department of Physics, Israel Institute of Technology -- Technion, \mbox{Haifa 32000, Israel}}

\date{\today}

\begin{abstract}
 The present is a companion paper to 
 ``A contemporary look at Hermann Hankel's 1861 pioneering work
 on Lagrangian fluid dynamics''  by Frisch, 
 Grimberg and Villone (2017).
 Here we present the English translation of the 1861 prize manuscript
 from  G\"ottingen University ``Zur allgemeinen Theorie der Bewegung der
 Fl\"ussigkeiten''  (On the general theory of the motion of the fluids)
 of Hermann Hankel (1839--1873),  which was originally
 submitted in Latin and then  translated into German by the Author for 
 publication. We also provide the English translation of two important  reports 
 on the manuscript, one written by Bernhard Riemann and the other by Wilhelm Eduard Weber,
 during the assessment process for the prize. Finally we give a short biography of 
 Hermann  Hankel  with his  complete bibliography.
\end{abstract}


\maketitle
\section{Introductory notes}
\label{s:intro}

\deffootnotemark{\textsuperscript{\thefootnotemark}}\deffootnote{2em}{1.6em}{${}^\thefootnotemark$\hspace{0.0319cm}\enskip}

\setcounter{footnote}{0}

Here we present, with some supplementary documents, a full translation from German of 
the winning essay, originally written in Latin,  by Hermann Hankel,  in response to the ``extraordinary mathematical prize''\footnote{The ordinary prize from the  Philosophical Faculty of the University of G\"ottingen was a philosophical one; extraordinary   prizes, however, could be launched in another discipline.} 
launched on 4th June 1860 by the Philosophical  Faculty of the University 
of G\"ottingen with a deadline of end of March 1861. By request of the Prize Committee to Hankel, the essay was then revised by him and finally published in 1861,  in German,   as a \textit{Preisschrift}.\footnote{Hankel, \hyperlink{Hankel}{1861.}  Hankel  asked for permission to have his revised  \textit{Preisschrift}  published in German.} The Latin manuscript has been very probably returned to the Author and now appears to be  irremediably lost.\footnote{This information has been conveyed to us through the G\"ottingen Archive. An attempt by us to obtain the Latin document from a descendant of Hankel was unsuccessful.} 
The G\"ottingen University Library possesses two copies of the 1861 German original edition of the \textit{Preisschrift}.
Hankel participated in this prize competition shortly after his arrival in the Spring 1860 in G\"ottingen as a 21-years old student in Mathematics, from the University of Leipzig. 
The formulation of the prize (see section~\ref{s:Preisschrift}) highlighted the problem of the equations of fluid motion in Lagrangian coordinates and was presented in memory of Peter Gustav Lejeune-Dirichlet. A \textit{post mortem} paper of his had outlined the advantages of the use of the Lagrangian approach, compared to  the Eulerian one, for the description of the fluid motion.\footnote{Lejeune-Dirichlet, \hyperlink{Dirichlet}{1860}; for more details, see the companion paper: Frisch, Grimberg and Villone, \hyperlink{FGV}{2017}.}

During the decision process for the winning essay, there was an exchange of German-written letters among the committee members.\footnote{G\"ottingen University Archive,   \hyperlink{G\"ottingen University Archive}{1860/1861}.}
Notably, among  all of these letters, the ones from  Wilhelm Eduard Weber (1804--1891)  and from Bernhard Riemann (1826--1866) have been decisive for the evaluation of the essay. Hankel had been the only one to  submit an {essay,  and the discussion among the committee members was devoted  to deciding whether or not} the essay  written by Hankel deserved to win the prize. 
For completeness, we  provide  an English translation of two of these letters, one signed by Bernhard Riemann, the other by Wilhelm Eduard Weber. 
As is thoroughly discussed in the companion paper \mbox{(Frisch, Grimberg and Villone, \hyperlink{FGV}{2017}),} Hankel's \textit {Preisschrift} reveals indeed a truly
deep understanding of Lagrangian fluid mechanics, with innovative use of  variational
methods and differential geometry. Until these days, this innovative work of Hankel has remained apparently poorly known among scholars, with some exceptions. For details and references, see the companion paper.

The paper is organised as follows. 
In section~\ref{s:Preisschrift}, we present the translation of the \textit{Preis-schrift}; this section begins with a preface signed by Hankel, containing the stated prize question together with the decision by the committee members.  For our translation we used   the digitized copy of the \textit{Preisschrift} from  the HathyTrust Digital Library (indicated as a link in the reference). It has been verified by the G\"ottingen University Library that the text that we used is indeed the digitized copy of the 1861 original printed copy.
Section~\ref {s:letters} contains the translation of  two written judgements on Hankel's essay in the procedure of assessment of the prize, one by Weber and the other by Riemann. 
In section~\ref{s:HHpapers} we provide some biographical notes,  together with a full publication list of Hankel.

Let us elucidate our conventions for author/translator footnotes, comments and equation numbering.
Footnotes by Hankel are denoted by an ``A.'' followed by a number (A stands for Author), enclosed in square brackets; translator footnotes
are treated identically except that the letter ``A'' is replaced with ``T'' (standing for translator).
For both author and translator footnotes, we apply a single number count, e.g., [T.1], [A.2], [A.3]. 
Very short translator comments are added directly in the text, and such comments are surrounded by square brackets.
Only a few equations have been numbered by Hankel (denoted by numbers in round brackets). 
To be able to refer to all equations, especially relevant for the companion paper, we have added additional 
equation numbers in the format [p.n], which means the nth equation of \S p.
Finally, we note that the abbreviations ``S.'' and ``Bd.'', which occur in the {\it Preisschrift},  refer respectively to the German words for ``page'' and ``volume''.

\deffootnotemark{\textsuperscript{[T.\thefootnotemark]}}\deffootnote{2em}{1.6em}{[T.\thefootnotemark]\enskip}

\setlength{\footnotesep}{0.24cm}

\section{Hankel's Preisschrift translation}
\label{s:Preisschrift}

\noindent
About the  prize question  set by the philosophical Faculty of Georgia Augusta  
on 4th June 1860\,:\,\footnote{The prize question is in  Latin in 
the {\it Preisschrift}\,:\,``Aequationes generales
motui fluidorum determinando inservientes
duobus modis exhiberi possunt, quorum alter Eulero, alter Lagrangio
debetur. Lagrangiani modi utilitates adhuc fere penitus neglecti
clarissimus Dirichlet indicavit in commentatione postuma 'de problemate
quodam hydrodynamico' inscripta; sed ab explicatione earum uberiore
morbo supremo impeditus esse videtur. Itaque postulat ordo theoriam
motus fluidorum aequationibus Lagrangianis superstructam eamque eo
saltem perductam, ut leges motus rotatorii a clarissimo  Helmholtz
alio modo erutae inde redundent.'' At that time it was common to state prizes
in Latin and to have the submitted essays in the same language.}
\\\\
\indent\indent
  The most useful equations for determining fluid motion may be presented
in two ways, one of  which  is Eulerian, the other one is
Lagrangian.
The illustrious  Dirichlet pointed out in the posthumous unpublished
paper
``On a problem of hydrodynamics'' the almost completely overlooked
advantages  of the
Lagrangian way, but he was  prevented from unfolding this way further  by a
fatal illness. So, this institution asks for a theory of fluid motion  based on
the equations of Lagrange, yielding, at least, the laws of vortex motion
already derived in another way by the illustrious
Helmholtz.\\\\
\noindent
Decision of the philosophical Faculty on  the present manuscript:  \\\\
\indent\indent
 The extraordinary mathematical-physical prize question
about the derivation of laws of fluid motion and, in  particular, of vortical motion, described by the so-called  Lagrangian  equations, 
was answered
by an essay carrying the motto: \,\,{\em The more signs express relationships  in nature, the more useful they are.}\footnote{This motto is in Latin in the  \textit{Preisschrift}: ``Tanto utiliores sunt notae, quanto magis exprimunt rerum relationes''. At that time it was  
common
that an Author  submitting his work for a prize  signed it  anonymously 
with a motto.}\hspace{0.1cm}
This manuscript gives  commendable evidence of the Author's diligence, of his knowledge and  ability in using the methods of computation recently developed  by  contemporary mathematicians.\hspace{0.4cm}In particular $\S\,6$\footnotemark \footnotetext{A number such as  ``$\S\,6$'' refers to a section number in the (lost) Latin version of the manuscript. Numbers are different in the revised German translation.}
contains an elegant method
to establish  the equations of motion for a flow in a fully arbitrary  coordinate system, from a point of view which is  commonly referred to as Lagrangian.\hspace{0.4cm}However, when developing the general laws of vortex motion, 
the Lagrangian approach is  
unnecessarily left aside, and,
 as a  consequence, the various
laws have to be found by 
quite  independent means.\hspace{0.4cm}Also the relation between the 
vortex laws and  the investigations of Clebsch, reported in $\S.\,14.\,15 $, is omitted  by  the Author.\hspace{0.4cm}Nonetheless, 
as his derivation  actually begins from the Lagrangian equations, one may consider  the prize-question as  fully answered by 
this manuscript.\hspace{0.4cm}Amongst the many good things to be found in this essay, the  evoked incompleteness and some 
mistakes due to rushing
which are easy to improve, do not prevent  this Faculty  from assigning the prize to the manuscript,  but the Author would be obliged
to submit a  revised copy of the manuscript, improved according to the suggestions made above before it goes into print.\\
\noindent
\centerline {\rule{5cm}{0.2mm}}
\\\\
\indent\indent On  my request, I have been permitted by the philosophical Faculty  to  have the manuscript, originally submitted in Latin,
 printed in German. --- 
\\\\
\indent\indent
The above mentioned $\S.\,6$ coincides  with $\S.\,5$ of the present essay.
The  $\S.\S\,14. 15$  included what is now in the note of  S.\,45 of  the text [now page~\pageref{footnoteClebsch}];\,\,these $\S.\S.$ were left out, on the one hand, because of  lack of space, and on the other hand because, in the present view,\,\,these $\S.\S.$ are not anymore connected to the rest of the essay.\\

\indent\indent Leipzig,\hspace{0.1cm} September 1861.\\
      $\phantom{C}$ \hfill                     Hermann Hankel. \quad\quad 

\vspace{2cm}

\centerline{\fett{$\S.\,1.$} } 
\vspace{0.3cm}
\noindent
The  conditions and  forces which underlie  most of the natural phenomena  are so  complicated and  various, that it is rarely 
possible  to take them into  account by fully analytical means.{\hspace{0.4cm}Therefore, one should,
for the time being, discard 
those forces and properties which evidently  have  little impact on the  
motion or the changes;
this is done by just  retaining the forces that are of essential and of fundamental importance.\hspace{0.4cm}
Only then, once   this first approximation has been made, 
one may reconsider  the previously disregarded forces and properties, and 
modify the underlying  hypotheses.\hspace{0.4cm}In the case of the  general theory of motion of liquid fluids, it seems therefore advisable to take into account solely the  continuity [of the fluid]  and the constancy of volume, and
to disregard  both viscosity and internal friction which are also present as suggested by  experience; 
it is  advisable  to take into account just the continuity and the elasticity, and consider  the pressure determined as a   function  of the density.
\hspace{0.1cm} Even if, in specific cases, the analytical 
methods have improved and may apply as well  to  more realistic hypotheses, these hypotheses are not yet suitable for a  fundamental  general theory.\\
\indent\indent
The  hypotheses of the general hydrodynamical equations, as  given  by 
\mbox{E\ws u\ws l\ws e\ws r,}  are the following. The fluid is considered as being
composed of an aggregate of molecules, 
which are so  small that one can  find  an infinitely large number of them  in an arbitrarily small space.\hspace{0.4cm}Therefore, the fluid  is considered  as  divided into infinitely small parts of dimensionality  of the first order;  each of these parts is filled with an infinite number of molecules, whose sizes have to be considered as infinitely small quantities of the second order.\hspace{0.4cm}These molecules fill space  continuously and move without friction against each other.\\
\indent \indent
The flow can be  set in motion either by accelerating forces from the individual molecules or  from external  pressure forces.\footnotemark \footnotetext{Literally, Hankel wrote  ``that emanates from the molecules''.} Considering  the nature of fluids, one  easily comes to the conclusion,   also confirmed  by experience,   that the pressure on each fluid element  of the external surface  acts normally  and proportionally to the size of that element. In order to have also a clear definition of the pressure at a given point within the fluid, let us think of an element of an arbitrary surface through this point: the pressure will be normal  to this surface element and proportional to its size, but independent on its direction.\hspace{0.4cm}The difference between liquid and elastic flows\footnote {Nowadays, in this context,  \textit {liquid}  flows would be called \textit{incompressible} flows.}{ is then that, 
for the former,  the density is a constant and,  
for the latter,  the density depends on the pressure, and, conversely,  the pressure depends on the density in a special way.\\
\indent\indent
We shall not insist here on a detailed discussion of these properties as they are usually discussed in the  better textbooks on mechanics.\\
\indent \indent
\deffootnotemark{\textsuperscript{[A.\thefootnotemark]}} \deffootnote{2em}{1.6em}{[A.\thefootnotemark]\enskip}
In order to study the above properties  analytically, two methods have been so far applied, both owed to \mbox{E\ws u\ws l\ws e\ws r.} \hspace{0.4cm}The first method\footnote {First published in the essay: Principes g\'en\'eraux du mouvement des fluides\,\, (Hist.\ de l'Acad.\ de Berlin, ann\'ee, \hyperlink{Euler226}{1755}).}
 considers the velocity of each point  of the fluid as a function of position and time. If  $u,v,w$ are  the velocity components  in  the orthogonal  coordinates $x,y,z$,  then $u,v,w$  are functions of  position $x,y,z$ and time $t$.\hspace{0.4cm}The velocity of the fluid in a given point  is thus the velocity of the fluid particles flowing through that point.\hspace{0.4cm}This method was  exclusively used for the study of   motion of fluids  until  \mbox{D\ws i\ws r\ws i\ws c\ws h\ws l\ws e\ws t}\footnotemark \footnotetext{Untersuchungen  \"uber ein Problem der Hydrodynamik.\,\,  (Crelle's Journal, Bd.\ 55, [actually it is Bd. 58], S.181.  [\hyperlink{Dirichlet2}{1861}]).}
observed that  this method had necessarily the drawback that the absolute space filled by the flow  in general changes over time  and, as a consequence, the coordinates $x,y,z$ are not  entirely independent variables. 
The  method just discussed  seems appropriate if the flow is always  filling the same space, i.e., in the case when the flow is filling the infinite space, or when the motion is stationary.\footnote{D\ws i\ws r\ws i\ws c\ws h\ws l\ws e\ws t  used for this case the first  \mbox{E\ws u\ws l\ws e\ws r}ian method: Ueber einige F\"alle, in denen sich  die Bewegung eines festen K\"orpers in einem incompressibelen fl\"ussigen Medium theoretisch bestimmen l\"asst [Ueber die Bewegung eines festen K\"orpers in einem incompressibeln fl\"ussigen Medium]. (Berichte der Berliner Akademie, \hyperlink{Dirichlet1}{1852}, S.12.).\,\,\, Also \mbox{R\ws i\ws e\ws m\ws a\ws n\ws n\ws}  used the first \mbox{E\ws u\ws l\ws e\ws rian} form in his {\em Essai}: Ueber die Fortpflanzung ebener Luftwellen von endlicher Schwingungsweite. (Bd. VIII d. Abhdlg. der G\"ott. Soc. \hyperlink{Riemann}{1860}).}\hspace{0.1cm}\\
\indent \indent The second, ingenious method of  \mbox{E\ws u\ws l\ws e\ws r}\footnote {De principiis motus fluidorum (Novi comm. acad. sc. Petropolitanae. Bd. XIV. Theil I. pro anno 1759 [\hyperlink{Euler1759}{1770}, in German].) im 6. Capitel: De motu fluidorum ex statu initiali definiendo. S.\,358.}} 
considers the coordinates $x,y,z$ of a flow particle, in any reference system, as a function of time $t$ and of  its position $a,b,c$  at initial time $t=0$.\hspace{0.4cm}This method,  by which the same  fluid particle is  followed during its motion, was reproduced, indeed in a slightly more elegant way, 
\deffootnotemark{\textsuperscript{[A.\thefootnotemark]}} \deffootnote{2em}{1.6em}{[A.\thefootnotemark]\enskip}
by  \mbox{L\ws a\ws g\ws r\ws a\ws n\ws g\ws e},\footnote{M\'ecanique analytique, \'ed III. par Bertrand. Bd.\,II.  S. 250--261.\,\, The first edition of the M\'ecanique is of  the year \hyperlink{Lagrange1}{1788}.} without giving any reference.\hspace{0.1cm} Since it appears that nowadays \mbox{E\ws u\ws l\ws e\ws r}'s work is
 rarely read in detail, this method is considered  due to  \mbox{L\ws a\ws g\ws r\ws a\ws n\ws g\ws e.} \hspace{0.4cm}However, the method was already present in its full completeness in E\ws u\ws l\ws e\ws r's work,  29 years before.\hspace{0.4cm}I owe this interesting, historical note to my honoured  Professor B.  \mbox{R\ws i\ws e\ws m\ws a\ws n\ws n.}\hspace{0.4cm} According to this method, one has thus
\begin{equation*}
x=\varphi_1  (a,b,c,t), \qquad y=\varphi_2  (a,b,c,t), \qquad z=\varphi_3  (a,b,c,t), \tag*{[1.1]}
\end{equation*}
 where $\varphi_1,\varphi_2,\varphi_3 $ are continuous functions of $a,b,c,t$.\,\, We have at initial time $t=0$ that
\begin{equation*} \tag*{[1.2]}
x=a, \hspace{2cm}y=b, \hspace{2cm}z=c
\end{equation*}
and, thus, the following conditions are valid for $\varphi_1,\,\varphi_2,\,\varphi_3$\,:
\be \tag*{[1.3]}
a=\varphi_1  (a,b,c,0), \qquad b=\varphi_2  (a,b,c,0), \qquad c=\varphi_3  (a,b,c,0) .
\ee
Obviously, one can take values of  $t$ so small,  that, according to Taylor's theorem,  $x,y,z$ 
can be expanded in  powers  of $t$.
Since at time $t=0$ we have $x=a,\,\, y=b,\,\, z=c$,  it follows evidently that
\begin{align*} \tag*{[1.4]}
  \begin{aligned}
x\,&= \,a\,+\, A_1t\,+\,A_2 t^2\,+\, \ldots     \\
y\,&= \,b\,+\, B_1t\,+\,B_2 t^2\,+ \hspace{0.09cm} \ldots  \\
z\,&= \,c\,+\, C_1t\,+\,C_2 t^2\,+ \hspace{0.11cm} \ldots 
 \end{aligned} 
\end{align*}
From these equations, one easily finds that at time $t=0$:
\begin{mynewequation} \tag*{[1.5]}
  \begin{aligned}
\frac{dx}{da} &= 1, \,\,\, \,\frac{dx}{db} = 0, \,\,\,\,\frac{dx}{dc} = 0 \\   
\frac{dy}{da} &= 0, \,\,\,\,\frac{dy}{db} = 1, \,\,\,\,\frac{dy}{dc}  = 0 \\
\frac{dz}{da} &= 0, \,\, \,\,\frac{dz}{db} = 0,\,\, \,\,\frac{dz}{dc}  = 1 
\end{aligned} 
\end{mynewequation}
\indent %
In the following, we will try  to give a presentation of the general theory of hydrodynamics based on the  second method.
It will turn out that the second [Lagrangian] form also merits to be preferred over the first one in some cases, because 
the fundamental equations of the former are more closely connected to the customary forms in mechanics.}\\[1cm]
\centerline{\fett{$\S.\,2.$}}\\[0.3cm]
\indent\indent
The two infinitely near particles\\
\centerline {$a,\,\,b,\,\,c$}
and\\
\centerline {$a+{\rm d}a, \,\,b+{\rm d}b, \,\,c+{\rm d}c$}\\\\
after some time $t$, will be at the two points\\\\
\centerline{$x,\,\, y, \,\,z$}
\noindent
and\\
\centerline {$x+{\rm d}x,\,\, y+{\rm d}y,\,\, z+{\rm d}z$,}\\\\
where one needs to put
\begin{mynewequation} \tag*{[2.1]}
 \begin{aligned}
{\rm d}x &=\frac{dx}{da} {\rm d}a + \frac{dx}{db} {\rm d}b + \frac{dx}{dc} {\rm d}c\\
{\rm d}y &=\frac{dy}{da} {\rm d}a + \frac{dy}{db} {\rm d}b + \frac{dy}{dc} {\rm d}c\\
{\rm d}z &=\frac{dz}{da} {\rm d}a + \frac{dz}{db} {\rm d}b + \frac{dz}{dc} {\rm d}c 
\end{aligned} 
\end{mynewequation}

%
\noindent
Let us think  of $a,\,b,\, c$ as  linear --- in general 
non-orthogonal
 --- coordinates; 
then,  we have that ${\rm d}a,\, {\rm d}b,\, {\rm d}c$
are the coordinates of a point  with respect to  a congruent  system of coordinates   $S_0$,  whose 
origin is at $a,\,b,\, c$.\hspace{0.4cm}As
 $x,\,y,\, z$ refer to the same coordinate system as  $a,\,b,\, c$  do, then   
${\rm d}x,\, {\rm d}y,\, {\rm d}z$  are the analogous coordinates with  respect to  a coordinate system, whose origin is at  $x,\,y,\, z$.\hspace{0.4cm}Now, let us think of an infinitely small surface
\be \tag*{[2.2]}
F({\rm d}a,\, {\rm d}b,\, {\rm d}c) = 0
\ee
with reference to $S_0$, then, at time  $t$,  the surface 
will have been moved into another surface 
\be \tag*{[2.3]}
F({\rm d}x,\, {\rm d}y,\, {\rm d}z) = 0
\ee
\noindent
{with reference to  $S$, or, \\
\be \tag*{[2.4]}
F\Big(\frac{dx}{da} {\rm d}a + \frac{dx}{db} {\rm d}b + \frac{dx}{dc} {\rm d}c, \frac{dy}{da} {\rm d}a + \frac{dy}{db} {\rm d}b + \frac{dy}{dc} {\rm d}c, \frac{dz}{da} {\rm d}a + \frac{dz}{db} {\rm d}b + \frac{dz}{dc} {\rm d}c \Big)=0
\ee
We see from this that an infinitely small algebraic surface of  degree $n$ always
 remains of the same degree.\\\\
\indent \indent
Hence, at any time,  infinitely close points on a plane, always stay on a \mbox{plane; and since a}
straight line may be thought of as an intersection of two planes, 
infinitely close points on a line  at a certain time  always stay on  a line.\\\\
\indent\indent
The points within an infinitely small  ellipsoid will always  stay in such an ellipsoid, because  a closed surface cannot transform into a non-closed
 surface and, with the exception of the ellipsoid, all the  surfaces of second degree are not closed.\hspace{0.4cm} The section of a plane and of an infinitely small ellipsoid will thus have to remain always the same.\hspace{0.4cm}Since that section always constitutes an ellipse, so an infinitely small ellipse will always 
remain the same.
\\\\
\indent\indent
The four points:
\begin{mynewequation} \notag
\begin{aligned}
a\,\,\,\,\,\,\,\,\,\,\,\,\,\,\,\,\,\,\,\,\,\,\,\,\,\,\,b\,\,\,\,\,\,\,\,\,\,\,\,\,\,\,\,\,\,\,\,\,\,\,\,\,\,c\,\,\,\,\,\,\,\,\,\\
a+m\,{\rm d}a\,\,\,\,\,\,\,\,\,\,\,b + n\,\,{\rm d}b\,\,\,\,\,\,\,\,c + p\,\,{\rm d}c\\
\,\,\,\,\,\, a+ m'\,{\rm d}a\,\,\,\,\,\,\,\,\,b + n'\,{\rm d}b\,\,\,\,\,\,\, c  +  p'\,{\rm d}c\\
a+m''{\rm d}a \,\,\,\,\,\,\,\,b + n''{\rm d}b\,\,\,\,\,\,\, c  + p''{\rm d}c\\
\end{aligned}
\end{mynewequation}
where $m,\,n,\,p,\,m',\,n',\,p',\,m'',\,n'',\,p''$ are finite numbers, 
at time $t$ will be at the position\\ \\
\centerline{\em \,\,\,\,\,\,x\,\,\,\,\,\,\,\,\,\,\,\,\,\,\,\,\,\,\,\,\,\,\,\,\,\,\,\,\,\,\,\,\,\,\,\,\,\,\,\,\,\,\,\,\,\,\,\,\,\,\,\,\,\,\,\,\,\,\,\,\,\,\,\,\,\,\,\,\,\,\,\,\,\,\,\,\,\,\,\,\,\,y\,\,\,\,\,\,\,\,\,\,\,\,\,\,\,\,\,\,\,\,\,\,\,\,\,\,\,\,\,\,\,\,\,\,\,\,\,\,\,\,\,\,\,\,\,\,\,\,\,\,\,\,\,\,\,\,\,\,\,\,\,\,\,\,\,\,\,\,\,\,\,\,\,\,\,\,\,\,\,\,\,\,\,\,\,\,\,z}
\begin{tinyequation} \nonumber
\begin {aligned}
&\!\!x + \frac{dx}{da} m\,\,{\rm d}a + \frac{dx}{db} n\,\,{\rm d}b + \frac{dx}{dc}  p\,\,{\rm d}c,\,\,\,\,\,\,\,\,\,\,y + \frac{dy}{da} m\,\,{\rm d}a + \frac{dy}{db} n\,\,{\rm d}b + \frac{dy}{dc} p\,\,{\rm d}c, \,\,\,\, z + \frac{dz}{da} m\,\,{\rm d}a + \frac{dz}{db} n\,\,{\rm d}b + \frac{dz}{dc} p\,{\rm d}c \\\\
\!\!&\!\!x+\frac{dx}{da} m'{\rm d}a + \frac{dx}{db} n'{\rm d}b + \frac{dx}{dc}  p'{\rm d}c,\,\,\,\,\,\,\,\,\,\,y + \frac{dy}{da} m'{\rm d}a + \frac{dy}{db} n'{\rm d}b + \frac{dy}{dc} p'{\rm d}c, \,\,\,\,\,\, z + \frac{dz}{da} m'{\rm d}a + \frac{dz}{db} n'{\rm d}b + \frac{dz}{dc} p'{\rm d}c\,\\\\
\!\!&\!\!x+\frac{dx}{da} m''{\rm d}a + \frac{dx}{db} n''{\rm d}b + \frac{dx}{dc}  p''{\rm d}c,\,\, y + \frac{dy}{da} m''{\rm d}a + \frac{dy}{db} n''{\rm d}b + \frac{dy}{dc} p''{\rm d}c, \,\,\, z + \frac{dz}{da} m''{\rm d}a + \frac{dz}{db} n''{\rm d}b + \frac{dz}{dc} p''{\rm d}c 
\end {aligned}
\end{tinyequation}
at time $t$.

\indent\indent
The volume of the tetrahedron $T_0$ whose vertices are  those points at time $t = 0$, is expressed by the determinant\,:\\
\begin{equation*} \tag*{[2.5]}
6T_0=
\begin{vmatrix} m& n & p\\m'& n' & p'\\m''&n''&p''\end{vmatrix} {\rm d}a\,{\rm d}b\,{\rm d}c
\end{equation*}
\indent\indent
The volume of the tetrahedron $T$ whose vertices  at time $t$   are formed by  the same particles, is\,:
\begin{smallequation} \tag*{[2.6]}
6T=
\begin{vmatrix}
\frac{dx}{da} m\,{\rm d}a + \frac{dx}{db} n\,{\rm d}b + \frac{dx}{dc}  p\,{\rm d}c,\,\,\, 
&\frac{dy}{da} m\,{\rm d}a + \frac{dy}{db} n\,{\rm d}b + \frac{dy}{dc} p\,{\rm d}c,\,\,\, &\frac{dz}{da} m\,{\rm d}a + \frac{dz}{db} n\,{\rm d}b + \frac{dz}{dc} p\,{\rm d}c \\\\
\frac{dx}{da}m'{\rm d}a + \frac{dx}{db}n'{\rm d}b + \frac{dx}{dc}p'{\rm d}c,\,\,\,
&\frac{dy}{da} m'{\rm d}a + \frac{dy}{db}n'{\rm d}b + \frac{dy}{dc} p'{\rm d}c, & \frac{dz}{da} m'{\rm d}a + \frac{dz}{db} n'{\rm d}b + \frac{dz}{dc} p'{\rm d}c\\\\
\frac{dx}{da}m''{\rm d}a + \frac{dx}{db} n'' {\rm d}b + \frac{dx}{dc}  p'' {\rm d}c,\,\,\,
&\frac{dy}{da} m'' {\rm d}a + \frac{dy}{db} n'' {\rm d}b + \frac{dy}{dc} p'' {\rm d}c, & \frac{dz}{da} m'' {\rm d}a + \frac{dz}{db} n'' {\rm d}b + \frac{dz}{dc} p'' {\rm d}c 
\end{vmatrix} 
\end{smallequation}
\indent\indent 
By known theorems, this  determinant  can however be 
written as a product of two determinants\,:
\begin{equation} \tag*{[2.7]}
6T=
\begin{vmatrix}
\frac{dx}{da}&\frac{dx}{db}&\frac{dx}{dc} \\\\
\frac{dy}{da}&\frac{dy}{db}&\frac{dy}{dc} \\\\
\frac{dz}{da}&\frac{dz}{db}&\frac{dz}{dc} 
\end{vmatrix}
\begin{vmatrix} m& n & p\\\\m'& n' & p'\\\\m''&n''&p''\end{vmatrix} {\rm d}a\,{\rm d}b\,{\rm d}c
\end{equation}
or
%
\begin{mynewequation} \tag*{[2.8]}
 \text{\small $T=T_0$} \begin{vmatrix}
\frac{dx}{da}\,\,&\frac{dx}{db}\,\,&\frac{dx}{dc} \\\\
\frac{dy}{da}\,\,&\frac{dy}{db}\,\,&\frac{dy}{dc} \\\\
\frac{dz}{da}\,\,&\frac{dz}{db}\,\,&\frac{dz}{dc} 
\end{vmatrix}
\end{mynewequation}
\indent\indent
It  results from the preceding considerations, that all particles which at time $t=0$ are in the  tetrahedron $T_0$,  
will be also in the tetrahedron $T$ at time $t$.\hspace{0.4cm}Let  $\varrho_0$ be
  the mean density of the tetrahedron $T_0$ at time $t=0$, and  $\varrho$ the mean density in $T$ at time $t$, so one has  $T:T_0=\varrho_0:\varrho$ and hence
\begin{mynewequation} \tag*{(1), [2.9]}
 \begin{vmatrix}
\frac{dx}{da}&\frac{dx}{db}&\frac{dx}{dc} \\\\
\frac{dy}{da}&\frac{dy}{db}&\frac{dy}{dc} \\\\
\frac{dz}{da}&\frac{dz}{db}&\frac{dz}{dc} 
\end{vmatrix}
  \text{\small $=$} \,\frac{\varrho_0}{\varrho} \,.
\end{mynewequation}
If the density is constant, the fluid is a liquid flow and thus
\deffootnotemark{\textsuperscript{[T.\thefootnotemark]}} \deffootnote{2em}{1.6em}{[T.\thefootnotemark]\enskip}
\begin{mynewequation} \tag*{(2), [2.10]}
 \begin{vmatrix}
\frac{dx}{da}&\frac{dx}{db}&\frac{dx}{dc} \\\\
\frac{dy}{da}&\frac{dy}{db}&\frac{dy}{dc} \\\\
\frac{dz}{da}&\frac{dz}{db}&\frac{dz}{dc} 
\end{vmatrix}
  \text{\small $=$} \,1.
\end{mynewequation}
\indent\indent
One can reasonably 
refer to these equations as the density equations, more particularly the last one 
as the equation of the constancy of the volume.\\\\
\indent\indent
The values of the functional determinant of $x, y, z$ with respect to $a,b,c$  may be used to develop  a set
of relationships, which are often needed  to pass from  the  \mbox{E\ws u\ws l\ws e\ws rian representation to the other [i.e.,  to the Lagrangian representation; or reciprocally].}
Indeed, if one solves the system of equations\,:
\begin{mynewequation} \tag*{[2.11]}
\begin{aligned}
{\rm d}x=\frac{dx}{da} {\rm d}a + \frac{dx}{db} {\rm d}b + \frac{dx}{dc} {\rm d}c\\
{\rm d}y=\frac{dy}{da} {\rm d}a + \frac{dy}{db} {\rm d}b + \frac{dy}{dc} {\rm d}c\\
{\rm d}z=\frac{dz}{da} {\rm d}a + \frac{dz}{db} {\rm d}b + \frac{dz}{dc} {\rm d}c
\end{aligned} 
\end{mynewequation}
 one has\,:
\begin{mynewequation} \tag*{[2.12]}
\begin{aligned}
\Big(\frac{dy}{db} \frac{dz}{dc} - \frac{dy}{dc}\frac{dz}{db}\Big) {\rm d}x +
\Big(\frac{dz}{db} \frac{dx}{dc} - \frac{dz}{dc}\frac{dx}{db}\Big) {\rm d}y +
\Big(\frac{dx}{db} \frac{dy}{dc} - \frac{dx}{dc}\frac{dy}{db}\Big) {\rm d}z = \frac{\varrho_0}{\varrho} {\rm d}a\\\\
\Big(\frac{dy}{dc} \frac{dz}{da} - \frac{dy}{da}\frac{dz}{dc}\Big) {\rm d}x +
\Big(\frac{dz}{dc} \frac{dx}{da} - \frac{dz}{da}\frac{dx}{dc}\Big) {\rm d}y +
\Big(\frac{dx}{dc} \frac{dy}{da} - \frac{dx}{da}\frac{dy}{dc}\Big) {\rm d}z = \frac{\varrho_0}{\varrho} {\rm d}b\\\\
\Big(\frac{dy}{da} \frac{dz}{db} - \frac{dy}{db}\frac{dz}{da}\Big) {\rm d}x +
\Big(\frac{dz}{da} \frac{dx}{db} - \frac{dz}{db}\frac{dx}{da}\Big) {\rm d}y +
\Big(\frac{dx}{da} \frac{dy}{db} - \frac{dx}{db}\frac{dy}{da}\Big) {\rm d}z = \frac{\varrho_0}{\varrho} {\rm d}c
\end{aligned} 
\end{mynewequation}
where we have substituted the value 
$\varrho_0/\varrho$ for the functional determinant.
\hspace{0.4cm}The comparison of these equations  with\,:
\begin{mynewequation} \tag*{[2.13]}
\begin{aligned}
\frac{da}{dx}\,{\rm d}x + \frac{da}{dy}\,{\rm d}y + \frac{da}{dz}\,{\rm d}z = {\rm d}a\\\\
\frac{db}{dx}\,{\rm d}x + \frac{db}{dy}\,{\rm d}y + \frac{db}{dz}\,{\rm d}z = {\rm d}b\\\\
\frac{dc}{dx}\, {\rm d}x + \frac{dc}{dy}\,{\rm d}y + \frac{dc}{dz}\,{\rm d}z = {\rm d}c
\end{aligned} 
\end{mynewequation}
gives this equation system\,:
\begin{mynewequation}\tag*{(3), [2.14]}
\left.
  \begin{aligned}
\frac{\varrho_0}{\varrho}\frac{da}{dx} & =  \frac{dy}{db} \frac{dz}{dc} - \frac{dy}{dc}\frac{dz}{db},\,
\frac{\varrho_0}{\varrho}\frac{da}{dy} =  \frac{dz}{db} \frac{dx}{dc} - \frac{dz}{dc}\frac{dx}{db}, \,
\frac{\varrho_0}{\varrho}\frac{da}{dz} = \frac{dx}{db} \frac{dy}{dc} - \frac{dx}{dc}\frac{dy}{db}
\\
\frac{\varrho_0}{\varrho}\frac{db}{dx} & = \frac{dy}{dc} \frac{dz}{da} - \frac{dy}{da}\frac{dz}{dc},\,
\frac{\varrho_0}{\varrho}\frac{db}{dy} =  \frac{dz}{dc} \frac{dx}{da} - \frac{dz}{da}\frac{dx}{dc}, \,
\frac{\varrho_0}{\varrho}\frac{db}{dz} = \frac{dx}{dc} \frac{dy}{da} - \frac{dx}{da}\frac{dy}{dc} 
\\
\frac{\varrho_0}{\varrho}\frac{dc}{dx} & =  \frac{dy}{da} \frac{dz}{db} - \frac{dy}{db}\frac{dz}{da},\,
\frac{\varrho_0}{\varrho}\frac{dc}{dy} =  \frac{dz}{da} \frac{dx}{db} - \frac{dz}{db}\frac{dx}{da}, \,
\frac{\varrho_0}{\varrho}\frac{dc}{dz} = \frac{dx}{da} \frac{dy}{db} - \frac{dx}{db}\frac{dy}{da}
 \end{aligned} \qquad  \right\}
\end{mynewequation}
\indent\indent
If the equation (1) is differentiated with respect to one of the independent variables $a,b,c$, and these differentiations are indicated 
with  $\delta$, one obtains
\begin{align} 
\Big(\frac{dy}{db} \frac{dz}{dc} - \frac{dy}{dc}\frac{dz}{db}\Big) \frac{d \delta x}{da} + 
\Big(\frac{dy}{dc} \frac{dz}{da} - \frac{dy}{da}\frac{dz}{dc}\Big) \frac{d \delta x}{db}  +
\Big(\frac{dy}{da} \frac{dz}{db} - \frac{dy}{db}\frac{dz}{da}\Big) \frac{d \delta x}{dc} \nonumber \\ \nonumber \\ 
+ \Big(\frac{dz}{db} \frac{dx}{dc} - \frac{dz}{dc}\frac{dx}{db}\Big)  \frac{d \delta y}{da} + 
\Big(\frac{dz}{dc} \frac{dx}{da} - \frac{dz}{da}\frac{dx}{dc}\Big)  \frac {d \delta y}{db} +
\Big(\frac{dz}{da} \frac{dx}{db} - \frac{dz}{db}\frac{dx}{da}\Big) \frac {d \delta y}{dc}  \nonumber \\ \nonumber \\
+ \Big( \frac{dx}{db} \frac{dy}{dc} - \frac{dx}{dc}\frac{dy}{db}\Big) \frac{d \delta z}{da} +
\Big(\frac{dx}{dc} \frac{dy}{da} - \frac{dx}{da}\frac{dy}{da}\Big) \frac{d \delta z}{db} +
\Big(\frac{dx}{da} \frac{dy}{db} - \frac{dx}{db}\frac{dy}{da}\Big) \frac{d \delta z}{dc}\,\,\,\,\,\,\,\hspace{-0.4cm} \nonumber \\ \nonumber \\
= -\frac{\varrho_0}{\varrho} \frac{\delta \varrho}{\varrho}\,\,\,\,\,\,\,\,\,\,\,\,\,\,\,\,\,\,\,\,\,\,\,\,\,\,\,\,\,\,\,\,\,\,\,\,\,\,\,\,\,\,\,\,\,\,\,\,\,\,\,\,\,\,\,\,\,\,\,\,\,\,\,\,\,\,\,\,\,\,\,\,\,\,\,\,\,\,\,\,\,\,\,\,\, \tag*{[2.15]}
\end{align}
and, by equations~(3):
\begin{equation*} \tag*{[2.16]}
\frac{d \delta x}{da}\frac{da}{dx} + \frac{d \delta x}{db}\frac{db}{dx} + \frac{d \delta x}{dc}\frac{dc}{dx} + \frac{d \delta y}{da}\frac{da}{dy}  + \frac{d \delta y}{db}\frac{db}{dy}  + \frac{d \delta y}{dc}\frac{dc}{dy}  + \frac{d \delta z}{da}\frac{da}{dz}  + \frac{d \delta z}{db}\frac{db}{dz}   + \frac{d \delta z}{dc}\frac{dc}{dz} + \frac{\delta \varrho}{\varrho}= 0
\end{equation*}
or 
\begin{equation*}\tag*{(4), [2.17]}
\frac{d \delta x}{dx} +  \frac{d \delta y}{dy}+ \frac{d \delta z}{dz}+  \frac{\delta \varrho}{\varrho}= 0
\end{equation*}
\indent
If by $\delta$, one understands the differentiation with respect to $t$,  one has:
\begin{equation*}  \tag*{[2.18]}
\varrho \Big(\frac{d \frac{dx}{dt}}{dx} + \frac{d \frac{dy}{dt}}{dy}+ \frac{d \frac{dz}{dt}}{dz}\Big) + \frac{{ \rm d} \varrho} {{\rm d} t}  = 0
\end{equation*}
Since $t$ appears not only explicitly in $\varrho$, but also implicitly  in $\varrho$ through its dependence  on $x$,$y$,$z$, one has
\begin{equation*} \tag*{[2.19]}
\frac{{ \rm d} \varrho} {{\rm d} t} = \frac{d\varrho}{dt}   + \frac{d\varrho}{dx} \frac{dx}{dt} +  \frac{d\varrho}{dy} \frac{dy}{dt} + \frac{d\varrho}{dz} \frac{dz}{dt} 
\end{equation*}
If one sets
\begin{equation*} \tag*{[2.20]}
u = \frac{dx}{dt }, \,\, v= \frac{dy}{dt }, \,\, w = \frac{dx}{dt }
\end{equation*}
one  thus has\,:
\begin{equation*} \tag*{[2.21]}
\varrho\Big( \frac{du}{dx}+ \frac{dv}{dy} + \frac{dw}{dz}\Big) +  { \frac{d\varrho}{dt} } +  \frac{d\varrho}{dx} u +  \frac{d\varrho}{dy} v + \frac{d\varrho}{dz} w = 0 
\end{equation*}
or
\begin{equation*}\tag*{(5), [2.22]}
\frac{d\varrho}{dt} +\frac {d (\varrho u)}{dx}  +  \frac {d (\varrho v)}{dy} +\frac {d (\varrho w)}{dz} = 0
\end{equation*}
This is the form in which   \mbox{E\ws u\ws l\ws e\ws r} first presented  the density equation.\hspace{0.4cm}For the case when $\varrho$ is constant in time,  
in this first form of dependence,  one obtains as the constancy-of-volume equation\,:
\deffootnotemark{\textsuperscript{[A.\thefootnotemark]}} \deffootnote{2em}{1.6em}{[A.\thefootnotemark]\enskip}
\begin{equation*}\tag*{(6), [2.23]}
\frac{du}{dx} + \frac{dv}{dy} + \frac{dw}{dz} = 0   
\end{equation*}
 L\ws a\ws g\ws r\ws a\ws n\ws g\ws e\footnote{M\'ecanique analytique, \,\, Bd.\ I., S. 179--183, Bd.\ II, S. 257--261. [First edition, \hyperlink{Lagrange1}{1788}].} 
treats the relation of these equations  with  (2) quite extensively; but this connection seems
to be a special case of%
\deffootnotemark{\textsuperscript{[A.\thefootnotemark]}} \deffootnote{2em}{1.6em}{[A.\thefootnotemark]\enskip}%
a theorem by  \mbox{J\ws a\ws c\ws o\ws b\ws i}},\footnote{C.\ G.\ J.\ Jacobi, Theoria novi multiplicatoris systemati aequationum differentialium vulgarium applicandi. \,\,\, Crelle's Journal, Bd.\ 27. S.\ 209 [\hyperlink{Jacobi}{1844}].}
 which for three variables is [actually] fully included in equations (1) and (5).\\
\indent\indent
If an integral:
\begin{equation*} 
\iiint f(x,y,z) \,\,\varrho \,\, {\rm d}x\,\, {\rm d}y\,\, {\rm d}z 
\end{equation*}
which is extended over the whole fluid mass, is transformed  into an integral over $a,\,b,\,c$,  one has:
\begin{equation*} \tag*{[2.24]}
\iiint f(x,y,z) \,\,\varrho \,\, {\rm d}x\,\, {\rm d}y\,\, {\rm d}z = \iiint f(a,b,c) \,\,\varrho\,\, {\rm d}a\,\, {\rm d}b\,\, {\rm d}c 
\begin{vmatrix}
\frac{dx}{da}\,\,&\frac{dx}{db}\,\,&\frac{dx}{dc} \\\\
\frac{dy}{da}\,\,&\frac{dy}{db}\,\,&\frac{dy}{dc} \\\\
\frac{dz}{da}\,\,&\frac{dz}{db}\,\,&\frac{dz}{dc}  
\end{vmatrix}
\end{equation*}
where also the second integral has to be  extended over all particles,  and $f(a,\,b,\,c)$  is the function into which
$f(x,\,y,\,z) $   transforms by substituting $a,\,b,\,c$ for  $x,\,y,\,z$.\hspace{0.4cm}
Thus,  from  the density equation (1)  follows\,:
\begin{equation*} \tag*{[2.25]}
\iiint f(x,y,z) \,\,\varrho \,\, {\rm d}x\,\, {\rm d}y\,\, {\rm d}z = \iiint f(a,b,c) \,\,\varrho_0\,\, {\rm d}a\,\, {\rm d}b\,\, {\rm d}c
\end{equation*}
---  an important transformation.
\\[1cm]
\centerline{\fett{$\S.\,3.$}}\\[0.3cm]
\indent\indent
Despite the  
complete analogy of these equations for the equilibrium and 
motion of liquid fluids with the equations for elastic fluids, 
\deffootnotemark{\textsuperscript{[T.\thefootnotemark]}} \deffootnote{2em}{1.6em}{[T.\thefootnotemark]\enskip}%
 there is still  an essential difference with regard to 
their derivation.
\hspace{0.4cm}For this reason, these cases have to be treated separately.\\\\
\indent\indent
Before developing the equations of motion, 
it will be convenient  to study  
more in detail the equilibrium conditions,
at first for liquid fluids.\\\\
\indent\indent
If $X, Y, Z$  are the accelerating forces  in the direction of the coordinates axes, acting on the  point $ x,\, y,\, z$, then it follows easily from the principle of virtual velocities that\,:
\begin{equation*} \tag*{[3.1]}
\iiint \big[\varrho\,\big(X \delta x + Y\delta y + Z \delta z \big) + p \delta L \big]\,\, {\rm d}x\,\,{\rm d}y\,\,{\rm d}z\,=\,0
\end{equation*}
in which $L=0$ gives the density equation for incompressible fluids; $\delta L$ is 
its relative variation 
 corresponding to the  variations of coordinates $\delta x, \delta  y, \delta z$, and $p$ is a not yet determined  quantity. \,\,The integral has to  be extended over all parts of the continuous flow.\hspace{0.4cm}From (4) in $\S.\,2$, since $\delta \varrho=0$, we have
\begin{equation*} \tag*{[3.2]}
\delta L = \frac{d \delta x}{dx} +   \frac{d \delta y}{dy} +  \frac{d \delta z}{dz}
\end{equation*}
And thus the previous integral becomes:
\begin{equation} \tag*{(1), [3.3]}
\iiint \Big[ \varrho \Big( X\delta x + Y \delta y + Z \delta z \Big) + p\Big(\frac{d \delta  x}{dx} + \frac{d \delta y} {dy} +\frac{d \delta z}{dz} \Big ) \Big] {\rm d}x\,{\rm d}y\,{\rm d}z = 0
\end{equation}
Integrating by parts, one finds\,:
\begin{align} \tag*{[3.4]}
\begin{aligned}
\iiint p \frac{d\delta x}{dx} {\rm d}x\,{\rm d}y\,{\rm d}z = \iint p\, \delta x\,{\rm d}y\,{\rm d}z  - \iiint \delta x \frac{dp}{dx}\, {\rm d}x\, {\rm d}y\, {\rm d}z \nonumber \\
\iiint p \frac{d\delta y}{dy} {\rm d}x\,{\rm d}y\,{\rm d}z = \iint p\, \delta y\,{\rm d}z\,{\rm d}x  - \iiint \delta y \frac{dp}{dy}\, {\rm d}x\, {\rm d}y\, {\rm d}z   \\
\iiint p \frac{d\delta z}{dz} {\rm d}x\,{\rm d}y\,{\rm d}z = \iint p\, \delta z\,{\rm d}x\,{\rm d}y  - \iiint \delta z \frac{dp}{dz}\, {\rm d}x\, {\rm d}y\, {\rm d}z \nonumber 
\end{aligned} 
\end{align}
where the double integrals extend over the surface of the fluid mass.\hspace{0.4cm}
Thus, one has for the  equation of  the principle of virtual velocity \,:
\begin{smallequation} \tag*{[3.5]}
0=\iiint\!\!\Big [\!\big(\varrho X \!-\! \frac{dp}{dx} \big) \delta x +  \big(\varrho Y - \frac{dp}{dy} \big) \delta y +  \big(\varrho Z - \frac{dp}{dz} \big) \delta z  \! \Big] {\rm d}x {\rm d}y {\rm d}z + \!\int\! p ({\rm d}x \cos \alpha + {\rm d}y \cos \beta + {\rm d}z \cos \gamma) {\rm d}\omega
\end{smallequation}
where ${\rm d}\omega$ is an element of the external surface, and $\alpha$,\,$\beta$,\,$\gamma$ indicate the angles between the normal to the element $d\omega$ and the coordinates axes.\hspace{0.4cm}From these equations, it  follows that the
equilibrium conditions are\,:
\begin{equation} \tag* {(2), [3.6]}
\int p (\delta x \,\cos \alpha + \delta y\, \cos \beta + \delta z\, \cos \gamma) {\rm d} \omega =0
\end{equation}
\begin{equation} \tag* {(3), [3.7]}
 \frac{dp}{dx} = \varrho X,  \,\,\, \frac{dp}{dy} = \varrho Y, \,\,\,  \frac{dp}{dz} = \varrho Z  
\end{equation}
 The last three equations 
require that the components $X,Y,Z$ be the 
differential quotients of 
an arbitrary function $V$ with respect to $x,y,z$;  thus
\begin{equation*} \tag*{[3.8]}
X=\frac{dV}{dx}, \,\, Y=\frac{dV}{dy},\,\,Z=\frac{dV}{dz}
\end{equation*}
so that one has\,:
\begin{equation*} \tag*{[3.9]}
p = \varrho V + c
\end{equation*}
where $p$ is determined up to  an arbitrary   constant $c$.\\\\
\indent\indent
Instead  of equations $(3)$, one can also write\,:
\begin{mynewequation} \tag*{[3.10]}
 \begin{aligned}
(p_{x + {\rm d}x}- p_x)\,{\rm d}y\, {\rm d}z-\varrho X\, {\rm d}x\, {\rm d}y\, {\rm d}z=0\\
(p_{y + {\rm d}y}- p_y)\,{\rm d}z\, {\rm d}x-\varrho Y\, {\rm d}x\, {\rm d}y\, {\rm d}z=0\\
(p_{z + {\rm d}z}-p_z)\, {\rm d}x\, {\rm d}y-\varrho Z\,\,{\rm d}x\, {\rm d}y\, {\rm d}z=0\\
\end{aligned} 
\end{mynewequation}
from which it is apparent
that  $p$ is the pressure at the point $x,y,z$, which acts against the given 
accelerating forces.\hspace{0.4cm}This pressure is determined up to an additive constant by $p=\varrho  V + c$\,; it is fully determined  if its value is given at any point.\hspace{0.4cm} 
Suppose now that in addition to the acceleration of the fluid particles there are also pressure forces acting on the external surface. We then find as an equilibrium condition that at each point of the external surface, these  pressure forces must  be equal and opposite  to the pressure  $p = \varrho V + c$.\\\\
\indent\indent
In  equation (2),  the variations $\delta x, \delta y, \delta z $ depend on  certain conditions which  result from the nature of the walls. If in special cases the actual meaning of  equation  (2) is specified more precisely, then  its mechanical necessity becomes manifest. Here we want to discuss just a few such cases.
\\\\
\indent\indent
Let us assume that  a part of the external surface be free and  that the same forces act on all parts of it. One can set  the pressure to zero in the points of 
that free surface, so that $p$, the difference between the pressure in a certain point and the pressure in a point of the free surface,  
be exactly determined in all other remaining points of the fluid.\\\\
\indent\indent
However, since $ \delta x, \delta y, \delta z$  are evidently arbitrary in a free surface, it follows that, if $(2)$ has to be satisfied, we must have $p=0$.\\\\
\indent\indent
For  fluid parts that are  lying on a fixed wall, it is evident that no  motion normal to the wall surface  can take place.\hspace{0.4cm}The normal component of the motion is obviously\,:
\begin{equation*}
\delta x\,\cos \alpha + \delta y \,\cos \beta  + \delta z\,\cos \gamma
\end{equation*}
and, since,  this must vanish,  equation (2) is indeed satisfied.\\\\
\indent\indent
In the points which are on a moving  wall, one can set
\begin{equation*} \tag*{[3.11]}
\delta x = \delta {x'} + \delta\xi,\,\, \delta y = \delta {y'} + \delta \eta, \,\,\delta z = \delta {z'} + \delta \zeta
\end{equation*}
where 
$\delta {x'}, \,\,\delta {y'} ,\,\,\delta {z'}$ are  the motions relative to the wall and $\delta\xi,\,\, \delta \eta,\,\,\delta\chi$ are the motions of the fluid particles simultaneously in motion with the wall.\hspace{0.4cm}Therefore, one has instead of equation (2)
\begin{equation*} \tag*{[3.12]}
0=\int p (\delta x' \,\, \cos\alpha + \delta y'\,\,  \cos\beta + \delta z'\,\,  \cos \gamma)  {\rm d} \omega+ \int p (\delta\xi\,\, \cos\alpha + \delta \eta\,\, \cos\beta + \delta \zeta\,\,	\cos \gamma)  {\rm d} \omega
\end{equation*}
but, since  no motion may happen against the wall, one has\,:
\begin{equation*} \tag*{[3.13]}
\delta x' \,\, \cos\alpha + \delta y'\,\, \cos\beta+ \delta z'\,\, \cos\gamma=0
\end{equation*}
 and therefore
\begin{equation*} \tag*{[3.14]}
\int p (\delta \xi \,\,\cos \alpha + \delta \eta\,\, \cos \beta + \delta \zeta \,\, \cos \gamma) {\rm d} \omega = 0
\end{equation*}
which amounts to  setting
\begin{equation*} \tag*{[3.15]}
\int \big (p_x \delta \xi + p_y \delta \eta + p_z \delta \zeta \,) {\rm d} \omega =0
\end{equation*}
provided that  $p_x, p_y, p_z$ are the components of the pressure with respect to  the coordinates axes.\hspace{0.4cm}However this  integral is the 
equilibrium condition of a body, on which  act external surface forces $p_x, p_y, p_z$, where 
$ \delta \xi, \delta \eta, \delta \zeta$ indicate  the variations, that the body can have, under the given circumstances.\hspace{0.4cm}Actually, also in this case, equation (2)  is needed by the nature of things.%
\deffootnotemark{\textsuperscript{[A.\thefootnotemark]}}\deffootnote{2em}{1.6em}{[A.\thefootnotemark]\enskip}%
\footnote {Cf.\ M\'ec.\ analy.\  Bd. I, S.\ 193--201. [First edition,  \hyperlink{Lagrange1}{1788}].}
\\\\
\indent\indent
We now discuss the case of elastic fluids; here the relation $L=0$ is not valid anymore.\hspace{0.4cm}In that case one has to 
consider also the forces due to elasticity of the fluid in addition  to the accelerating and  external pressure forces.\hspace{0.4cm}Let  $p$ be the pressure at a point $x, y, z$, then this tends to reduce the 
 volume of the element ${\rm d}x,\, {\rm d}y,\, {\rm d}z$\,; 
the momentum by this force is thus  $p \delta  ({\rm d}x \, {\rm d}y\, {\rm d}z)$; in order to get a different expression for  $ \delta ({\rm d}x \, {\rm d}y \,{\rm d}z)$, we note that $\rho\,{\rm d}x\,{\rm d}y\,{\rm d}z$, being the mass of an element, is always the same; thus we have $\delta (\rho\,{\rm d}x\, {\rm d}y\,{\rm d}z)=0$;  herefrom it follows that\:
\begin{equation*} \tag*{[3.16]}
\varrho \delta ( {\rm d} x \,\, {\rm d}y\,\,  {\rm d}z) +  {\rm d}x\,\,  {\rm d}y\,\,  {\rm d}z\,\,  \delta\varrho=0
\end{equation*}
and thus 
\begin{equation*} \tag*{[3.17]}
\delta ( {\rm d}x\,\, {\rm d}y\,\,  {\rm d}z) = - \frac{\delta \varrho}{\varrho}  {\rm d}x\,\, {\rm d}y\,\,  {\rm d} z
\end{equation*}
or, by  equation $(4)$ in \S$.2$,
\begin{equation*} \tag*{[3.18]}
\delta ( {\rm d}x\,\, {\rm d}y\,\,  {\rm d}z) \,\,=\,\,\Big(\frac{d \delta x}{dx}\,\,+  \frac{d \delta y}{dy}\,\,+\frac{d \delta z}{dz}\,\,\Big)\,\  {\rm d}x\,\,{\rm d}y\,\, {\rm d}z; 
\end{equation*}
Therefore the equilibrium conditions  are\,:
\begin{equation*} \tag*{[3.19]}
\iiint \Big[ \varrho \Big( X \delta x\,\,+\,\,Y\delta y\,\,+ Z\delta z \Big) \,\, + p \Big(\frac{d \delta x}{dx}\,\,+  \frac{d \delta y}{dy}\,\,+\frac{d \delta z}{dz}\,\,\Big) \Big]\,\  {\rm d}x\,\,  {\rm d}y\,\,  {\rm d}z = 0
\end{equation*}
\indent\indent
Since this equation is identical with $(1)$, the equilibrium  equations for  elastic and  liquid  flows are formally the same.\hspace{0.4cm} Also here, 
in accordance with equation $(2)$, we have\,:
\begin{equation} \tag*{[3.20]}
\frac{dp}{dx} = \varrho X, \,\,\frac{dp}{dy} = \varrho Y,\,\,\frac{dp}{dz} = \varrho Z 
\end{equation}
For elastic flows, $\rho$ is a 
given function of $p$, 
say, $\varrho = \varphi (p)$.\hspace{0.4cm} Let us put
\begin{equation*} \tag*{[3.21]}
f(p) = \int \frac{{\rm d}p}{\varphi(p)}
\end{equation*}
from which it follows obviously that
\begin{math} \displaystyle \frac{1}{\varrho} \frac{dp}{dx} = \frac{1}{\varphi(p)} \frac{dp}{dx} = \frac{df(p)}{dp} \frac{dp}{dx} = \frac{df(p)}{dx}, \end{math}
therefore, the three equations for the equilibrium condition become\,:
\begin{equation*} \tag*{[3.22]}
\frac{df(p)}{dx} = X, \,\,\, \frac{df(p)}{dy} = Y, \,\,\, \frac{df(p)}{dx} = Z\,,
 \end{equation*}
so that, also for elastic fluids in equilibrium,  $X,Y,Z$  have to be 
partial differential quotients
of  the same function with respect to $x,y,z$.\hspace{0.4cm}As $\varphi(p)$, and consequently also $f(p)$ is known,  $p$ 
can always be expressed through  $X,Y,Z$.\\[1cm]
\centerline{\fett{$\S.\,4.$}}\\[0.3cm]
\indent\indent
As follows from the considerations of
  $\S.\,3$, the principle of virtual velocities and lost forces
for the motion of liquid and elastic fluids, implies\,:
\leqnomode
\begin{smallequation} \notag
  \,\,\,\,0=\iiint \biggl\{ \varrho \biggl[ \Big( X - \frac{d^2 x}{d t^2} \Big) \delta x + \Big( Y - \frac{d^2 y}{d t^2}\Big) \delta y + \Big(Z - \frac{d^2 z}{d t^2}\Big) \delta z \biggr]  + p \biggl[ \frac{d \delta x}{dx} + 
\frac{d \delta y}{dy} +  \frac{d \delta z}{dz} \biggr] \biggr\} {\rm d}x {\rm d}y {\rm d}z  
\end{smallequation}
\reqnomode 
\vskip-1.15cm\begin{align}\tag*{[4.1]}
\end{align}
\vskip-1.07cm \mbox{\!\!\!\!\!\!(1)}\vskip0.8cm
\noindent from which   follows  firstly
equation $(2)$ of $\S.\,3$\,:
\leqnomode
\begin{align} \tag{2}
 \int p \big(\delta x \cos \alpha  + \delta y \cos \beta + \delta z \cos \gamma \big) {\rm d} \omega = 0,
\end{align}
\reqnomode 
\vskip-1.5cm\begin{align}\tag*{[4.2]}
\end{align}
which concerns only the external surface;
secondly, we have for the fundamental  equations of liquid or elastic fluids, 
\leqnomode
\begin{align}\tag{3}
  \varrho \Big ( \frac{d^2 x}{d t^2} - X \Big) + \frac{dp}{dx} =0,\,\, \varrho \Big ( \frac{d^2 y}{d t^2} - Y \Big) + \frac{dp}{dy} =0, \,\,\varrho \Big ( \frac{d^2 z}{d t^2} - Z \Big) + \frac{dp}{dz} =0  \,,
\end{align}
\reqnomode 
\vskip-1.5cm\begin{align}\tag*{[4.3]}
\end{align}
\noindent 
where $p$ indicates the pressure in each point.
\\\\
\indent \indent According to the first  \mbox{E\ws u\ws l\ws e\ws r\ws ian} method, the  components $u,\,v,\,\omega$  are  considered as function of time and space $x,\,y,\,z$. Therefore,
\begin{mynewequation}\tag*{[4.4]}
 \begin{aligned}
\frac{d^2 x}{dt^2} &=\frac{du}{dt} + \frac{du}{dx} u + \frac{du}{dy} v + \frac{du}{dz} w\,\,\,\\
\frac{d^2 y}{dt^2} &=\frac{dv}{dt} + \frac{dv}{dx} u + \frac{dv}{dy} v + \frac{dv}{dz} w\,\,\,\,\\
\frac{d^2 z}{dt^2} &=\frac{dw}{dt} + \frac{dw}{dx} u + \frac{dw}{dy} v + \frac{dw}{dz} w
\end{aligned} 
\end{mynewequation}
so that the fundamental equations in this form are the following\,:
\begin{mynewequation}\tag*{(4), [4.5]}
\left.
\begin{aligned}
\frac{du}{dt}\, +\, \frac{du}{dx} u + \frac{du}{dy} v + \frac{du}{dz} w - X + \frac{1}{\varrho}\frac{dp}{dx}=0\\\\
\frac{dv}{dt}\, +\,\, \frac{dv}{dx} u + \frac{dv}{dy} v + \frac{dv}{dz} w - Y + \frac{1}{\varrho}\frac{dp}{dy}=0\\\\
\frac{dw}{dt}\, + \frac{dw}{dx} u + \frac{dw}{dy} v + \frac{dw}{dz} w - Z + \frac{1}{\varrho}\frac{dp}{dz}=0
\end{aligned} \qquad \right\}
\end{mynewequation}
which may be also written in such a way that each equation is obtained from the other by cyclic permutation, namely\,:
\begin{mynewequation}\tag*{(5), [4.6]}
\left.
\begin{aligned}
\frac{du}{dt}\, +\, \frac{du}{dx} u + \frac{du}{dy} v + \frac{du}{dz} w - X + \frac{1}{\varrho}\frac{dp}{dx}=0\\\\
\frac{dv}{dt}\, +\,\, \frac{dv}{dy} v + \frac{dv}{dz}  w + \frac{dv}{dx} u - Y + \frac{1}{\varrho}\frac{dp}{dy}=0\\\\
\frac{dw}{dt}\, + \frac{dw}{dz} w + \frac{dw}{dx} u + \frac{dw}{dy} v - Z + \frac{1}{\varrho}\frac{dp}{dz}=0
\end{aligned} \qquad \right\}
\end{mynewequation} 
In addition to these formulae, there is also the density equation~(5) of the $\S.\,2$:
\be \tag*{[4.7]}
\frac{d\varrho}{dt} + \frac{d(\varrho u )}{dx} + \frac{d(\varrho v )}{dy} + \frac{d(\varrho w )}{dz} =0
\ee
or, 
in particular, for liquid flows 
\be \tag*{[4.8]}
 \frac{du}{dx} + \frac{dv}{dy} + \frac{dw}{dz}=0
\ee
We see that  these four equations are sufficient to determine the four unknowns $u,v,w$ and $p$ as functions of $x,y,z$ and $t$;\,  $\varrho$ is either a known function of $p$ or a constant.\\\\
\indent\indent
By the second 
E\ws u\ws l\ws erian  representation [Lagrangian representation], equations~$(3)$ 
are respectively multiplied by 
$$
\frac{dx}{da},\,\,\frac{dy}{da},\,\,\frac{dz}{da}$$
and summed, then, similarly, by
$$
\frac{dx}{db},\,\,\frac{dy}{db},\,\,\frac{dz}{db}$$
and
$$
\frac{dx}{dc},\,\,\frac{dy}{dc},\,\,\frac{dz}{dc}$$
\noindent
Thus, one obtains for  the three fundamental equations:
\begin{equation*}\tag*{(6), [4.9]}
\left.
\begin{aligned}
 \Big( \frac{d^2 x}{d t^2} - X \Big)\frac{dx}{da}  +  \Big( \frac{d^2 y}{d t^2} - Y \Big) \frac{dy}{da}+\Big( \frac{d^2 z}{d t^2} - Z \Big) \frac{dz}{da} + \frac{1}{\varrho} \frac{dp}{da} =0\\\
\Big( \frac{d^2 x}{d t^2} - X \Big)\frac{dx}{db}  +  \Big( \frac{d^2 y}{d t^2} - Y \Big) \frac{dy}{db}+\Big( \frac{d^2 z}{d t^2} - Z \Big) \frac{dz}{db} + \frac{1}{\varrho} \frac{dp}{db} =0\\\
\Big( \frac{d^2 x}{d t^2} - X \Big)\frac{dx}{dc}  +  \Big( \frac{d^2 y}{d t^2} - Y \Big) \frac{dy}{dc}+\Big( \frac{d^2 z}{d t^2} - Z \Big) \frac{dz}{dc} + \frac{1}{\varrho} \frac{dp}{dc} =0
\end{aligned} \qquad \right\}
\end{equation*}
and, in addition, there is the density equation (1) of $\S.\,2$
\begin{myequation} \tag*{[4.10]}
 \begin{vmatrix}
\frac{dx}{da}&\frac{dx}{db}&\frac{dx}{dc} \\\\
\frac{dy}{da}&\frac{dy}{db}&\frac{dy}{dc} \\\\
\frac{dz}{da}&\frac{dz}{db}&\frac{dz}{dc} 
\end{vmatrix}
=\frac{\varrho_0}{\varrho}
\end{myequation}
From these four equations $x,y,z$ and $p$ are found as functions of the initial 
location $a,b,c$ and time $t$.\\

\indent\indent Evidently, the  solutions of these partial differential equations must  contain  arbitrary functions, which have to be determined from  initial conditions and are in accordance with  the nature of the walls and the  flow boundaries.\\

\indent\indent These last equations, which  are usually called after \mbox{L\ws a\ws g\ws r\ws a\ws n\ws g\ws e,} significantly simplify their form by 
setting\,:
\be \tag*{[4.11]}
 X= \frac{dV}{dx}, \,\,Y=\frac{dV}{dy}, \,\, Z=\frac{dV}{dz}
\ee 
so they become\,:
\begin{equation*}\tag*{(7), [4.12]}
\left.
\begin{aligned}
 \frac{d^2 x}{d t^2}\frac{dx}{da}  +  \frac{d^2 y}{d t^2} \frac{dy}{da}+\frac{d^2 z}{d t^2} \frac{dz}{da} - \frac{dV}{da}  +  \frac{1}{\varrho}\frac{dp}{da}=0\\\\
 \frac{d^2 x}{d t^2}\frac{dx}{db}  +  \frac{d^2 y}{d t^2} \frac{dy}{db}+\frac{d^2 z}{d t^2} \frac{dz}{db} - \frac{dV}{db}  +  \frac{1}{\varrho}\frac{dp}{db}=0\\\\ 
\frac{d^2 x}{d t^2}\frac{dx}{dc}  +  \frac{d^2 y}{d t^2} \frac{dy}{dc}+\frac{d^2 z}{d t^2} \frac{dz}{dc} - \frac{dV}{dc}  +  \frac{1}{\varrho}\frac{dp}{dc}=0
\end{aligned} \qquad \right\}
\end{equation*}
We will limit ourselves to this 
assumption about  $X,Y,Z$ which, apart from the  boundary conditions, 
\noindent
coincides with the one that necessarily has to hold in the equilibrium state of the fluid.
\\
\centerline{\fett{$\S.\,5.$}}\\[0.3cm]
\indent\indent
Partially integrating the last three terms of equation (1) of $\S.\,4$, 
ignoring the boundary contributions to the double integrals and using as already stated
\be \tag*{[5.1]}
  X= \frac{dV}{dx}, \,\,Y=\frac{dV}{dy}, \,\, Z=\frac{dV}{dz} \,,
\ee 
one obtains the following equation\,: 
\be \tag*{[5.2]}
  0 = \iiint \varrho\,{\rm d}x\,{\rm d}y\, {\rm d}z\ \bigg\{\Big(\frac{d^2 x}{d t^2} - \frac{dV}{dx} + \frac{1}{\varrho}\frac{dp}{dx}\Big) \delta x + \Big(\frac{d^2 y}{d t^2} - \frac{dV}{dy} + \frac{1}{\varrho}\frac{dp}{dy}\Big) \delta y + \Big(\frac{d^2 z}{d t^2} - \frac{dV}{dz} + \frac{1}{\varrho}\frac{dp}{dz}\Big) \delta z  \bigg\}
\ee
If one puts, as  in $\S.\,3$,  $\varrho= \varphi (p )$ and:
\be  \tag*{[5.3]}
 f(p) =\int\frac{ {\rm d}p}{\varphi(p)}
\ee
one has, 
\begin{smallequation} \tag*{[5.4]}
 \!\!\!\!0 = \iiint \varrho_0 \, {\rm d}a\,  {\rm d}b\,  {\rm d}c \bigg\{\Big(\frac{d^2 x}{d t^2} - \frac{dV}{dx} + \frac{df(p)}{dx}\Big) \delta x + \Big(\frac{d^2 y}{d t^2} - \frac{dV}{dy} + \frac{df(p)}{dy}\Big) \delta y + \Big(\frac{d^2 z}{d t^2} - \frac{dV}{dz} + \frac{df(p)}{dz}\Big) \delta z  \bigg\} \,,
\end{smallequation}
provided that the transformation of the integral  is made according to  \S.\,2.
Now, if one sets 
\be \tag*{[5.5]} \label{pointerOmega}
 V - f(p) = \Omega,
\ee 
 one can write\,:
\be \tag*{[5.6]}
 0 = \iiint \varrho_0 \, {\rm d}a\,  {\rm d}b\,  {\rm d}c \bigg[ \frac{d^2 x}{d t^2} \delta x  +  \frac{d^2 y}{d t^2} \delta y + \frac{d^2 z}{d t^2} \delta z -\delta \Omega  \bigg]
\ee
If one now integrates under the triple integral 
 with respect to 
the variable $t$ which is independent of $a,b,c$, one obtains
\begin{equation} \tag*{[5.7]}
  \begin{aligned}
\int \frac{d^2 x}{dt^2} \delta x\, {\rm d}t = \Big[\frac{dx}{dt} \delta x\Big] - \int \frac{dx}{dt}
\frac{d \delta x}{dt} {\rm d} t\\\\
\int\frac{d^2 y}{dt^2} \delta y\, {\rm d}t= \Big[\frac{dy}{dt} \delta y\Big] - \int\frac{dy}{dt}\frac{d \delta y}{dt} {\rm d }t\\\\
\int\frac{d^2 z}{dt^2} \delta z\, {\rm d}t = \Big[\frac{dz}{dt} \delta z\Big] - \int\frac{dz}{dt}\frac{d \delta z}{dt} {\rm d} t 
\end{aligned} 
\end{equation}
\noindent Dropping the triple integrals, 
one has the equation\,:
\be \tag*{[5.8]}
 0 = \iiint \varrho_0 {\rm d}a\, {\rm d}b\,  {\rm d}c \int {\rm d}t \bigg\{\frac{dx}{dt} \delta \frac{dx}{dt} + \frac{dy}{dt} \delta \frac{dy}{dt}  +\frac{dz}{dt} \delta \frac{dz}{dt} + \delta \Omega\bigg\}
\ee 
which  coincides with the following\,:\hypertarget{eq1}{}
\begin{equation} \tag*{(1), [5.9]} 
0 = \delta \iiint \varrho_0\, {\rm d}a\, {\rm d}b\, {\rm d}c \int {\rm d}t \Big [\Big(\frac{ds}{dt}\Big)^2 + 2 \Omega \Big] \,.
\end{equation}
The first three hydrodynamical fundamental equations are satisfied. 
Hence
\be \notag
 \iiint \varrho_0\, {\rm d}a\, {\rm d}b\, {\rm d}c \int {\rm d}t \Big [\Big(\frac{ds}{dt}\Big)^2 + 2 \Omega \Big]
\ee
disappears; conversely, if the first variation of this integral  with respect to  $x,y,z$ is set to  zero, 
then one obtains these first three equations. \\\\

\indent\indent
This theorem, which can be easily considered, in view of  the meaning of $\displaystyle \frac{ds}{dt}$  and $\Omega$,
as a mechanical principle, has a certain analogy with the principle of least action.
For us it possesses an analytical importance, since it gives an 
extremely simple tool for transforming the hydrodynamical equations.
\hspace{0.4cm}
Indeed, in order to introduce in these equations new coordinates, instead of  $x,y,z$,  it is only  necessary to write the arc element 
\be \tag*{[5.10]}
{\rm d}s^2= {\rm d}x^2 + {\rm d}y^2 + {\rm d}z^2
\ee
in terms of the new coordinates and, then, 
apply the simple operation of variation, using the integral in the new coordinates\,:
\be \notag
  \iiint \varrho_0\,{\rm d} a\,  {\rm d} b\, {\rm d}c \int {\rm d} t \Big[\Big(\frac{ds}{dt}\Big)^2 + 2 \Omega \Big]
\ee
By setting the coefficients of the three variations to zero, we thereby obtain three equations in a similar form as equations (3) in $\S.\,4$\,; in order to write them into the first or second E\,u\,l\,e\,r\,ian form, one has to apply analogous procedures as was done in $\S\,4$.\\
\indent\indent 
If the new coordinates\,:
\be \tag*{[5.11]}
 \varrho_1=f_1(x,\,\,y,\,\,z), \,\,\,\,\varrho_2=f_2(x,\,\,y,\,\,z),\,\,\,\,\varrho_3=f_3(x,\,\,y,\,\,z)
\ee
are used instead of $x,\,y,\,z$,  one obtains:
\begin{mynewequation} \tag*{[5.12]}
 \begin{aligned}
 {\rm d}x = \frac{dx}{d \varrho_1}{\rm d} \varrho_1 + \frac{dx}{d\varrho_2}{\rm d} \varrho_2  + \frac{dx}{d \varrho_3} {\rm d} \varrho_3\\
{\rm d}y = \frac{dy}{d \varrho_1}{\rm d} \varrho_1 + \frac{dy}{d\varrho_2}{\rm d} \varrho_2  + \frac{dy}{d \varrho_3} {\rm d} \varrho_3\\
{\rm d}z = \frac{dz}{d \varrho_1}{\rm d} \varrho_1 + \frac{dz}{d\varrho_2}{\rm d} \varrho_2  + \frac{dz}{d \varrho_3} {\rm d} \varrho_3\\
\end{aligned} 
\end{mynewequation}
therefore
\be \tag*{[5.13]}
 {\rm d}s^2 = {\rm d}x^2+{\rm d}y^2+{\rm d}z^2=N_1 {\rm d} \varrho_1^2 + N_2 {\rm d} \varrho_2^2 +N_3 {\rm d} \varrho_3^2+ 2 n_3  {\rm d} \varrho_1  {\rm d}\varrho_2 + 2n_1  {\rm d} \varrho_2 {\rm d} \varrho_3 + 2n_2 {\rm d} \varrho_3 {\rm d}\varrho_1 
\ee
where
\begin{align} \tag*{[5.14]}
\begin{aligned}
N_1=  \Big(\frac{dx}{d \varrho_1}\Big)^2 + \Big(\frac{dy}{d \varrho_1}\Big)^2 + \Big(\frac{dz}{d \varrho_1}\Big)^2\\
N_2= \Big(\frac{dx}{d \varrho_2}\Big) ^2 + \Big(\frac{dy}{d \varrho_2}\Big)^2 + \Big(\frac{dz}{d \varrho_2}\Big)^2\\
N_3= \Big(\frac{dx}{d \varrho_3}\Big) ^2 + \Big(\frac{dy}{d \varrho_3}\Big)^2 + \Big(\frac{dz}{d \varrho_3}\Big)^2
\end{aligned}  \\[0.3cm] \tag*{[5.15]}
\begin{aligned}
n_3= \frac{dx}{d \varrho_1} \frac{dx}{d \varrho_2} +  \frac{dy}{d \varrho_1} \frac{dy}{d \varrho_2} +  \frac{dz}{d \varrho_1} \frac{dz}{d \varrho_2}\\\\
n_1= \frac{dx}{d \varrho_2} \frac{dx}{d \varrho_3} +  \frac{dy}{d \varrho_2} \frac{dy}{d \varrho_3} +  \frac{dz}{d \varrho_2} \frac{dz}{d \varrho_3}\\\\
n_2= \frac{dx}{d \varrho_3} \frac{dx}{d \varrho_1} +  \frac{dy}{d \varrho_3} \frac{dy}{d \varrho_1} +  \frac{dz}{d \varrho_3} \frac{dz}{d \varrho_1}
\end{aligned} 
\end{align}
\indent\indent
Once $N_1, N_2, N_3,n_1,n_2,n_3$ are expressed  in terms of the new variables $\varrho_1, \varrho_2, \varrho_3$, one has to  vary the integral
\begin{smallequation} \notag
\iint\!\!\varrho_0  \, {\rm d}a\,  {\rm d}b\, {\rm d}c \int\!\! {\rm d}t \Big[N_1 \Big(\frac{d \varrho_1}{dt} \Big)^2\! + N_2 \Big(\frac{d \varrho_2}{dt} \Big)^2 \!+N_3 \Big(\frac{d \varrho_3}{dt} \Big)^2\! + 2 n_3 \frac{d\varrho_1}{dt} \frac{d \varrho_2}{dt} + 2 n_1 \frac{d\varrho_2}{dt} \frac{d \varrho_3}{dt}  + 2 n_2 \frac{d\varrho_3}{dt} \frac{d \varrho_1}{dt} + 2 \Omega \Big]
\end{smallequation}
with respect to these new variables. 
Then, after integration by parts with respect to $t$,  one removes  from the quadruple  integral the quantities in which appear the time derivatives of the variations $\delta \varrho_1, \delta \varrho_2, \delta \varrho_3$.
Then one has to set the coefficients of $\delta \varrho_1, \delta \varrho_2, \delta \varrho_3$ equal to zero.\hspace{0.4cm} 
After that one obtains the first three hydrodynamical fundamental equations, 
which are completed by a fourth one, the density equation.\hspace{0.6cm} In  order to express also the density equation in the new coordinates, we notice that 
\noindent the volume element ${\rm d}x\,{\rm d}y\,{\rm d}z$  may be expressed in such coordinates very easily, namely\,:
\begin{equation}\tag*{[5.16]}
\begin{aligned}
 {\rm d}x\,{\rm d}y\,{\rm d}z=  {\rm d}\varrho_1\, {\rm d}\varrho_2\, {\rm d}\varrho_3 \begin{vmatrix}
\frac{dx}{d\varrho_1}&\frac{dx}{d\varrho_2}&\frac{dx}{d\varrho_3} \\\\
\frac{dy}{d\varrho_1}&\frac{dy}{d\varrho_2}&\frac{dy}{d\varrho_3} \\\\
\frac{dz}{d\varrho_3}&\frac{dz}{d\varrho_3}&\frac{dz}{d\varrho_3} 
\end{vmatrix}
\end{aligned}
\end{equation}
herefrom follows
\begin{equation} \notag
({\rm d} x\,{\rm d}y\,{\rm d}z)^2= ({\rm d}\varrho_1\,{\rm d}\varrho_2\,{\rm d}\varrho_3)^2 {\bf {  \times}} \\[-0.5cm]
\end{equation}
\begin{equation}\tag*{[5.17]}
\begin{vmatrix}
(\frac{dx}{d\varrho_1})^2  +(\frac{dy}{d\varrho_1})^2+(\frac{dz}{d\varrho_1})^2,&
\frac{dx}{d\varrho_1} \frac{dx}{d\varrho_2} + \frac{dy}{d\varrho_1} \frac{dy}{d\varrho_2}+\frac{dz}{d\varrho_1} \frac{dz}{d\varrho_2},& 
\frac{dx}{d\varrho_3} \frac{dx}{d \varrho_1} + \frac{dy}{d\varrho_3} \frac{dy}{d \varrho_1}+\frac{dz}{d\varrho_3}\frac{dz}{d\varrho_1}\\\\
\frac{dx}{d\varrho_1 }\frac{dx}{d\varrho_2} + \frac{dy}{d\varrho_1}\frac{dy}{d\varrho_2}+\frac{dz}{d\varrho_1} \frac{dz}{d\varrho_2},&
(\frac{dx}{d\varrho_2})^2  +(\frac{dy}{d\varrho_2})^2+(\frac{dz}{d\varrho_2})^2,&
\frac{dx}{d\varrho_2} \frac{dx}{d\varrho_3} + \frac{dy}{d\varrho_2} \frac{dy}{d\varrho_3}+\frac{dz}{d\varrho_2} \frac{dz}{d\varrho_3}\\\\
\frac{dx}{d\varrho_3} \frac{dx}{d\varrho_1} + \frac{dy}{d\varrho_3} \frac{dy}{d\varrho_1}+\frac{dz}{d\varrho_3} \frac{dz}{d\varrho_1},&
\frac{dx}{d\varrho_2} \frac{dx}{d\varrho_3} + \frac{dy}{d\varrho_2} \frac{dy}{d\varrho_3}+\frac{dz}{d\varrho_2} \frac{dz}{d\varrho_3},&
(\frac{dx}{d\varrho_3})^2  +(\frac{dy}{d\varrho_3})^2+(\frac{dz}{d\varrho_3})^2&
\end{vmatrix}
\end{equation}
or, with the notation as defined above\,:
\begin{mynewequation} \tag*{[5.18]}
({\rm d}x\,{\rm d}y\,{\rm d}z)^2= ({\rm d}\varrho_1\,{\rm d}\varrho_2\,{\rm d}\varrho_3)^2 \begin{vmatrix}
N_1&n_3&n_2\\\\
n_3&N_2&n_1\\\\
n_2&n_1&N_3
\end{vmatrix}
\end{mynewequation}
\indent\indent
 Let us denote the values of  $\varrho_1, \varrho_2,\varrho_3, N_1, N_2, N_3, n_1,n_2, n_3 $ at time $t=0$ with $ \varrho_1^0, \varrho_2^0,\varrho_3^0, N_1^0, N_2^0, N_3^0, n_1^0,n_2^0, n_3^0 $,  then we get from the previous equation for $t=0$
\begin{mynewequation} \tag*{[5.19]}
({\rm d}a\,{\rm d}b\,{\rm d}c)^2= ({\rm d}\varrho_1^0\,{\rm d}\varrho_2^0\,{\rm d}\varrho_3^0)^2 \begin{vmatrix}
N_1^0&n_3^0&n_2^0\\\\
n_3^0&N_2^0&n_1^0\\\\
n_2^0&n_1^0&N_3^0
\end{vmatrix}
\end{mynewequation}
One also has to think of the ensuing value of ${\rm d}a\,{\rm d}b\,{\rm d}c$ as substituted into the integral to  be varied. But now  the density equation is,  for the general case\,:
\be \tag*{[5.20]}
\frac{ {\rm d}x\,{\rm d}y\,{\rm d}z}{{\rm d}a\,{\rm d}b\,{\rm d}c} = \frac{\varrho_0}{\varrho} \,.
\ee
Then, dividing $({\rm d}x\,\,{\rm d}y\,\,{\rm d}z)^2$ by $({\rm d}a\,\,{\rm d}b\,\,{\rm d}c)^2$, one obtains the density equation in the new variables\,:
\begin{mynewequation} \tag*{[5.21]}
\Big(\frac{{\rm d} \varrho_1}{{\rm d}\varrho_1^0} \frac{{\rm d} \varrho_2}{{\rm d}\varrho_2^0} \frac{{\rm d} \varrho_3}{{\rm d} \varrho_3^0} \Big)^2\,\,\,\Big(\frac{\varrho}{\varrho_0}\Big)^2 = \begin{vmatrix} N_1^0& n_3^0&n_2^0\\
n_3^0&N_2^0&n_1^0\\
n_2^0&n_1^0&N_3^0 
\end{vmatrix} 
:
\begin{vmatrix} N_1& n_3&n_2\\
n_3&N_2&n_1\\
n_2&n_1&N_3 
\end{vmatrix}
\end{mynewequation}
or, as it is well-known
\begin{mynewequation} \tag*{[5.22]}
  \text{\small ${\rm d}\varrho_1\,{\rm d}\varrho_2\,{\rm d}\varrho_3 = {\rm d}\varrho_1^0\,{\rm d}\varrho_2^0\,{\rm d}\varrho_3^0$}
\begin{vmatrix}
\frac{d\varrho_1}{d\varrho_1^0} &\frac{d\varrho_1}{d\varrho_2^0} &\frac{d\varrho_1}{d\varrho_3^0} \\\\
\frac{d\varrho_2}{d\varrho_1^0} &\frac{d\varrho_2}{d\varrho_2^0} &\frac{d\varrho_2}{d\varrho_3^0} \\\\
\frac{d\varrho_3}{d\varrho_1^0} &\frac{d\varrho_3}{d\varrho_2^0} &\frac{d\varrho_3}{d\varrho_3^0} 
\end{vmatrix} \,,
\end{mynewequation}
one finally obtains\,:
\begin{mynewequation}\tag*{[5.23]}
\begin{vmatrix}
\frac{d\varrho_1}{d\varrho_1^0} &\frac{d\varrho_1}{d\varrho_2^0} &\frac{d\varrho_1}{d\varrho_3^0} \\\\
\frac{d\varrho_2}{d\varrho_1^0} &\frac{d\varrho_2}{d\varrho_2^0} &\frac{d\varrho_2}{d\varrho_3^0} \\\\
\frac{d\varrho_3}{d\varrho_1^0} &\frac{d\varrho_3}{d\varrho_2^0} &\frac{d\varrho_3}{d\varrho_3^0} 
\end{vmatrix} \cdot
\text{\Large $\frac{\varrho}{\varrho_0}$} ={\boldsymbol {\surd }}\begin{vmatrix} N_1^0& n_3^0&n_2^0\\
n_3^0&N_2^0&n_1^0\\
n_2^0&n_1^0&N_3^0 
\end{vmatrix} 
:
\begin{vmatrix} N_1& n_3&n_2\\
n_3&N_2^0&n_1\\
n_2&n_1&N_3 
\end{vmatrix}
\end{mynewequation}
\indent\indent
All quantities appearing in this transformed  density equation are already known through the transformation of the arc element.\hspace{0.4cm}Here, the problem  of the transformation of the four hydrodynamical equations in an arbitrary coordinate system is reduced to the problem of the transformation of the arc element.\\
\indent\indent
It is obvious that the applicability of this procedure does not depend on the number of  variables  and  that, in the same way, by variation of the integral with respect to $x_1,x_2,...,x_n$: 
\begin{equation*} \notag
\iiint \varrho_0\, {\rm d}a_1 {\rm d}a_2\,...\,{\rm d}a_n \int {\rm d}t \biggl\{ \Big(\frac{ds}{dt}\Big)^2 + 2 \Omega\biggr\}
\end{equation*}
one obtains
\begin{mynewequation}\tag*{[5.24]}
\left.
\begin{aligned}
\frac{d^2 x_1}{d t^2} \frac{dx_1}{da_1}  +  \frac{d^2 x_2}{d t^2} \frac{dx_2}{da_1}+...+  \frac{d^2 x_n}{d t^2} \frac{dx_n}{da_1}- \frac{dV}{da_1}+ \frac{1}{\varrho}\, \frac{dp}{da_1} =0\\
\,\,.\,\,\,\,.\,\,\,\,.\,\,\,\,.\,\,\,\,.\,\,\,\,.\,\,\,\,.\,\,\,\,.\,\,\,\,.\,\qquad \qquad \qquad \qquad \qquad  \\
\frac{d^2 x_1}{d t^2} \frac{dx_1}{da_n}  +  \frac{d^2 x_2}{d t^2} \frac{dx_2}{da_n}+...+  \frac{d^2 x_n}{d t^2} \frac{dx_n}{da_n}- \frac{dV}{da_n}+ \frac{1}{\varrho} \frac{dp}{da_n} =0\\
\end{aligned} \qquad \right\}
\end{mynewequation}
\noindent
where $a_1, a_2, .\,. \,a_n$ are the values of $x_1, x_2,.\,. \,x_n$  at time $t=0$.\hspace{0.4cm}Then, the transformation of these equations happens exactly in the same way.\hspace{0.4cm}For the sake of brevity, we here limit ourselves to three variables.\\
\indent\indent
This transformation happens to be very simple when the new variables form an orthogonal system,  a case, which  apart from this simplification, is also  very interesting since all  frequently used coordinate systems are included therein.\\
\indent\indent
The points where  $\varrho_1$ takes a given value will in general
form a surface, whose equation with respect to the axes of $x,\,y,\,z$  is
\be \tag*{[5.25]}
 \varrho_1 = f_1 (x,\, y,\,z)
\ee
The cosines of the angles formed by the normal to the point  $x, \,y, \,z$ of this surface and the coordinates axes, are\,:
\be \tag*{[5.26]}
\frac{1}{\Delta_1} \,\frac{d \varrho_1}{dx}, \,\,\,\,\,\,\,\frac{1}{\Delta_1}\,\frac{d \varrho_1}{dy}, \,\,\,\,\,\,\frac{1}{\Delta_1}\,
\frac{d \varrho_1}{dz}, \,\,\,\,\,\,\,\,\,\Delta_1^2= \Big(\frac{d\varrho_1}{dx}\Big)^2 + \Big(\frac{d\varrho_1}{dy}\Big)^2 +\Big(\frac{d\varrho_1}{dz}\Big)^2 
\ee
The analogous cosines for the normal to the surface
\be \tag*{[5.27]}
 \varrho_2 = f_2 (x,\, y,\, z)
\ee
are\,:
\be \tag*{[5.28]}
 \frac{1}{\Delta_2} \,\frac{d \varrho_2}{dx}, \frac{1}{\Delta_2}\,\frac{d \varrho_2}{dy}, \frac{1}{\Delta_2}\,
\frac{d \varrho_2}{dz}, \,\,\,\,\,\,\,\,\,\Delta_2^2= \Big(\frac{d\varrho_2}{dx}\Big)^2 + \Big(\frac{d\varrho_2}{dy}\Big)^2 +\Big(\frac{d\varrho_2}{dz}\Big)^2; 
\ee
for the surface
\be \tag*{[5.29]}
 \varrho_3 = f_3 (x, y, z)
\ee
the analogous cosines will be\,:
\be \tag*{[5.30]}
 \frac{1}{\Delta_3} \,\frac{d \varrho_3}{dx}, \frac{1}{\Delta_3}\,\frac{d \varrho_3}{dy}, \frac{1}{\Delta_3}\,
\frac{d \varrho_3}{dz}, \,\,\,\,\,\,\,\,\,\Delta_3^2= \Big(\frac{d\varrho_3}{dx}\Big)^2 + \Big(\frac{d\varrho_3}{dy}\Big)^2 +\Big(\frac{d\varrho_3}{dz}\Big)^2 
\ee
\indent\indent
Suppose now that $\varrho_1, \, \varrho_2, \, \varrho_3$ form an orthogonal system; 
then, in the point of intersection the normals to the three surfaces $\varrho_1, \, \varrho_2, \, \varrho_3$ are mutually orthogonal to each other. The conditions that the axes of $x,\,y,\,z$ form an orthogonal system, when referring their positions to the normals of the surfaces $\varrho_1, \varrho_2,\,\varrho_3$, will be the following\,:
\begin{mynewequation} \tag*{[5.31]}
\left.
\begin{aligned}
 \frac{1}{\Delta_1^2} \Big(\frac{d \varrho_1}{dx}\Big)^2  + \frac{1}{\Delta_2^2} \Big(\frac{d \varrho_2}{dx}\Big)^2 + \frac{1}{\Delta_3^2} \Big(\frac{d \varrho_3}{dx}\Big)^2  =1\\\\
\frac{1}{\Delta_1^2} \Big(\frac{d \varrho_1}{dy}\Big)^2  + \frac{1}{\Delta_2^2} \Big(\frac{d \varrho_2}{dy}\Big)^2 + \frac{1}{\Delta_3^2} \Big(\frac{d \varrho_3}{dy}\Big)^2  =1\\\\
\frac{1}{\Delta_1^2} \Big(\frac{d \varrho_1}{dz}\Big)^2  + \frac{1}{\Delta_2^2} \Big(\frac{d \varrho_2}{dz}\Big)^2 + \frac{1}{\Delta_3^2} \Big(\frac{d \varrho_3}{dz}\Big)^2  =1
\end{aligned} \qquad\right\}
\end{mynewequation}
\begin{mynewequation} \tag*{[5.32]}
\left.
\begin{aligned}
 \frac{1}{\Delta_1^2} \frac{d \varrho_1}{dx} \frac{d \varrho_1}{dy}  + \frac{1}{\Delta_2^2} \frac{d \varrho_2}{dx} \frac{d \varrho_2}{dy}+ \frac{1}{\Delta_3^2} \frac{d \varrho_3}{dx} \frac{d \varrho_3}{dy} =0\\\\
\frac{1}{\Delta_1^2} \frac{d \varrho_1}{dy} \frac{d \varrho_1}{dz}  + \frac{1}{\Delta_2^2} \frac{d \varrho_2}{dy} \frac{d \varrho_2}{dz}+ \frac{1}{\Delta_3^2} \frac{d \varrho_3}{dy} \frac{d \varrho_3}{dz} =0\\\\
\frac{1}{\Delta_1^2} \frac{d \varrho_1}{dz} \frac{d \varrho_1}{dx}  + \frac{1}{\Delta_2^2} \frac{d \varrho_2}{dz} \frac{d \varrho_2}{dx}+ \frac{1}{\Delta_3^2} \frac{d \varrho_3}{dz} \frac{d \varrho_3}{dx} =0
\end{aligned} \qquad \right\}
\end{mynewequation} 
\indent \indent Now we have:
\begin{mynewequation} \tag*{[5.33]}
\begin{aligned}
 \frac{\rm {d} \varrho_1}{\Delta_1} =\frac{1}{\Delta_1} \frac{d \varrho_1}{dx} {\rm {d}x} + \frac{1}{\Delta_1} \frac{d \varrho_1}{dy} {\rm {d}y} + \frac{1}{\Delta_1} \frac{d \varrho_1}{dz} {\rm {d}z} \\\\
 \frac{\rm {d} \varrho_2}{\Delta_2} =\frac{1}{\Delta_2} \frac{d \varrho_2}{dx} {\rm {d}x} + \frac{1}{\Delta_2} \frac{d \varrho_2}{dy} {\rm {d}y} + \frac{1}{\Delta_2} \frac{d \varrho_2}{dz} {\rm {d}z} \\\\
 \frac{\rm {d} \varrho_3}{\Delta_3} =\frac{1}{\Delta_3} \frac{d \varrho_3}{dx} {\rm {d}x} + \frac{1}{\Delta_3} \frac{d \varrho_3}{dy} {\rm {d}y} + \frac{1}{\Delta_2} \frac{d \varrho_3}{dz}{\rm {d}z}
\end{aligned} 
\end{mynewequation}
\indent\indent
By squaring  and adding these equations, using  the given relations,  one obtains\,: 
\begin{mynewequation} \tag*{[5.34]}
\Big(\frac{\rm {d} \varrho_1}{\Delta_1}\Big)^2 + \Big(\frac{\rm {d} \varrho_2}{\Delta_2}\Big)^2 + \Big(\frac{\rm {d} \varrho_3}{\Delta_3}\Big)^2 \text{\small $={\rm d}x^2+{\rm d}y^2+{\rm d}z^2= {\rm d}s^2$}
\end{mynewequation} 
The comparison of this expression with the general one yields\,:
\begin{mynewequation} \tag*{[5.35]}
N_1 = \frac{1}{A_1^2},\,\,\,N_2 = \frac{1}{A_2^2},\,\,\,N_3 = \frac{1}{A_3^2},\,\,\,n_1=n_2=n_3=0
\end{mynewequation} 
If one sets\,:
\begin{mynewequation} \tag*{[5.36]}
 \begin{aligned}
N_1=\Big(\frac{dx}{d\varrho_1}\Big)^2 + \Big(\frac{dy}{d\varrho_1}\Big)^2 + \Big(\frac{dz}{d\varrho_1}\Big)^2=
\frac{1}{(\frac{d \varrho_1}{dx})^2 + (\frac{d \varrho_1}{dy})^2 + (\frac{d \varrho_1}{dz})^2 }\\
N_2=\Big(\frac{dx}{d\varrho_2}\Big)^2 + \Big(\frac{dy}{d\varrho_2}\Big)^2 + \Big(\frac{dz}{d\varrho_3}\Big)^2=
\frac{1}{(\frac{d \varrho_2}{dx})^2 + (\frac{d \varrho_2}{dy})^2 + (\frac{d \varrho_2}{dz})^2 }\\
N_3=\Big(\frac{dx}{d\varrho_3}\Big)^2 + \Big(\frac{dy}{d\varrho_3}\Big)^2 + \Big(\frac{dz}{d\varrho_3}\Big)^2=
\frac{1}{(\frac{d \varrho_3}{dx})^2 + (\frac{d \varrho_3}{dy})^2 + (\frac{d \varrho_3}{dz})^2 }
 \end{aligned} 
\end{mynewequation}
so the equation for the density becomes\,:
\begin{myequation} \tag*{(2), [5.37]}
\begin{vmatrix}
\frac{d\varrho_1}{d\varrho_1^0}\,\,\,\, &\frac{d\varrho_1}{d\varrho_2^0}\,\,\,\, &\frac{d\varrho_1}{d\varrho_3^0} \\\\
\frac{d\varrho_2}{d\varrho_1^0}\,\,\,\, &\frac{d\varrho_2}{d\varrho_2^0}\,\,\,\, &\frac{d\varrho_2}{d\varrho_3^0} \\\\
\frac{d\varrho_3}{d\varrho_1^0}\,\,\,\, &\frac{d\varrho_3}{d\varrho_2^0}\,\,\,\, &\frac{d\varrho_3}{d\varrho_3^0} 
\end{vmatrix}
\cdot  \frac{\varrho}{\varrho_0} =\sqrt {\frac{ N_1^0\, N_2^0\, N_3^0}{N_1\, N_2\, N_3}}
\end{myequation}
and the equation  [``expression'' is here meant]
to be varied is\,:
\begin{mynewequation} \notag
\int {\rm d} t \int {\rm d} T \left \{ N_1 \left(\frac{d\varrho_1}{dt} \right)^2 +  N_2 \left(\frac{d \varrho_2}{dt} \right)^2 +  N_3 \left(\frac{d \varrho_3}{dt} \right)^2 + 2 \Omega \right \}
\end{mynewequation}
where ${\rm d}T$ indicates the new element of mass $\varrho_0\,{\rm d}a\,{\rm d}b\,{\rm d}c$ written in the new coordinates.\hspace{0.4cm}The part  
dependent on $\deltavarrho_1$  of the variation of  this integral is\,: 
\begin{smallequation} \notag
\int {\rm d} t \int {\rm d} T \left \{ 2 N_1  \frac{d \varrho_1}{dt} \frac{d \deltavarrho_1}{dt} + \left( \frac{d \varrho_1}{dt} \right)^2 \frac{d N_1}{d \varrho_1} \deltavarrho_1 + \left( \frac{d \varrho_2}{dt} \right)^2 \frac{d N_2}{d \varrho_1} \deltavarrho_1 + \left( \frac{d \varrho_3}{dt} \right)^2 \frac{d N_3}{d \varrho_1} \deltavarrho_1  + 2 \frac{d \Omega}{d \varrho_1} \deltavarrho_1 \right \} \,.
\end{smallequation}
When the first member of this expression is integrated by parts in $t$, all members have the factor $\deltavarrho_1$. After it is set to zero, one obtains\,:
\begin{myequation} \tag*{[5.38]}
{\text{\footnotesize $2$}}  \frac{d \Omega}{d\varrho_1} = \text{\footnotesize $2$}  \frac{d \left(N_1  \frac{d \varrho_1}{dt}\right)} {dt} -\left(  \frac{d \varrho_1}{dt}  \right)^2 \frac{dN_1}{d \varrho_1} -\left(  \frac{d \varrho_2}{dt}  \right)^2 \frac{dN_2}{d \varrho_1} 
-\left(  \frac{d \varrho_3}{dt}  \right)^2 \frac{dN_3}{d \varrho_1} \,,
\end{myequation}
similarly, one has\,:
\begin{equation*} \tag*{(3), [5.39]}
2\frac{d \Omega}{d\varrho_2} = 2 \frac{d \left(N_2  \frac{d \varrho_2}{dt}\right)} {dt} -\left(  \frac{d \varrho_1}{dt}  \right)^2 \frac{dN_1}{d \varrho_2} -\left(  \frac{d \varrho_2}{dt}  \right)^2 \frac{dN_2}{d \varrho_2} 
-\left(  \frac{d \varrho_3}{dt}  \right)^2 \frac{dN_3}{d \varrho_2} 
\end{equation*}
\be \tag*{[5.40]}
2\frac{d \Omega}{d\varrho_3} = 2 \frac{d \left(N_3  \frac{d \varrho_3}{dt}\right)} {dt} -\left(  \frac{d \varrho_1}{dt}  \right)^2 \frac{dN_1}{d \varrho_3} -\left(  \frac{d \varrho_2}{dt}  \right)^2 \frac{dN_2}{d \varrho_3} 
-\left(  \frac{d \varrho_3}{dt}  \right)^2 \frac{dN_3}{d \varrho_3} 
\ee
These are equations which are built analogously to (3) of  \S\,4\,;  in order for  
these equations to take the second \mbox{E\ws u\ws l\ws e\ws rian} form,  
the so-called  \mbox{L\ws a\ws g\ws r\ws a\ws n\ws g\ws ian}  form, we multiply in turn the previous equations by\,:
$$
\frac{d \varrho_1}{d \varrho_1^0}, \,\,\,\,\,\, \frac{d \varrho_2}{d \varrho_1^0}, \,\,\,\,\,\, \frac{d \varrho_3}{d \varrho_1^0},
$$
and we add them; then we multiply by
$$
\frac{d \varrho_1}{d \varrho_2^0}, \,\,\,\,\,\, \frac{d \varrho_2}{d \varrho_2^0}, \,\,\,\,\,\, \frac{d \varrho_3}{d \varrho_2^0},
$$
and we add them as well; finally, we multiply by
$$
\frac{d \varrho_1}{d \varrho_3^0}, \,\,\,\,\,\, \frac{d \varrho_2}{d \varrho_3^0}, \,\,\,\,\,\, \frac{d \varrho_3}{d \varrho_3^0},
$$
and we add them too.\hspace{0.4cm}In this way, we get the following equations\,:
\begin{mynewequation} \tag* {(4), [5.41]}
\left.
\begin{aligned}
 2 \frac{d \Omega}{d {\varrho_1}^0} =  2 \frac{d\left(N_1 \frac{d\varrho_1}{dt} \right)}{dt}  \frac{d \varrho_1}{d \varrho_1^0}  + 2 \frac{d\left(N_2 \frac{d\varrho_2}{dt} \right)}{dt}  \frac{d \varrho_2}{d \varrho_1^0} + 2 \frac{d\left(N_3 \frac{d\varrho_3}{dt} \right)}{dt}  \frac{d \varrho_3}{d \varrho_1^0} \\
\,\,\,\,\,\,\,\,\,\,\,\,\,\,\,\,\,\,  -\left(  \frac{d \varrho_1}{dt}  \right)^2 \frac{dN_1}{d \varrho_1^0} -\left(  \frac{d \varrho_2}{dt}  \right)^2 \frac{dN_2}{d \varrho_1^0} 
-\left(  \frac{d \varrho_3}{dt}  \right)^2 \frac{dN_3}{d \varrho_1^0} \\\\
2 \frac{d \Omega}{d {\varrho_2}^0} =  2 \frac{d\left(N_1 \frac{d\varrho_1}{dt} \right)}{dt}  \frac{d \varrho_1}{d \varrho_2^0}  + 2 \frac{d\left(N_2 \frac{d\varrho_2}{dt} \right)}{dt}  \frac{d \varrho_2}{d \varrho_2^0} + 2 \frac{d\left(N_3 \frac{d\varrho_3}{dt} \right)}{dt}  \frac{d \varrho_3}{d \varrho_2^0} \\
\,\,\,\,\,\,\,\,\,\,\,\,\,\,\,\,\,\,  -\left(  \frac{d \varrho_1}{dt}  \right)^2 \frac{dN_1}{d \varrho_2^0} -\left(  \frac{d \varrho_2}{dt}  \right)^2 \frac{dN_2}{d \varrho_2^0} 
-\left(  \frac{d \varrho_3}{dt}  \right)^2 \frac{dN_3}{d \varrho_2^0} \\\\
2 \frac{d \Omega}{d {\varrho_3}^0} =  2 \frac{d\left(N_1 \frac{d\varrho_1}{dt} \right)}{dt}  \frac{d \varrho_1}{d \varrho_3^0}  + 2 \frac{d\left(N_2 \frac{d\varrho_2}{dt} \right)}{dt}  \frac{d \varrho_2}{d \varrho_3^0} + 2 \frac{d\left(N_3 \frac{d\varrho_3}{dt} \right)}{dt}  \frac{d \varrho_3}{d \varrho_3^0} \\
\,\,\,\,\,\,\,\,\,\,\,\,\,\,\,\,\,\,  -\left(  \frac{d \varrho_1}{dt}  \right)^2 \frac{dN_1}{d \varrho_3^0} -\left(  \frac{d \varrho_2}{dt}  \right)^2 \frac{dN_2}{d \varrho_3^0} 
-\left(  \frac{d \varrho_3}{dt}  \right)^2 \frac{dN_3}{d \varrho_3^0} 
\end{aligned} \qquad \right\} 
\end{mynewequation}
\indent\indent
A very elegant example of an orthogonal system are the elliptical coordinates $\varrho_1, \varrho_2,\varrho_3$ which can be defined as the roots of the equation with respect to  $\varepsilon$\,: 
\be \tag*{[5.42]}
\frac{x^2}{\alpha ^2 - \varepsilon^2} + \frac{y^2}{\beta^2 - \varepsilon^2} + \frac{z^2}{\gamma^2 - \varepsilon^2} =1
\ee
and so taken  that one has\,:
\be \tag*{[5.43]}
 \alpha  > \varrho_1 > \beta > \varrho_2>\gamma >\varrho_3 >0
\ee
\noindent
From the identity obtained after partial fraction decomposition\,: 
\begin{align} 
\frac {\left( \varepsilon^2 - \varrho_1^2 \right)  \left( \varepsilon^2 - \varrho_2^2 \right)  \left( \varepsilon^2 - \varrho_3^2 \right)}
{\left( \varepsilon^2 - \alpha^2 \right) \left( \varepsilon^2 - \beta^2 \right)  \left( \varepsilon^2 - \gamma^2\right)}  = 1 - \frac{\left(\alpha^2 - \varrho_1^2 \right) \left(\alpha^2 - \varrho_2^2 \right) \left(\alpha^2 - \varrho_3^2 \right)  }{ \left (\alpha^2 - \beta^2 \right)\left(\alpha^2 - \gamma^2 \right)} \cdot \frac{1}{\alpha^2 - \varepsilon^2}  \hspace{2cm} \nonumber \\
 - \frac{\left(\beta^2 - \varrho_1^2 \right) \left(\beta^2 - \varrho_2^2 \right) \left(\beta^2 - \varrho_3^2 \right)  }{ \left (\beta^2 - \alpha^2 \right)\left(\beta^2 - \gamma^2 \right)} \frac{1}{\beta^2 - \varepsilon^2} - \frac{\left(\gamma^2 - \varrho_1^2 \right) \left(\gamma^2 - \varrho_2^2 \right) \left(\gamma^2 - \varrho_3^2 \right)  }{ \left (\gamma^2 - \beta^2 \right)\left(\gamma^2 - \alpha^2 \right)} \frac{1}{\gamma^2 - \varepsilon^2} \,. \tag*{[5.44]}
\end{align}
When $\varepsilon$ is put equal to  one of the roots 
$\varrho_1,\varrho_2,\varrho_3$ of the equation\,:
\be \tag*{[5.45]}
\frac{x^2}{\alpha^2 - \varepsilon^2} + \frac{y^2}{\beta^2 - \varepsilon^2} + \frac{z^2}{\gamma^2 -\varepsilon^2}=1
\ee
 one obtains\,:
\begin{align} 
 \frac{\left(\alpha^2 - \varrho_1^2 \right) \left(\alpha^2 - \varrho_2^2 \right) \left(\alpha^2 - \varrho_3^2 \right)  }{ \left (\alpha^2 - \beta^2 \right)\left(\alpha^2 - \gamma^2 \right)}  \frac{1}{\alpha^2 - \varepsilon^2} + \frac{\left(\beta^2 - \varrho_1^2 \right) \left(\beta^2 - \varrho_2^2 \right) \left(\beta^2 - \varrho_3^2 \right)  }{ \left (\beta^2 - \alpha^2 \right)\left(\beta^2 - \gamma^2 \right)} \frac{1}{\beta^2 - \varepsilon^2} \nonumber \\ + \frac{\left(\gamma^2 - \varrho_1^2 \right) \left(\gamma^2 - \varrho_2^2 \right) \left(\gamma^2 - \varrho_3^2 \right)  }{ \left (\gamma^2 - \beta^2 \right)\left(\gamma^2 - \alpha^2 \right)} \frac{1}{\gamma^2 - \varepsilon^2} =1  \tag*{[5.46]}
\end{align}
which compared with the previous equations gives the relations\,:
\begin{align}  \tag*{[5.47]}
 \begin{aligned}
x^2 = \frac {\left( \alpha^2 - \varrho_1^2 \right)  \left( \alpha^2 - \varrho_2^2 \right)  \left( \alpha^2 - \varrho_3^2 \right)}
{\left( \alpha^2 - \beta^2 \right) \left( \alpha^2 - \gamma^2 \right)} \\  \\
y^2 = \frac {\left( \beta^2 - \varrho_1^2 \right)  \left( \beta^2 - \varrho_2^2 \right)  \left( \beta^2 - \varrho_3^2 \right)}
{\left( \beta^2 - \alpha^2 \right) \left( \beta^2 - \gamma^2 \right)} \\ \\
z^2 = \frac {\left( \gamma^2 - \varrho_1^2 \right)  \left( \gamma^2 - \varrho_2^2 \right)  \left( \gamma^2 - \varrho_3^2 \right)}
{\left( \gamma^2 - \alpha^2 \right) \left( \gamma^2 - \beta^2 \right)} 
\end{aligned} 
\end{align}
so that, if $\varepsilon$ is an arbitrary variable, one has\,:
\be \tag*{[5.48]}
\frac {\left( \varepsilon^2 - \varrho_1^2 \right)  \left( \varepsilon^2 - \varrho_2^2 \right)  \left( \varepsilon^2 - \varrho_3^2 \right)}
{\left( \varepsilon^2 - \alpha^2 \right) \left( \varepsilon^2 - \beta^2 \right)  \left( \varepsilon^2 - \gamma^2\right)}  = 1 - \frac{x^2}{\alpha^2 - \varepsilon^2} - \frac{y^2}{\beta^2 - \varepsilon^2} - \frac{z^2}{\gamma^2 -\varepsilon^2}
\ee
When this equation is  differentiated with respect to $\varepsilon$ and  after setting $\varepsilon=\varrho_1, \varrho_2, \varrho_3$,  one finds\,:
\begin{mynewequation} \tag*{[5.49]}
 \begin{aligned}
- \frac{\left(\varrho_1^2 - \varrho_2^2 \right) \left(\varrho_1^2 - \varrho_3^2 \right)}{ \left (\varrho_1^2 - \alpha^2 \right)\left(\varrho_1^2 - \beta^2 \right)\left(\varrho_1^2 - \gamma^2 \right)} =\left( \frac{x}{\alpha^2 - \varrho_1^2}\right)^2 + \left( \frac{y}{\beta^2 - \varrho_1^2}\right)^2 + \left( \frac{z}{\gamma^2 - \varrho_1^2}\right)^2 \\
- \frac{\left(\varrho_2^2 - \varrho_1^2 \right) \left(\varrho_2^2 - \varrho_3^2 \right)}{ \left (\varrho_2^2 - \alpha^2 \right)\left(\varrho_2^2 - \beta^2 \right)\left(\varrho_2^2 - \gamma^2 \right)} =\left( \frac{x}{\alpha^2 - \varrho_2^2}\right)^2 + \left( \frac{y}{\beta^2 - \varrho_2^2}\right)^2 + \left( \frac{z}{\gamma^2 - \varrho_2^2}\right)^2 \\
- \frac{\left(\varrho_3^2 - \varrho_1^2 \right) \left(\varrho_3^2 - \varrho_2^2 \right)}{ \left (\varrho_3^2 - \alpha^2 \right)\left(\varrho_3^2 - \beta^2 \right)\left(\varrho_3^2 - \gamma^2 \right)}
=\left( \frac{x}{\alpha^2 - \varrho_3^2}\right)^2 + \left( \frac{y}{\beta^2 - \varrho_3^2}\right)^2 + \left( \frac{z}{\gamma^2 - \varrho_3^2}\right)^2 \\
\end{aligned} 
\end{mynewequation}
Logarithmical differentiation of the equations by which $x^2,y^2, z^2$ are represented as functions of $\varrho_1^2,\varrho_2^2, \varrho_3^2$ gives\,:
\begin{mynewequation} \tag*{[5.50]}
 \begin{aligned}
- \frac{dx}{d \varrho_1} = \frac{x \varrho_1}{\alpha^2 - \varrho_1^2}, \,\,\,  - \frac{dx}{d \varrho_2} = \frac{x \varrho_2}{\alpha^2 - \varrho_2^2},  \,\,\,  -\frac{dx}{d \varrho_3} = \frac{x \varrho_3}{\alpha^2 - \varrho_3^2}\\
- \frac{dy}{d \varrho_1} = \frac{y \varrho_1}{\beta^2 - \varrho_1^2}, \,\,\,  - \frac{dy}{d \varrho_2} = \frac{y \varrho_2}{\beta^2 - \varrho_2^2},  \,\,\,  -\frac{dy}{d \varrho_3} = \frac{y \varrho_3}{\beta^2 - \varrho_3^2}\\
 \frac{dz}{d \varrho_1} = \frac{z \varrho_1}{\gamma^2 - \varrho_1^2}, \,\,\,  - \frac{dz}{d \varrho_2} = \frac{z \varrho_2}{\gamma^2 - \varrho_2^2},  \,\,\,  -\frac{dz}{d \varrho_3} = \frac{z \varrho_3}{\gamma^2 - \varrho_3^2}
\end{aligned} 
\end{mynewequation}
so that one has\,: 
\begin{mynewequation} \tag*{[5.51]}
  \begin{aligned}
N_1 = \varrho_1^2 \left\{ \left( \frac{x}{\alpha^2 - \varrho_1^2}\right)^2 + \left( \frac{y}{\beta^2 - \varrho_1^2}\right)^2 + \left( \frac{z}{\gamma^2 - \varrho_1^2}\right)^2 \right \}\\
N_2 = \varrho_2^2 \left \{   \left( \frac{x}{\alpha^2 - \varrho_2^2}\right)^2 + \left( \frac{y}{\beta^2 - \varrho_2^2}\right)^2 + \left( \frac{z}{\gamma^2 - \varrho_2^2}\right)^2  \right \}\\
N_3 = \varrho_3^2 \left \{  \left( \frac{x}{\alpha^2 - \varrho_3^2}\right)^2 + \left( \frac{y}{\beta^2 - \varrho_3^2}\right)^2 + \left( \frac{z}{\gamma^2 - \varrho_3^2}\right)^2      \right \}
\end{aligned} 
\end{mynewequation}
From these  relations, one can therefore write\,:
\begin{mynewequation} \tag*{[5.52]}
 \begin{aligned}
N_1 = -\varrho_1^2 \frac{\left(\varrho_1^2 - \varrho_2^2 \right) \left(\varrho_1^2 - \varrho_3^2 \right)}{ \left (\varrho_1^2 - \alpha^2 \right)\left(\varrho_1^2 - \beta^2 \right)\left(\varrho_1^2 - \gamma^2 \right)}\\\\
N_2 = -\varrho_2^2   \frac{\left(\varrho_2^2 - \varrho_1^2 \right) \left(\varrho_2^2 - \varrho_3^2 \right)}{ \left (\varrho_2^2 - \alpha^2 \right)\left(\varrho_2^2 - \beta^2 \right)\left(\varrho_2^2 - \gamma^2 \right)}\\\\
N_3 = - \varrho_3^2 \frac{\left(\varrho_3^2 - \varrho_1^2 \right) \left(\varrho_3^2 - \varrho_2^2 \right)}{ \left (\varrho_3^2 - \alpha^2 \right)\left(\varrho_3^2 - \beta^2 \right)\left(\varrho_3^2 - \gamma^2 \right)}  
\end{aligned} 
\end{mynewequation}
\indent\indent 
These three quantities
are recognised as positive because of
\begin{mynewequation} \tag*{[5.53]}
\alpha  > \varrho_1 > \beta > \varrho_2>\gamma >\varrho_3 >0
\end{mynewequation}
have only to be substituted into equations~(4) 
to obtain the hydrodynamical fundamental equations for elliptical coordinates.\\
\indent\indent
The polar coordinate system $r, \theta, \varphi$, which is determined by\,:
\begin{mynewequation} \tag*{[5.54]}
x = r \cos \theta, \, \, y = r \sin \theta \cos \varphi, \, \,  z= r \sin \theta \sin \varphi
\end{mynewequation}
is orthogonal\,;\hspace{0.4cm}in fact one has\,:
\begin{mynewequation} \tag*{[5.55]}
{\rm d}s^2= {\rm d}x^2 + {\rm d} y^2 + {\rm d} z^2 = {\rm d} r^2 + r^2 {\rm d} \theta^2 + r^2 \sin^2 \theta {\rm d} \varphi
\end{mynewequation}
so that one has
\begin{mynewequation} \tag*{[5.56]}
N_1 = 1 , \,\,\,N_2 = r^2,\,\, N_3 = r^2 \sin^2 \theta .
\end{mynewequation}
\noindent
The density equation (2) becomes\,:
\begin{myequation} \tag*{[5.57]}
\frac{\varrho}{\varrho_0} \begin{vmatrix}
\frac{dr}{dr_0}\,\,\,\,&\frac{d\theta}{dr_0}\,\,\,\,&\frac{d\varphi}{d r_0} \\\\
\frac{dr}{d\theta_0}\,\,\,\,&\frac{d \theta}{d\theta_0}\,\,\,\,&\frac{d \varphi}{d\theta_0} \\\\
\frac{dr}{d\varphi_0}\,\,\,\,&\frac{d \theta}{d\varphi_0}\,\,\,\,&\frac{d \varphi}{d\varphi_0} 
\end{vmatrix}
= \frac{r_0^2 \sin \theta_0}{r^2 \sin \theta}
\end{myequation}
In this case, by setting
\begin{align} \tag*{[5.58]}
 \begin{aligned}
\Phi_1&= \frac{d^2 r}{dt^2} - r \left( \frac{d \theta}{dt}\right)^2  - r \sin ^2 \theta \left(\frac{d \varphi}{dt} \right )^2  \\
\Phi_2&= \frac{d \left( r^2 \frac{ d\theta}{dt}\right)}{dt} - \left( \frac{d \varphi}{dt} \right )^2 r^2 \sin \theta \cos \theta  \, \, \, \,\,  \\
\Phi_3 &= \frac{d \left(r^2 \sin^2 \theta \frac{d \varphi}{dt}\right)}{dt} 
\end{aligned} 
\end{align}
equations~(3) become very simply\,:
\be \tag*{[5.59]}
\frac{d \Omega}{dr} = \Phi_1, \,\,\, \frac{d \Omega}{d \theta} = \Phi_2, \,\,\,\frac{d \Omega}{d \varphi} = \Phi_3, \,\,\, 
\ee
and equations~(4) become\,:
\begin{mynewequation} \tag*{[5.60]}
 \begin{aligned}
\Phi_1 \frac{dr}{dr_0} + \Phi_2 \frac{d \theta}{d r_0} + \Phi_3 \frac{d \varphi}{dr_0} - \frac{d \Omega}{d r_0}\,\,\,=0\\
\Phi_1\,\, \frac{dr}{d\theta_0} + \Phi_2 \frac{d \theta}{d\theta_0} + \Phi_3 \frac{d \varphi}{d\theta_0} - \frac{d \Omega}{d \theta_0}\hspace{0.09cm}=0\\
\Phi_1 \frac{dr}{d \varphi_0} + \Phi_2 \frac{d \theta_0}{d \varphi_0} + \Phi_3 \frac{d \varphi}{d\varphi_0} - 
\frac{d \Omega}{d \varphi_0}=0\\
\end{aligned} 
\end{mynewequation}
\indent\indent
The same transformation can be directly performed in the following way.\hspace{0,4cm} One observes that\,:
\begin{smallequation} \tag*{[5.61]}
\begin{aligned}
\!\!\frac{d^2 x}{dt^2} &= \cos \theta \frac{d^2 r }{d t^2} - r \sin \theta \frac{d^2 \theta}{d t^2} - r \cos \theta \left(\frac{d \theta}{dt} \right )^2 - 2 \sin \theta \frac{dr}{dt } \frac{d \theta}{dt}\,\,\hspace{5.6cm}  \\  \\
\!\!\frac{d^2 y}{dt^2} &= \sin \theta \cos \varphi \frac{d^2 r}{d t^2} + r \cos \theta \cos \varphi \frac{d^2 \theta}{d t^2} - r \sin \theta \sin \varphi \frac{d^2 \varphi}{d t^2} - r \sin \theta \sin \varphi  \frac{d^2  \varphi }{d t^2} - r \sin \theta \cos \varphi \left(\frac{d \theta}{dt} \right )^2  \\
&-r  \sin \theta \cos \varphi \left(\frac{d \varphi}{dt} \right )^2 + 2 \cos \theta \cos \varphi  \frac{dr}{dt } \frac{d \theta}{dt} - 2 \sin \theta \sin \varphi  \frac{dr}{dt } \frac{d \varphi}{dt} - 2 r \cos \theta  \sin \varphi \frac{d \theta}{dt } \frac{d \varphi}{dt}\hspace{1.5cm}  \\  \\
\!\!\frac{d^2 z}{dt^2} &=\sin \theta \sin \varphi \frac{d^2 r }{d t^2} + r \cos \theta  \sin \varphi \frac{d^2 \theta }{d t^2}  + r \sin \theta \cos \varphi \frac{d^2 \varphi}{d t^2} - r \sin \theta \sin \varphi \left(\frac{d \theta}{dt} \right )^2\,\hspace{2.7cm}  \\
&- r \sin \theta \sin \varphi \left(\frac{d \varphi}{dt} \right )^2 + 2 \cos \theta \sin \varphi \frac{dr}{dt } \frac{d \theta}{dt}  + 2 \sin \theta \cos \varphi \frac{dr}{dt } \frac{d \varphi}{dt}  + 2 r  \cos \theta \cos \varphi \frac{d \theta}{dt } \frac{d \varphi}{dt} \hspace{1.3cm}  
\end{aligned}  
\end{smallequation}   
Furthermore one has\,:
\begin{align} \tag*{[5.62]}
\begin{aligned}
\frac{dx}{da} &= \cos \theta \frac{dr}{da}     - r  \sin \theta \frac{d \theta}{da}\hspace{4.2cm}  \\ \\
\frac{dy}{da} &= \sin \theta \cos \varphi \frac{dr}{da } + r \cos \theta \cos \varphi \frac{d \theta}{da } - r \sin \theta \sin \varphi \frac{d \varphi}{da } \\  \\
\frac{dz}{da} &= \sin \theta \sin \varphi \frac{dr}{da } + r \cos \theta \sin \varphi \frac{d \theta}{da } + r \sin \theta \sin \varphi \frac{d \varphi}{da }\,; 
\end{aligned} 
\end{align}  
then the equation\,:
\begin{equation*} \tag*{[5.63]}
\frac{d^2 x}{dt^2} \frac{dx}{da} + \frac{d^2 y}{dt^2} \frac{dy}{da} + \frac{d^2 z}{dt^2} \frac{dz}{da} - \frac{d \Omega}{da}=0
\end{equation*}
becomes\,:
\be \tag*{[5.64]}
 \Phi_1 \frac{dr}{da} + \Phi_2 \frac{d \theta}{da} + \Phi_3 \frac{d \varphi}{da} - \frac{d \Omega}{da}=0 
\ee
where  $\Phi_1,\Phi_2,\Phi_3$ have the meaning as written above.\hspace{0.4cm}In addition to this equation, there are two others\,:
%
\be \tag*{[5.65]}
\begin{aligned}
\Phi_1\frac{dr}{db} +  \Phi_2 \frac{d \theta}{db} + \Phi_3 \frac{d \varphi}{db} - \frac{d \Omega}{db}=0\\\\
\Phi_1\frac{dr}{dc} +  \Phi_2 \frac{d \theta}{dc} + \Phi_3 \frac{d \varphi}{dc} - \frac{d \Omega}{dc}=0
\end{aligned}
\ee
 where $a, b, c$ depend on $r_0, \theta_0, \varphi_0$ through the relations\,:
\be \tag*{[5.66]}
 a=r_0 \cos \theta_0,\,\, b= r_0 \sin \theta_0\cos \varphi_0,\,\, c=r_0 \sin \theta_0 \sin \varphi_0
\ee
The change of variables from 
$a,\,b,\,c$ 
to 
$r_0, \,\theta_0,\,\varphi_0$
into 
 the hydrodynamical equations can be easily carried out by multiplying with the appropriate factors and adding the equations. 
Then, one arrives at the above formulae in a different way. \\
\indent\indent The transformation of the fundamental equations into cylindrical coordinates is extremely simple.\hspace{0.4cm}Namely, 
if one sets\,:
\be \tag*{[5.67]}
x = r \cos \theta, \,\,y = r \sin \theta, \,\,z=z
\ee
then one has\,:
\be \tag*{[5.68]}
{\rm d}s^2= {\rm d}x^2 + {\rm d}y^2 + {\rm d}z^2= {\rm d}r^2 + r^2 {\rm d} \theta^2 + {\rm d}z^2
\ee
so that
\be \tag*{[5.69]}
 N_1 = 1,\,\, N_2= r^2,\,\, N_3=1
\ee
hence\,:
\be \tag*{[5.70]}
\begin{aligned}
\frac{d^2 r}{dt^2} - r  \left(\frac{d \theta}{dt} \right)^2 &= \frac{d \Omega}{dr}\\
\frac{d \left( r^2 \frac{d \theta}{dt} \right) }{dt} &= \frac{d \Omega}{d \theta}\\
\frac{d^2 z}{dt^2}  &= \frac{d \Omega}{dz}
\end{aligned} 
\ee
\indent\indent
In case $r,\, \theta,\, z$ need to be expressed as functions of the initial values $r_0, \,\, \theta_0, \,\, z_0$,  one obtains the equations\,:
\be \tag*{[5.71]}
\begin{aligned}
\left( \frac {d^2 r}{dt^2} - r \left( \frac{d \theta}{dt} \right)^2  \right) \frac{dr}{dr_0} + \frac{d \left(r^2 \frac{d \theta}{dt}\right)}{dt} \frac{d \theta}{d r_0} + \frac{d^2 z}{dt^2} \frac{dz}{d r_0} &= \frac{d \Omega}{dr_0}\\\\
\left( \frac {d^2 r}{dt^2} - r \left( \frac{d \theta}{dt} \right)^2  \right) \frac{dr}{d \theta_0} + \frac{d \left(r^2 \frac{d \theta}{dt}\right)}{dt} \frac{d \theta}{d \theta_0} + \frac{d^2 z}{dt^2} \frac{dz}{d \theta_0} &= \frac{d \Omega}{d \theta_0} \\\\
\left( \frac {d^2 r}{dt^2} - r \left( \frac{d \theta}{dt} \right)^2  \right) \frac{dr}{dz_0} + \frac{d \left(r^2 \frac{d \theta}{dt}\right)}{dt} \frac{d \theta}{d z_0} + \frac{d^2 z}{dt^2} \frac{dz}{d z_0} &= \frac{d \Omega}{dz_0}
\end{aligned} 
\ee
The density equation becomes\,:
\begin{mynewequation}\tag*{[5.72]}
 \left| \begin{aligned}
\frac{dr}{dr_0}\,\,\,\,\,\,&\frac{dr}{d\theta_0}&\frac{dr}{d\varphi_0} \\\\
\frac{d\theta}{dr_0}\,\,\,\,\,&\frac{d\theta}{d\theta_0}&\frac{d\theta}{d\varphi_0} \\\\
\frac{d\varphi}{dr_0}\,\,\,\,\,&\frac{d\varphi}{d\theta_0}&\frac{d\varphi}{d\varphi_0} 
\end{aligned} \right| =
\text{\Large $\frac{r_0}{r}$}
\end{mynewequation}
If we  take the initial conditions and the accelerating forces 
to be symmetric with respect to the $z$ axis, these equations become advantageous\,; 
then we have
$\frac{d \Omega}{d \theta} = 0$, 
therefore
 $\frac{d \left( r^2 \frac{d \theta}{dt}\right)}{dt}=0,$ and thus\,:
\be \tag*{[5.73]}
\frac{d \theta}{dt}=\frac{H}{r^2},
\ee
where $H$ is a time-independent constant which has always the same value for a certain particle, but varies  from particle to particle and has to be determined  by the initial conditions.\hspace{0.4cm}Since $\frac{d \theta}{dt}$ is the rotational velocity of a particle around the $z$ axis,  
the rotational velocity of one and the same particle around the symmetry axis is inversely proportional to  the  
relative distance squared  of the particle to the axis.\hspace{0.4cm}We see from this that no particle initially rotating ceases to rotate under the influence of forces generated by a potential 
and, conversely, no particle begins to rotate if 
it is not initially in rotation.
\deffootnotemark{\textsuperscript{[A.\thefootnotemark]}} \deffootnote{2em}{1.6em}{[A.\thefootnotemark]\enskip}
\hspace{0.4cm}This elegant theorem is thanks to \mbox{S\ws v\ws a\ws n\ws b\ws e\ws r\ws g,} {who obtained it in  cylindrical coordinates from the first \mbox{E\ws u\ws l\ws e\ws r\ws ian} equations which, actually, are for this purpose slightly less convenient than the second \mbox{E\ws u\ws l\ws e\ws r\ws ian} equations.
\footnote{On fluides r\"orelse. Kongl. Vetenskaps-Academiens Handlingar f\"or {\aa}r 1839, S.\,139. Stockholm. Also cf: Sur le mouvement des fluides.\,\, Crelle's Journal Bd. 24, S.\,157, [\hyperlink{Svanberg}{1842}].}

\vspace{0.6cm}
\centerline{\fett{$\S.\,6.$}}
\vspace{0.3cm}
\indent\indent
{Using the $\Omega$ function introduced above (cf. S.\,18),  the hydrodynamical fundamental equations can be written as\,:
\be \tag*{[6.1]}
 \begin{aligned}
\frac{d u}{d t}\frac{dx}{da}  +  \frac{d v}{d t} \frac{dy}{da}+\frac{dw }{d t} \frac{dz}{da} - \frac{d \Omega}{da} =0\\\\
 \frac{d u}{d t}\frac{dx}{db}  +  \frac{d v}{d t} \frac{dy}{db}+\frac{dw }{d t} \frac{dz}{db} - \frac{d \Omega}{db} =0\\\\\frac{d u}{d t}\frac{dx}{dc}  +  \frac{d v}{d t} \frac{dy}{dc}+\frac{dw }{d t} \frac{dz}{dc} - \frac{d \Omega}{dc} =0
\end{aligned} 
\ee
\indent\indent
From these equations, one can easily eliminate  the function $\Omega$ and obtain equations which represent all possible fluid motions under the influence of potential forces.\hspace{0.4cm}The elimination is easily done through differentiations with respect to $a,b,c$ and subtractions 
from which one obtains\,:
\be \tag*{[6.2]}
 \begin{aligned}
\frac{d^2 u}{d t dc}\frac{dx}{db} - \frac{d^2 u}{d t db}\frac{dx}{dc}  +  \frac{d^2 v}{d t dc} \frac{dy}{db} -\frac{d^2 v}{d t db} \frac{dy}{dc}  +\frac{d^2 w }{d t dc} \frac{dz}{db}  - \frac{d^2 w }{d t db} \frac{dz}{dc}    =0\\\\
\frac{d^2 u}{d t da}\frac{dx}{dc} - \frac{d^2 u}{d t dc}\frac{dx}{da}  +  \frac{d^2 v}{d t da} \frac{dy}{dc} -\frac{d^2 v}{d t dc} \frac{dy}{da}  +\frac{d^2 w }{d t da} \frac{dz}{dc}  - \frac{d^2 w }{d t dc} \frac{dz}{da}    =0\\\\
\frac{d^2 u}{d t db}\frac{dx}{da} - \frac{d^2 u}{d t da}\frac{dx}{db}  +  \frac{d^2 v}{d t db} \frac{dy}{da} -\frac{d^2 v}{d t da} \frac{dy}{db}  +\frac{d^2 w }{d t db} \frac{dz}{da}  - \frac{d^2 w }{d t da} \frac{dz}{db}    =0
\end{aligned} 
\ee
One can readily integrate these equations with respect  to time 
by writing each of the three differences in these equations
 as an exact time derivative.\hspace{0.4cm} Denoting the time-independent integration constants as $2A,2B,2C$, one finds\,: 
\begin{equation*}\tag*{(1), [6.3]}
\left.
\begin{aligned}
     \frac{du}{d c}\frac{dx}{db}  - \frac{du}{d b}\frac{dx}{dc}  +   \frac{dv}{d c}\frac{dy}{db}  - \frac{dv}{d b}\frac{dy}{dc} + \frac{dw}{d c}\frac{dz}{db}  - \frac{dw}{d b}\frac{dz}{dc}  =2A\\\\
  \frac{du}{d a}\frac{dx}{dc}  - \frac{du}{d c}\frac{dx}{da}  +   \frac{dv}{d a}\frac{dy}{dc}  - \frac{dv}{d c}\frac{dy}{da} + \frac{dw}{d a}\frac{dz}{dc}  - \frac{dw}{d c}\frac{dz}{da}  =2B\\\\
\frac{du}{d b}\frac{dx}{da}  - \frac{du}{d a}\frac{dx}{db}  +   \frac{dv}{d b}\frac{dy}{da}  - \frac{dv}{d a}\frac{dy}{db} + \frac{dw}{d b}\frac{dz}{da}  - \frac{dw}{d a}\frac{dz}{db}  =2C
\end{aligned} \qquad \right\} 
\end{equation*}
The left sides of these integral equations 
can be seen as the difference of two differential quotients [a curl is meant].
\hspace{0.4cm}Namely, defining\,: \\
\begin{equation*} \tag*{(2), [6.4]}
\left.
\begin{aligned}
\alpha = u \frac{dx}{da} + v \frac{dy}{da} + w \frac{dz}{da}\\\\
\beta = u \frac{dx}{db} + v \frac{dy}{db} + w \frac{dz}{db}\\\\
\gamma = u \frac{dx}{dc} + v \frac{dy}{dc} + w \frac{dz}{dc} 
 \end{aligned} \qquad  \right\} 
\end{equation*}
 one obtains, instead of $(1)$\,: 
\begin{equation*} \tag*{(3), [6.5]} \label{pointer2}
\frac{d \beta}{dc} - \frac{d \gamma}{db} = 2A,\,\,\, \frac{d \gamma}{da} - \frac{d \alpha}{dc} = 2B, \,\,\,\frac{d \alpha}{db} - \frac{d \beta}{da} = 2C
\end{equation*}
These interesting relations are the analogues of the equations that \mbox{C\ws a\ws u\ws c\ws h\ws y}\footnote{%
In an  essay prized  by the Paris Academy: M\'emoire sur la th\'eorie de la propagation des ondes \`a la surface d'un fluide pesant d'une profondeur infinie (M\'em. sav. \'etran.\ Bd.\ 1) [\hyperlink{Cauchy1}{1827}].} 
found already in 1816 [actually, already in 1815] for the first  E\ws u\ws l\ws e\ws r\ws ian dependence.\hspace{0.4cm}They  attained their actual importance 
only when  \mbox{H\ws e\ws l\ws m\ws h\ws o\ws l\ws t\ws z}\footnote{%
Ueber Integrale der hydrodynamischen Gleichungen, welche den Wirbelbewegungen entsprechen, Crelle's Journal, Bd. 55, S.105, [\hyperlink{Helmholtz1}{1858}]. [\textit {English translation}:  On Integrals of the hydrodynamic equations that correspond to vortex motions.
Philos. Mag. {\bf 4}, Vol. 33, 485--511, [\hyperlink{Helmholtz1}{1868}].}
realised their mechanical significance%
\deffootnotemark{\textsuperscript{[T.\thefootnotemark]}}\deffootnote{2em}{1.6em}{[\hspace{0.01cm}T.\thefootnotemark]\enskip}%
\footnote{It appears likely that Helmholtz was not aware of Cauchy's equations but stressed the importance of vortex dynamics.} [of these equations], and thereby laid the foundations for a peculiar treatment of hydrodynamics.\hspace{0.4cm}It will be our next task to investigate this significance using a method adapted to the second  \mbox{E\ws u\ws l\ws e\ws r} dependence; for this
we first need to develop an appropriate theorem.\\\\
\centerline{\fett{$\S.\,7.$}}\\[0.3cm]
\indent\indent
For this purpose
we start from the known theorem that if $\xi$ and $\eta$ are arbitrary  
continuous functions of $x$ and $y$, one has the relation\,:
\be \tag*{[7.1]}
\int \left( \xi {\rm d}x + \eta  {\rm d}y \right) =\iint \left(  \frac{d \xi}{dy} - \frac{d \eta}{dx}\right) {\rm d}x \, {\rm d}y
\ee
where the double integral has to be extended over all elements of a domain on the $xy$ plane, and the simple integral is over the boundary of the domain suitably oriented.%
\deffootnotemark{\textsuperscript{[A.\thefootnotemark]}}\deffootnote{2em}{1.6em}{[A.\thefootnotemark]\enskip}%
\footnote{Cf.\ B.\ Riemann, Lehrs\"atze aus der analysis situs f\"ur die Theorie der Integrale von zweigliedrigen vollst\"andigen Differentialien. Crelle's Journal Bd.\ 54, S.\ 105, [\hyperlink{Riemann1}{1857}].}\\
\indent\indent
One can  generalise this theorem in the following way\,:
Let there be an arbitrary
closed curve in space and consider the path integral over all elements of this curve\,:
\be \notag
\int \left( \xi\, {\rm d}x + \eta\, {\rm d}y + \zeta\, {\rm d}z \right)
\ee
where  $\xi,\,\eta,\,\chi$ are  arbitrary continuous functions of $x,y,z$.\hspace{0.4cm}Let us think  of an arbitrary  connected surface limited by this curve, so that one can consider on this whole surface  $z$ as a function of $x$ and $y$ and set\,:
\be \tag*{[7.2]}
{\rm d}z = \frac{dz}{dx} {\rm d}x + \frac{dz}{dy} {\rm d}y, 
\ee
from which the given integral transforms into\,:
\be \notag
\int \left\{ \left (\xi + \frac{dz}{dx} \zeta \right) {\rm d}x   + \left (\eta + \frac{dz}{dy} \zeta \right) {\rm d}y  \right \}
\ee
By the {stated theorem for [the case of] two independent variables, this integral becomes\,:
\be \notag
\iint \left\{ \frac{ d\left(\xi + \frac{dz}{dx} \zeta \right)}{dy} - \frac{\left(\eta + \frac{dz}{dy} \zeta \right)}{dx}\right\} {\rm d}x {\rm d}y
\ee
Since  $z$  is a function of $x$ and $y$, one has\,:
\be \tag*{[7.3]}
 \begin{aligned}
\frac{ d\left(\xi + \frac{dz}{dx} \zeta \right)}{dy} &= \frac{d \xi}{dy} + \frac{d \xi}{dz} \frac{d z}{dy} + \frac{d \zeta}{dy} \frac{d z}{dx} + \frac{d \zeta}{dz} \frac{d z}{dx}\frac{d z}{dy} + \frac{d^2 z}{dx dy} \zeta\\
\frac{ d\left(\eta + \frac{dz}{dy} \zeta \right)}{dx} &= \frac{d \eta}{dx} + \frac{d \eta}{dz} \frac{d z}{dx} + \frac{d \zeta}{dx} \frac{d z}{dy} + \frac{d \zeta}{dz} \frac{d z}{dy}\frac{d z}{dx} + \frac{d^2 z}{dy dx} \zeta
\end{aligned} 
\ee
and so one has the equation\,:
\begin{smallequation} \tag*{[7.4]}
\int \left( \xi {\rm d}x + \eta {\rm d}y + \zeta {\rm d}z \right) = \iint \left \{ \left(\frac{d \xi} {dy} - \frac{d \eta}{dx} \right) + \left(\frac{d \zeta} {dy} - \frac{d \eta}{dz} \right)\frac{dz}{dx} +  \left( \frac{d \xi}{dz} - \frac{d \zeta}{dx}\right) \frac{dz}{dy}  \right\}{\rm d}x {\rm d}y,
\end{smallequation}
where now the double integral has to be extended over all the elements of the surface
limited by the curve.\hspace{0.4cm}This otherwise arbitrary surface through the curve has just  to satisfy the condition that the part embedded in  the curve is not multiply-connected  and that the curve forms a complete boundary to it.\hspace{0.4cm}Let $\lambda,\mu, \nu$ be the angles of the normal drawn to the surface with the coordinates axes,  one has\,: 
\be \tag*{[7.5]}
\frac{dz}{dx} = - \frac{\cos \lambda}{\cos \nu}, \,\,\, \frac{dz}{dy} = - \frac{\cos \mu}{\cos \nu};
\ee
from which follows that one has
\begin{equation*}\tag*{(1), [7.6]}
\begin{aligned}
 &\hspace{4cm}\int \left( \xi {\rm d}x + \eta {\rm d}y +\zeta {\rm d}z \right)=\\
  &\int \left \{ \left(\frac{d \eta} {dz} -  \frac{d \zeta}{dy} \right)\cos \lambda + \left(\frac{d \zeta} {dx} - \frac{d \xi}{dz} \right) \cos \mu +  \left( \frac{d \xi}{dy} - \frac{d \eta}{dx}\right) \cos \nu \right\} {\rm d} \sigma
\end{aligned}
\end{equation*}
where $\displaystyle \frac{ {\rm d}x\,{\rm d}y}{\cos \nu} =  {\rm d} \sigma$ denotes the surface element. \\
\indent\indent 
For our purposes, we can 
bring
both of the above integrals into more convenient  forms\,:\\\\
\indent\indent
If one determines three angles  $\lambda', \, \, \mu', \,\,\nu'$  so that\,:
\be \tag*{[7.7]}
\cos \lambda': \cos \mu': \cos \nu'=\left(\frac{d \eta} {dz} -  \frac{d \zeta}{dy} \right) : \left(\frac{d \zeta} {dx} - \frac{d \xi}{dz} \right) : \left( \frac{d \xi}{dy} - \frac{d \eta}{dx} \right) 
\ee
and, at the same time, assumes that they are the angles of a definite direction with the  coordinates axes, so that\,:
\be \tag*{[7.8]}
\cos^2 \lambda' +  \cos^2 \mu' + \cos^2 \nu' =1
\ee
one finds\,:
\be \tag*{[7.9]}
 2 \Delta \cos \lambda' =\frac{d \eta} {dz} -  \frac{d \zeta}{dy},\,\,
 2 \Delta \cos \mu' =  \frac{d \zeta} {dx} - \frac{d \xi}{dz}, \,\,
 2 \Delta \cos \nu' = \frac{d \xi}{dy} - \frac{d \eta}{dx} 
\ee
where\,:
\be \tag*{[7.10]}
 4 \Delta^2 =  \left(\frac{d \eta} {dz} -  \frac{d \zeta}{dy} \right) ^2  + \left(\frac{d \zeta} {dx} - \frac{d \xi}{dz} \right) ^2 + \left( \frac{d \xi}{dy} - \frac{d \eta}{dx} \right)^2  
\ee
The integral above, whose element is ${\rm d}\sigma$, now goes into 
\be \notag
2 \int \Delta \left(\cos \lambda \cos \lambda'  + \cos \mu \cos \mu' + \cos \nu \cos \nu' \right) {\rm d}\sigma
\ee
and one obtains\,:
\begin{equation*} \tag*{(2), [7.11]}
\int \left\{\left( \frac{d \eta}{dz}  - \frac{d \zeta}{dy} \right) \cos \lambda
+ \left( \frac{d \zeta}{dx}  - \frac{d \xi}{dz} \right) \cos \mu + \left( \frac{d \xi}{dy}  - \frac{d \eta}{dx} \right) \cos \nu \right \}{\rm d}\sigma = 2\int \Delta \cos \theta {\rm d}\sigma \,,
\end{equation*}
where $\theta$ is the angle between the   directions determined by $\lambda, \mu, r$ and
$\lambda',\mu', r' $ .\\
\indent\indent
Indicating by $\varrho, \sigma, \tau$ the angles of the coordinate axes with the tangent to the curve at the point $x,y,z$,  one has\,:
\be \tag*{[7.12]}
{\rm d}x= \cos \varrho {\rm d}s, \,\,\, {\rm d}y = \cos \sigma {\rm d}s, \, \, \,  {\rm d}z= \cos \tau {\rm d}s, 
\ee
where ${\rm d}s$ indicates the arc element.\hspace{0.4cm}Let now $ \varrho', \sigma', \tau'$ be the angles of a 
direction with the coordinate axes, 
so that\,:
\be \tag*{[7.13]}
\cos^2 \varrho' +  \cos^2 \sigma' + \cos^2 \tau' = 1,
\ee
and furthermore [require that]
\be \tag*{[7.14]}
\cos \varrho': \cos \sigma': \cos \tau'= \xi:\eta:\zeta,
\ee
so one has\,:
\be \tag*{[7.15]}
U^2 = \xi^2 + \eta^2 + \zeta^2
\ee
and
\be \tag*{[7.16]}
U \cos \rho'=\xi,\,\, U \cos \sigma' = \eta, \,\,U \cos \tau' = \zeta.
\ee
These values of $\xi,\,\eta, \,\zeta$ and ${\rm d}x$, ${\rm d}y$, ${\rm d}z$ substituted into the simple integral, give\,:
\be \tag*{[7.17]}
\int \left( \xi {\rm d} x +  \eta {\rm d}y + \zeta  {\rm d}z\right)\, = \int U (\cos \varrho \cos \varrho' +\cos \sigma \cos \sigma' + \cos \tau \cos \tau' ) {\rm d}s
\ee
or
\begin{equation*}\tag*{(3), [7.18]}
\int\left( \xi {\rm d} x +  {\eta\rm d} y + \zeta  {\rm d}z \right)\, =\,\int U \cos \theta'  {\rm d}s
\end{equation*}
where $\theta'$ indicates the angle between the directions determined by $\varrho,\sigma,\tau$ and $\varrho',\sigma',\tau'$. \\[1cm]
\centerline{\fett{$\S.\,8.$}}\\[0.3cm]
\indent\indent
After development of this lemma we come back to the task 
described at the end of $\S\,6$\,:\\
\indent\indent
According to $(1)$ of $\S.\,7.$, one can write the equation\,:
\begin{smallequation} \tag*{[8.1]}
\int\left( \alpha  {\rm d}a + \beta  {\rm d}b + \gamma  {\rm d}c \right)\, =\, \int \left \{ \left( \frac{d \beta}{dc} -
\frac {d \gamma}{db} \right)\cos \lambda  + \left( \frac{d \gamma}{da} -
\frac {d \alpha}{dc} \right)\cos \mu + \left(\frac{d \alpha}{db} -
\frac {d \beta}{da} \right) \cos \nu \right \}  {\rm d} \sigma_0
\end{smallequation}
where the first [l.h.s.] integral is over a closed curve, the second [r.h.s.] over an arbitrary simply-connected surface bounded by that curve, and where $\lambda, \mu, \nu$,  are the angles of the normal to that surface with the coordinate axes,  
${\rm d} \sigma_0$ is the surface element  and $\alpha,\beta,\gamma  $ are arbitrary  functions of $ a,b,c $.\hspace{0.4cm}According to $(2)$ of the $\S\,7$, the second  integral  can be transformed, so that\,:

\begin{equation*} \tag*{(1), [8.2]}
\int\left( \alpha\,{\rm d}a + \beta\, {\rm d}b + \gamma \, {\rm d}c \right)\, =\, 2 \int \Delta_0 \cos \theta_0 {\rm d} \sigma_0 \,.
\end{equation*}
If one now takes $\alpha, \beta,\gamma $ to be the quantities defined  in $(2)$ of $\S\,6$ and makes use of equations~(3) of $\S\,6$, one obtains\,:
\be \tag*{[8.3]}
\Delta_0 = \sqrt{A^2 + B^2 + C^2} \,.
\ee
As to $\theta_0$, it is the angle  between the directions determined by 
 $\lambda, \mu, \nu$ and  $\lambda', \mu', \nu'$, the latter being determined by\,:
\be \tag*{[8.4]}
\Delta_0 \cos \lambda' = \Delta, \, \, \Delta_0 \cos \mu' = B, \, \, \Delta_0 \cos \nu' = C \,.
\ee
The integral on the right-hand side of
the preceding  equation $(1)$ does not depend on time,  
thus
one may choose an arbitrary value of $t$ on its left-hand side;
if one sets $t=0$, 
then
$\alpha, \beta, \gamma $
 go into the initial values of $u,\,v,\,w$
 that we denote  with $u_0, v_0, w_0$. 
One then obtains\,:
\be \tag*{[8.5]}
\int \left ( u_0 {\rm d}a + v_0 {\rm d}b + w_0 {\rm d}c \right) = 2 \int \Delta_0 \cos \theta_0 {\rm d} \sigma_0
\ee
According to $(3)$ of the $\S.\,7.$, the first integral  may be transformed  and thus one finds\,:                      
\begin{equation*} \tag*{(2)\,, [8.6]}
\int U_0  \cos \theta'_0 {\rm d} s_0 = 2 \int \Delta_0 \cos \theta_0 {\rm d} \sigma_0
\end{equation*}
where 
\be \tag*{[8.7]}
U_0 = \sqrt{u_0 ^2 + v_0^2 +w_0^2}
\ee
is the initial velocity, $\theta'_0$ is  the angle between $U_0$ and the element ${\rm d}s_0$ of the closed curve.\\
\indent\indent
Let us assume, that the given closed curve 
is a circle with an infinitely small radius $r$ 
and the surface spanned by the curve is the disk of the circle. So it is clear that $\Delta_0 \cos \theta_0 $ for  different points of the circle surface changes only by infinitely small quantities.\hspace{0.4cm} Then, we have  $\int \Delta_0\, \cos \theta_0 \,{\rm d} \sigma  =\pi r^2 \cdot \Delta_0 \cos \theta_0$.\hspace{0.4cm}The component of the initial velocities with respect to the tangent to this circle is $U_0 \cos \theta_0'$, a quantity, that for differents points of the circumference will take different values.\hspace{0.4cm}However,   $U_0 \cos \theta_0'$ in each point may be considered as the sum of two velocities  $T_0 + T_0'$, where $T_0'$ is the progressive motion projected to  a tangent, which is common to all points of the circle and $T_0$ the tangential velocity by the rotation around the  center of the infinitely small circle.\hspace{0.4cm}The progressive motion of all particles is the same\,: as a consequence its component $T'_0$ with respect to the  tangent of the circle will be the same  for two particles diametrically opposite, but of opposite signs, so that $\int T_0' {\rm d}s_0=0$, when one integrates over the whole circle.\hspace{0.4cm}The relative velocity of the particles will be the same except for infinitely small quantities provided the velocities are  continous functions of the position, which  we always suppose.\hspace{0.4cm}So, it  follows that $\int T_0 {\rm d}s_0=T_0 \cdot \int {\rm d}s_0=2\pi r \cdot T_0$.\hspace{0.4cm}From the equation (2) it  follows now $2\pi r \cdot T_0= 2 \cdot  \pi r^2 \cdot \Delta_0 \cos \theta_0$ or\,:
\be \tag*{[8.8]}
\Delta_0 \cos \theta_0= \frac{T_0}{r}
\ee
Since $T_0$ is the tangential velocity, so $T_0\,: r$ is the rotational velocity  around the nfinitely small distant center of the circle, or --- as one can say, in order to take into account, also the location and orientation of the circle --- is the rotational velocity around the normal to  the surface of the circle taken as axis.
Since $\theta_0$ is the angle between the normal to the  surface element and the direction determined by the angles $\lambda',\mu',\nu'$, then $\Delta_0$ in a point  becomes the rotational velocity  around an axis passing by this point and oriented in the direction  $\lambda',\mu',\nu'$.\\
\indent\indent
Furthermore  $A,\,\,B,\,\,C$ are the components of the rotational velocity of a particle $a,\,b,\,c$
around axes, which are parallel to the coordinates axes through the point $a,\,b,\,c$. 
\\[1cm]
\centerline{\fett{$\S\,9.$}}\\[0.3cm]
\indent\indent
Analogously to $(2)$ of $\S.\,8$,  one has at each time $t$\,:
\begin{equation*} \tag*{(1), [9.1]}
\int U \cos \theta'  {\rm d} s = 2 \int \Delta \cos \theta  {\rm d}\sigma
\end{equation*}
where  now $ U \cos \theta' $ is  the component of the velocity of a particle $x, y, z$ with respect to the tangent to the closed curve, over whose element ${\rm d}s$ the first integral has  to be extended, and where $\Delta \cos \theta$ is  the component of the rotational velocity expressed with respect to the normal to  the element ${\rm d} \sigma$ of the surface, which is bounded by that curve.\\
\indent\indent
According to $(3)$ of $\S.\,7$, one has
\be \tag*{[9.2]}
\int U \cos \theta' {\rm d}s = \int \left( u{\rm d}x + v {\rm d}y + w {\rm d}z  \right)
\ee
but, since  from $ (2) $ of $\S.\,6$ one easily infers that
\be \tag*{[9.3]}
\alpha {\rm d} a + \beta {\rm d} b + \gamma {\rm d} c = u {\rm d}x + v {\rm d}y + w {\rm d}z,
\ee
one has\,:
\be \tag*{[9.4]}
\int U \cos \theta' {\rm d}s = \int \left(\alpha {\rm d}a + \beta {\rm d}b + \gamma {\rm d}c  \right)
\ee
The second integral is calculated according to $(1)$  of 
$\S.\,8$\,:
\be \tag*{[9.5]}
\int U \cos \theta' {\rm d}s = 2 \int \Delta_0 \cos \theta_0 {\rm d} \sigma_0 \,.
\ee
Hence, because of equation $(1)$, one obtains\,:
\begin{equation*} \tag*{(2), [9.6]}
\int \Delta \cos \theta {\rm d}\sigma = \int \Delta_0 \cos \theta_0  {\rm d} \sigma_0
\end{equation*}
Herefrom follows that $\int \Delta \cos \theta' {\rm d}\sigma$ is constant with respect to time, provided the integral is always extended over a surface moving with the flow
and consisting of the same particles.
Following \mbox{H\ws e\ws l\ws m\ws h\ws o\ws l\ws t\ws z,} if one designates by \mbox{r\ws o\ws t\ws a\ws t\ws i\ws o\ws n\ws a\ws l} \,\mbox{i\ws n\ws t\ws e\ws n\ws s\ws i\ws t\ws y} the product of the rotational velocity around the normal to the surface as an axis, times the size of the surface element,  one obtains the following result\,:\\
\indent\indent
The integral of the rotational intensity  over a surface always formed by the same particles remains unchanged in time.
\\
\indent\indent
Since this theorem is valid, 
irrespective of how small  the surface may be,  it is also valid for each single surface element\,:  the rotational intensity of a surface element always stays
the same.\hspace{0.4cm}Because  such a particle cannot spread out infinitely, its rotational velocity cannot decrease infinitely. It follows herefrom that no particle once put into  rotational motion 
can stop  rotating\,; and on the other hand, one easily sees that no particle which at initial time is not rotating,  may ever begin to rotate.\\
\indent\indent
One must remark that these results are obtained under the assumption that the accelerating forces acting on the fluid are partial derivatives of a potential function. If, however, the accelerating forces do not possess this property,  these theorems do not apply.\hspace{0.4cm}Herefrom one obtains a criterion to know whether  accelerating forces acting on a fluid without pressure forces, have or do not have a potential.{ In the latter case, the rotational intensity of single particles is not conserved; in general new particles begin rotating, and rotating particles will lose this characteristic motion.\\
\indent\indent  
So far we have considered only rotating surface elements, but let us take into account a mass element which is contained in a cylinder whose axis is  the rotation axis. Its constant mass is the product of  its transverse section, its length and its density.\hspace{0.4cm}Since the product of the rotational velocity by the  transverse section is constant,  one sees that  for each element  the  ratio of its rotational velocity to   the product of the distance measured in the direction of its rotation axis by the density is constant.\label{pointer}\hspace{0.4cm} Therefore, if the fluid is liquid, i.e., its density  considered as constant,  the ratio of the rotational velocity to  the length of the particle is constant.\\[1cm]
\centerline{\fett{$\S.\,10.$}}\\[0.3cm]
\indent\indent
In his theory of rotational motion, 
\mbox{H\ws e\ws l\ws m\ws h\ws o\ws l\ws t\ws z}  introduced this following important principle\,: instead of considering the whole rotating mass,  one should fragment it into \mbox{v\ws o\ws r\ws t\ws e\ws x}  \mbox{l\ws i\ws n\ws e\ws s.} Here, a vortex line is a line lying  in the  flow so that its direction will stay always parallel to the instantaneous rotation axis.}
By v\ws o\ws r\ws t\ws e\ws x   f\ws i\ws l\ws a\ws m\ws e\ws n\ws t, we understand the infinitely thin cylinder which, wrapped around the vortex line,  includes the rotating particles.\\
\indent \indent
Denoting by ${\rm d}a, {\rm d}b, {\rm d}c$ [the three components of]} an element of such a vortex line at time $t=0$, one obviously has\,:
\begin{equation*} \tag*{(1), [10.1]}
{\rm d}a : {\rm d}b :{\rm d}c = A\,:\,B\,:\,C
\end{equation*}
Let $\varphi$ and $\psi$ be  functions of $a,\, b,\,c$, such that the vortex lines 
at time $t=0$ are obtained by setting these two functions to constant values. 
Then the following conditions must hold\,:
\be \tag*{[10.2]}
 \begin{aligned}
\frac{d \varphi}{da} {\rm d}a + \frac{d \varphi}{db} {\rm d}b +  \frac{{\rm d}\varphi}{dc} {\rm d}c = 0\\\\
\frac{d \psi}{da} {\rm d}a + \frac{d \psi}{db} {\rm d}b + \frac{d \psi}{dc} {\rm d}c = 0\\
\end{aligned} 
\ee
or, according to $(1)$\,:
\be \tag*{[10.3]}
 \begin{aligned}
\frac{d \varphi}{da} A + \frac{d \varphi}{db} B +  \frac{d \varphi}{dc} C &= 0\\
\frac{d \psi}{da} A + \frac{d \psi}{db} B + \frac{d \psi}{dc} C &= 0 \,.\\
\end{aligned} 
\ee
If $A,\,B,\,C$ were} known, one could find  $\varphi$ and $\psi$ from these equations by integration. 
\hspace {0.4cm}One can easily observe that $\varphi, \,\psi$ must be such that, one may write\,:
\begin{equation*} \tag*{(2), [10.4]}
-2 A = \begin{vmatrix} \frac{d \varphi}{db}& \frac{d \varphi}{dc} \\\\
\frac{d \psi}{db}& \frac{d \psi}{dc} \end{vmatrix},\\
\,\,-2 B = \begin{vmatrix} \frac{d \varphi}{dc}& \frac{d \varphi}{da} \\\\
\frac{d \psi}{dc}& \frac{d \psi}{da} \end{vmatrix},\\
\,\,-2 C = \begin{vmatrix} \frac{d \varphi}{da}& \frac{d \varphi}{db} \\\\
\frac{d \psi}{da}& \frac{d \psi}{db} \end{vmatrix} \,.\\
\end{equation*}
Indeed, from the preceding  equations  follows that $A, B, C$ are proportional to  these determinants\,; however, by substitution of~(2) into~(3) of $\S.\,6$,
 the following relation is identically satisfied\,:
\begin {equation*} \tag*{(3), [10.5]}
\frac{dA}{da} + \frac{dB}{db} +\frac{dC}{dc} =0 \,.
\end{equation*}
Thus, $\varphi,\, \psi$ are always determined by integration in such a way that equations~(2) are satisfied.\\
\indent\indent 
From these considerations it follows
 that rotating particles may always be considered as arranged in vortex filaments. \hspace{0.4cm} When such a vortex filament moves along  with the fluid, always the same  particles will belong to it, as  no particle belonging to the vortex filament may lose its rotational  motion and, furthermore,  each element parallel to the rotation axis of the vortex filament always remains parallel.
Namely, a surface element, which, at some time, is parallel to the rotation axis, will have a vanishing rotational intensity with respect to the direction of its normal; since the rotational intensity always stays the same, the normal  to the surface element always stays perpendicular to the rotation axis, because the particles themselves always maintain their rotational motion around the axis of the vortex filament.\hspace{0.4cm}Therefore, one obtains the equations of the vortex lines at time $t$, if one expresses the values of $a,\,b,\,c$ in $\varphi$ and $\psi$ through $x,\,y,\,z,\,t$.\hspace{0.4cm}The equations of the vortex lines at time $t$ must be such
that $\displaystyle \frac{{\rm d} \varphi}{{\rm d}t} =0$ and $\displaystyle  \frac{{\rm d}\psi}{{\rm d}t} =0$ 
are satisfied. 
Here one has to differentiate with respect to $t$, whether it appears explicitly 
or implicitly, through $x,\,y,\,z$.
\hspace{0.4cm}
In motion, the vortex filament will sometimes grow in density and size and sometimes decrease, namely in such a way that the length of a vortex-tube element multiplied by the density remains proportional to the rotational velocity. 
\mbox{(Cf.\ S.\,40 [now page~\pageref{pointer}].)}\\
\indent\indent
The rotational intensity 
 in each  section of this vortex filament  will remain unchanged over time.\hspace{0.2cm} Moreover,  it  is also the same  for all  sections. 
\hspace{0.4cm}Obviously we have $\int \Delta \cos \theta {{\rm d} \sigma}=0$ if the integral extends over a closed surface.\hspace{0.4cm} 
Namely, let us think of a closed curve  drawn on this surface. In order to extend  $\int \Delta \cos \theta {{\rm d}\sigma}$ over {the two portions of the surface
separated by the curve, one has to integrate $\int U \cos \theta' {\rm d}s$ twice 
but in the opposite direction over this curve, so that one has  $\int \Delta \cos \theta {{\rm d}\sigma}=0$.\hspace{0.4cm}
Let this closed surface be formed by two sections of a vortex filament and the part of the filament surface  between them,
then $\int \Delta \cos \theta {{\rm d}\sigma}$ 
will disappear for the latter part, and so 
will the sum of the integrals for both sections. 
Indeed, since  $\cos \theta$, in a section, is $+1$ [and]  $-1$ in the other one, the rotational intensity with regard to the axis of the vortex filament is the same for each section. \\
\indent\indent
The results of the preceding investigations about these particular rotational motions are essentially due to  \mbox{H\ws e\ws l\ws m\ws h\ws o\ws l\ws t\ws z,} who, by restricting himself to liquid flows, derived them  in a classical essay\footnote{Crelle's Journal, Bd.\ 55, 
S.\,33 and ff
 [\hyperlink{Helmholtz1}{1858}].} from the \mbox{E\ws u\ws l\ws e\ws r\ws }ian first form of the equations.
%

\indent\indent
 A special case of the theorem on the constancy of the rotational velocity has been already  given  by \mbox{S\ws v\ws a\ws n\ws b\ws e\ws r\ws g\,:}\\
\indent\indent
Namely, when all motions happen symmetrically around an axis and circular vortex  filaments centered on the symmetry axis} are given in the flow with  radius $r$ and  infinitely small section $\omega$, 
such a vortex filament may be considered as a cylinder of height $2 \pi r$ and basis $\omega$  so that its content will be $ 2 \pi  r \omega$.\hspace{0.4cm}Since the rotational motion of the particles is conserved,  the content $2 \pi r \omega$ for each vortex filament  must be constant in time.\hspace{0.4cm}Then $\Delta \omega r$ must be proportional  to the rotational velocity $\Delta$.\hspace{0.4cm}Since $\Delta \omega$, as rotational intensity, is constant, so  the rotational velocity $\Delta$ for each vortex ring at different times must be proportional to its radius.\\
\indent\indent
This is the meaning  of the formulae derived by \mbox{S\ws v\ws a\ws n\ws b\ws e\ws r\ws g} from the first \mbox{E\ws u\ws l\ws e\ws r\,ian} form\footnote {Crelle's Journal Bd.\ 24. S.\ 159. Nro.\ 31 [\hyperlink{Svanberg}{1842}].} of the fundamental equations.\\[1cm]
\centerline{\fett{$\S \,11.$}}
\\[0.3cm]
\indent\indent
The concept of vortex lines gives rise to an interesting transformation of the hydrodynamical equations in their first form: Namely, by introducing as dependent variables, functions which are in a certain relation to $u,v,w$.\\
\indent\indent From equations (2) of  $\S \,10$  follows indeed\,:
\be \tag*{[11.1]}
 \begin{aligned}
\frac{d\gamma}{db} - \frac{d\beta }{dc}= \frac{d\varphi}{db}\frac{d\psi}{dc} - \frac{d\varphi}{dc}\frac{d\psi}{db}\\\\
\frac{d\alpha}{dc} - \frac{d\gamma }{da}= \frac{d\varphi}{dc}\frac{d\psi}{da} - \frac{d\varphi}{da}\frac{d \psi}{dc}\\\\
\frac{d \beta}{da} - \frac{d\alpha }{db} = \frac{d\varphi}{da}\frac{d\psi}{db} - \frac{d\varphi}{db}\frac{d\psi}{dc}
\end{aligned} 
\ee
provided that
$\varphi$ and $\psi$ 
are
functions of $a,b,c$ that represent the vortex line, when set equal to constants.
Herefrom follows that one can represent $\alpha, \beta, \gamma$ in the 
 form\,:
\begin{equation*} \tag*{(1), [11.2]}
\alpha = \frac{dF}{d a} + \varphi \frac{d \psi}{d a}, \,\,\,\, \beta= \frac{dF}{d b} + \varphi \frac{d \psi}{d b}, \,\,\,\,\gamma = \frac{dF}{d c} + \varphi \frac{d \psi}{d c}
\end{equation*}
where in general $F$  will be  function of $a,\,b,\,c,\,t$.\hspace{0.4cm}Herefrom one finds\,:
\be \tag*{[11.3]}
 \begin{aligned}
\alpha \frac{da}{dx} + \beta \frac{db}{dx} + \gamma \frac{dc}{dx} = \frac{dF}{dx} + \varphi \frac{d \psi}{dx}\\\\
\alpha \frac{da}{dy} + \beta \frac{db}{dy} + \gamma \frac{dc}{dy} =  \frac{dF}{dy} + \varphi \frac{d \psi}{dy}\\\\
\alpha \frac{da}{dz} + \beta \frac{db}{dz} + \gamma \frac{dc}{dz} =   \frac{dF}{dx} + \varphi \frac{d \psi}{dz}
\end{aligned} 
\ee
\hspace{0.6cm}By solving  the system (2) of $\S\,6$ using the relations $(3)$ of $\S\,2$, one finds\,:
\be \tag*{[11.4]}
 \begin{aligned}
\alpha \frac{da}{dx} + \beta \frac{db}{dx} + \gamma \frac{dc}{dx} &= u\\\\
\alpha \frac{da}{dy} + \beta \frac{db}{dy} + \gamma \frac{dc}{dy} &=  v\\\\
\alpha \frac{da}{dz} + \beta \frac{db}{dz} + \gamma \frac{dc}{dz} &=  w
\end{aligned} 
\ee
so that one can always set\,:
\begin{equation*} \tag*{(2), [11.5]}
u = \frac{dF}{dx} + \varphi \frac{d \psi} {dx}, \,\,\, v= \frac{dF}{dy} + \varphi \frac{d \psi}{dy},\,\,\, w= \frac{dF}{dz} + \varphi \frac{d \psi}{dz}
\end{equation*}
where $\varphi= Const$ and $\psi=Const$ are the equations of the vortex lines.\hspace{0.4cm} If one replaces the quantities $a,\,b,\,c$ in $\varphi$ and $\psi$ in terms of 
$x,y,z$ [at time] $t$, one obtains\;:
\begin{equation*} \tag*{[11.6]}
\frac{{\rm d} \varphi}{{\rm d}t} = 0, \qquad \frac{{\rm d} \psi}{{\rm d}t}=0
\end{equation*}
whereby differentiation has to be performed when $t$ appears explicitly as well as implicitly in $x,y,z$.
Hence we have\,:
\be \tag*{[11.7]}
 \begin{aligned} 
\frac{d \varphi}{dt} + \frac{d \varphi}{dx} u + \frac{d \varphi}{dy} v +\frac{d \varphi}{dz} w &=0\\\\
\frac{d \psi}{dt} + \frac{d \psi}{dx} u + \frac{d \psi}{dy} v +\frac{d \psi}{dz} w &=0 \
\end{aligned} \,.
\ee
Substituting the values of $u,\,v,\,w$ [taken from~(2)], we obtain\,:
\begin{equation*}\tag*{(3), [11.8]}
\left.
\begin{aligned} \frac{d \varphi}{dt} + \frac{d \varphi}{dx} \left (\frac{dF}{dx}+ \varphi \frac{d \psi}{dx}\right)  + \frac{d \varphi}{dy} \left (\frac{dF}{dy}+ \varphi \frac{d \psi}{dy}\right) + \frac{d \varphi}{dz} \left (\frac{dF}{dz}+ \varphi \frac{d \psi}{dz}\right)  =0 \\\\
\frac{d \psi}{dt} + \frac{d \psi}{dx} \left (\frac{dF}{dx}+ \varphi \frac{d \psi}{dx}\right)  + \frac{d \psi}{dy} \left (\frac{dF}{dy}+ \varphi \frac{d \psi}{dy}\right) + \frac{d \psi}{dz} \left (\frac{dF}{dz}+ \varphi \frac{d \psi}{dz}\right)  =0
 \end{aligned} \qquad \right\}
\end{equation*}
In addition to these relations there is also the density equation $(5)$ of  $\S.\,2.$
\begin{equation*} \tag*{[11.9]}
\frac{d \varrho}{dt} + \frac{d }{dx} \varrho \left (\frac{dF}{dx}+ \varphi \frac{d \psi}{dx}\right)  + \frac{d }{dy}  \varrho \left (\frac{dF}{dy}+ \varphi \frac{d \psi}{dy}\right) + \frac{d }{dz}  \varrho \left (\frac{dF}{dz}+ \varphi \frac{d \psi}{dz}\right) =0 
\end{equation*}
We observe that, in general, these  three equations are not sufficient to determine  the four unknown functions $F,\varphi,\psi$ and  $\varrho$. Only  in the case of liquid  fluids, where  $\varrho$ is constant, 
are these sufficient, insofar as the latter expression transforms into\,:
\begin{equation*} \tag*{(4), [11.10]}
\frac{d}{dx} \left(\frac{dF}{dx}+ \varphi \frac{d \psi}{dx} \right)  + \frac{d}{dy} \left(\frac{dF}{dy}+ \varphi \frac{d \psi}{dy} \right) + \frac{d}{dz} \left(\frac{dF}{dz}+ \varphi \frac{d \psi}{dz} \right)  =0
\end{equation*}
These are the transformed equations in the elegant form that \mbox{C\ws l\ws e\ws b\ws s\ws c\ws h} found by another way,
without giving the meaning of $\varphi$ and $\psi$.\footnote{%
Ueber die Integration der hydrodynamischen Gleichungen  (Crelle's Journal Bd.\ 56, S.\ 1 [\hyperlink{Clebsch1}{1859}]). 
In this essay   \mbox{C\ws l\ws e\ws b\ws s\ws c\ws h} starts from the first form of the \mbox{E\ws u\ws l\ws e\ws r} equation to perform  the above transformation. On the one hand,  the method applied by him  is susceptible  of a bigger simplification; on the other hand,
one can give a procedure quite analogous to that used above to prove that
$\displaystyle\frac{{\rm d} \varphi} {dt}= 0$ and $\displaystyle\frac{{\rm d} \psi} {dt}= 0.$\hspace{0.4cm} 
\mbox{C\ws l\ws e\ws b\ws s\ws c\ws h} showed (in Theorem 2) that the equations transformed by the introduction of $F, \varphi, \psi$ are the conditions  for the variation of a quadruple integral to disappear;
this is [somewhat] similar to what  \mbox{C\ws l\ws e\ws b\ws s\ws c\ws h} did for steady-state flows in his paper\,: Ueber eine allgemeine Transformation der hydrodynamischen Gleichungen  (Crelle's Journal Bd.\ 54, S.\ 301 [\hyperlink{Clebsch2}{1857}]).
\hspace{0.4cm}In another form, this is the Theorem, I derived directly  from the principle of the virtual velocities $\S.\,5$.\,eq.\ (\hyperlink{eq1}{1}).\hspace{0.4cm}One can easily show that from this result one  arrives at the  introduction  of the functions $a$ for the case  of the stationary motion in a simpler way  than  \mbox{C\ws l\ws e\ws b\ws s\ws c\ws h} does in the latter essays\,; for a lack of space, we are not giving a more accurate argumentation here.
}\label{footnoteClebsch}\\[1cm]
\centerline {\fett{$\S.\,12$}}
\\[0.3cm]
If one puts  into  the equations $(3)$ of  $\S.\,6$, $t=0$, one finds\,:
\be \tag*{[12.1]}
\frac{d v_0}{dc} - \frac{d w_0}{db}  = 2A, \,\,\,\, \frac{d w_0}{da} - \frac{d u_0}{dc}  = 2B,\,\,\,\, \frac{d u_0}{db} - \frac{d v_0}{da}  = 2 C
\ee
where $A,\,\,B,\,\,C$ are the rotational velocities with respect to the coordinates axes in the point  $a,\,\,b,\,\,c$ at time $t=0$.\hspace{0.4cm}Quite similarly, at time $t$ we have\,:
\begin{equation*} \tag*{(1), [12.2]}
\frac{d v}{dz} - \frac{d w}{dy}  = 2X, \,\, \frac{d w}{dx} - \frac{d u}{dz}  = 2Y,\,\,\, \frac{du}{dy} - \frac{d v}{dz}  = 2 Z
\end{equation*}
These are the equations cited at page 34 [now page~\pageref{pointer2}] found by  \mbox{C\ws a\ws u\ws c\ws h\ws y,}  where $X,\,Y,\,Z$ mean the rotational velocities with respect to the coordinates $x,\,\,y,\,\,z$  at time $t$ taken as axes.\hspace{0.4cm}If $A,\,B,\,C$ are zero overall, i.e.\ at the beginning of the motion,  there are no rotating particles  and thus
$X,\,Y,\,Z$ will stay zero and hence, one can write
\begin{equation*} \tag*{(2), [12.3]}
u= \frac{dF}{dx}, \,\, v= \frac{dF}{dy}, \,\, w = \frac{dF}{dz}.
\end{equation*}
Restricting ourselves to liquid flows, by this substitution  the density equation (6) of $\S.\,2.$ becomes\,:
\begin{equation*} \tag*{(3), [12.4]} 
 \frac{d^2 F}{dx^2} + \frac{d^2 F}{dy^2} +  \frac{d^2 F}{dz^2} =0\,;
\end{equation*}
therefore, 
$F$ satisfies \mbox{L\ws a\ws p\ws l\ws a\ws c\ws e}'s differential equation for all fluid particles\,;
for this reason the function $F$ has been designated by \mbox{H\ws e\ws l\ws m\ws h\ws o\ws l\ws t\ws z}  as the   \mbox{v\ws e\ws l\ws o\ws c\ws i\ws t\ws y} \mbox{p\ws o\ws t\ws e\ws n\ws t\ws i\ws a\ws l.}\hspace{0.4cm} 
From the equation~(3) derives an interesting analogy between the fluid motions in a simply connected space and the effects  of magnetic masses\,: 
the velocities are equal and aligned to the forces exerted by a certain distribution of magnetic masses on the surface on a magnetic particle in the interior. \\
\indent \indent
If the function $F$ is determined by~(3) 
by suitable boundary conditions on the surface,
it is still necessary to determine the pressure $p$ in order to get a complete solution to the problem.
The required equation for that is easily obtained, 
by using the values of $u, \, v, \, w$ from~(2) 
 in~(4) of $\S.\,4$.
Then, one obtains\,:
\begin{equation*} \tag*{[12.5]}
\left.
\begin{aligned} 
\frac{d^2 F}{dt dx} + \frac{d^2 F}{ dx^2} \frac{dF}{dx }+  \frac{d^2 F}{ dx dy} \frac{dF}{dy  } + \frac{d^2 F}{ dx dz} \frac{dF}{dz } - \frac{d V}{dx} + \frac{1}{\varrho} \frac{dp}{dz} = 0\\\\
\frac{d^2 F}{dt dy} + \frac{d^2 F}{ dx dy} \frac{dF}{dx }+  \frac{d^2 F}{ dy^2} \frac{dF}{dy  } + \frac{d^2 F}{ dy dz} \frac{dF}{dz } - \frac{d V}{dy} + \frac{1}{\varrho} \frac{dp}{dy} = 0\\\\
\frac{d^2 F}{dt dz} + \frac{d^2 F}{ dz dx} \frac{dF}{dx }+  \frac{d^2 F}{ dz dy} \frac{dF}{dy  } + \frac{d^2 F} {dz^2} \frac{dF}{dz } - \frac{d V}{dz} + \frac{1}{\varrho} \frac{dp}{dz}=0
 \end{aligned} \qquad \right\} 
\end{equation*}
Through multiplication with ${\rm d} x, {\rm d} y, {\rm d} z$ and 
summation, one gets\,:
\be \tag*{[12.6]}
 {\rm d} \left(\frac{d F}{dt} \right) + \frac{1}{2}  {\rm d}\left( \frac{dF}{dx}\right)^2 + \frac{1}{2}  {\rm d} \left( \frac{dF}{dy}\right)^2 + \frac{1}{2}  {\rm d} \left( \frac{dF}{dz}\right)^2 -  {\rm d} \Omega = 0
\ee
where $\Omega$ is the function defined on S.\,18 [eq.\,\ref{pointerOmega}, now on page~\pageref{pointerOmega}].\hspace{0.4cm}
It follows herefrom by integration\,:
\be \tag*{[12.7]}
\frac{d F}{dt}  +\frac{1}{2} \left[ \left( \frac{dF}{dx}\right)^2 + \left( \frac{dF}{dy}\right)^2   + \left( \frac{dF}{dz}\right)^2  \right] - \Omega =0 
\ee
where the additive integration constant, which is a function of $t$, can be included in $F$.
\hspace{0.4cm}
Apart from the unknown $p$ the known function $V$ is contained in $\Omega$, from which one can easily determine $p$ once $F$ is obtained. \\[0.5cm]
\centerline{\fett{$\S.\,13.$}}

\vspace{0.3cm}
\indent\indent
From the relations $(1)$ of $\S.\,12$, 
\begin{equation*} \tag*{(1), [13.1]}
\frac{d v}{dz} - \frac{d w}{dy}  = 2X, \,\,\frac{d w}{dx} - \frac{d u}{dz} =2Y,\,\,\, \frac{d u}{dy} - \frac{d v}{dz}  = 2 Z 
\end{equation*}
and from the density equation $(6)$ of $\S.\,2$ for liquid 
fluids
\be \tag*{[13.2]}
 \frac{du}{dx} + \frac{dv}{dy} +  \frac{d w}{dz} =0 \,,
\ee
$u,\,v,\,w$ 
{can be determined
as functions of $X,\,Y,\,Z$.
\hspace{0.4cm} Indeed, one finds easily from these equations\,:
\be \tag*{[13.3]}
 \begin{aligned}
\frac{d^2 u}{d x^2}\,\, + \frac{d^2 u}{d y^2}\, + \frac{d^2 u}{d z^2}\,  =  2 \left( \frac{d Z}{dy}\,- \frac{d Y}{dz}\right)\\\\
\frac{d^2 v}{d x^2}\,\, + \frac{d^2 v}{d y^2}\, + \frac{d^2 v}{d z^2}\,  =  2 \left( \frac{d X}{dz}\, - \frac{d Z}{dx}\right)\\\\
\frac{d^2 w}{d x^2} + \frac{d^2 w}{d y^2} + \frac{d^2 w}{d z^2}\,  =  2 \left( \frac{d Y}{dx} - \frac{d X}{dy}\right)
\end{aligned} 
\ee
The integration of these partial differential equations follows from known theorems\,: $u,\,v,\,w$ appear as potential functions of fictitious, attracting masses, 
which are distributed through the fluid-filled space with the density,
\be \tag*{[13.4]}
-\frac{1}{2 \pi} \left( \frac{d Z}{dy} - \frac{d Y}{dz}\right), \,\, -\frac{1}{2 \pi} \left( \frac{d X}{dz} - \frac{d Z}{dx}\right),\,\, -\frac{1}{2 \pi} \left( \frac{d Y}{dz} - \frac{d X}{dy}\right)  \,.
\ee
Denoting by $u_1,\, v_1, \, w_1$ the velocity components at point $x_1, \, y_1,\, z_1$, and by $r$ the distance of this point to $x,\,y\,z$, one has\,:
\be \tag*{[13.5]}
 \begin{aligned}
u_1 =\frac{1}{2 \pi} \iiint \left( \frac {dY}{dz} - \frac{dZ}{dy}\right) \frac{{\rm d}x {\rm d}y {\rm d}z}{r}\\\\
v_1 =\frac{1}{2 \pi} \iiint \left( \frac {dZ}{dx} - \frac{dX}{dz}\right) \frac{{\rm d}x {\rm d}y {\rm d}z}{r}\\\\
w_1 =\frac{1}{2 \pi} \iiint \left( \frac {dX}{dy} - \frac{dY}{dx}\right) \frac{{\rm d}x {\rm d}y {\rm d}z}{r}
\end{aligned} 
\ee
where the integrals have to be extended over all points of the continuous fluid.\hspace{0.4cm}These are not yet the most general values of $u_1, \, v_1, \, w_1$\,: 
since these values  $u_1, \, v_1, \, w_1$ are determined, 
 one can always add to them successively
\be \notag
\begin{aligned}
\frac{dP_1}{dx_1},\,\,\,\frac{dP_1}{dy_1}, \,\,\,\frac{dP_1}{dz_1} \,,
\end{aligned}
\ee
while equations~(1) are still satisfied 
as well as the density equation, provided that\,:
\be \tag*{[13.6]}
\frac{d^2 P}{dx^2} + \frac{d^2P}{dy^2} + \frac{d^2P}{dz^2} =0
\ee
is assumed for all points of the fluid.\hspace{0.4cm}One can consider here $P$ as the potential function of attracting masses which are outside of the space filled with the fluid, and must be dermined so that the conditions for $u_1,\, v_1, \, w_1$ on the  fluid's surface are satisfied.\\\\
\indent\indent
The values found  for $u_1,\, v_1, \, w_1$   in this way can be transformed through integration by parts.\hspace{0.4cm}Since\,:
\be \tag*{[13.7]}
r^2=\left(x -x_1 \right)^2 + \left(y - y_1 \right)^2+\left(z -z_1 \right)^2
\ee 
one finds\,:
\be \tag*{[13.8]}
 \begin{aligned}
\iiint \frac{dY}{dz} \frac{{\rm d}x\,{\rm d}y\, {\rm d}z}{r}= \iint Y \frac{{\rm d}x\, {\rm d}y}{r} + \iiint Y \frac {z-z_1}{r^3} {\rm d}x\, {\rm d}y\, {\rm d}z\\\\
\iiint \frac{dZ}{dy} \frac {{\rm d}x\, {\rm d}y\, {\rm d}z}{r}= \iint Z \frac{{\rm d}x\, {\rm d}z}{r} + \iiint Z \frac {y-y_1}{r^3} {\rm d}x\, {\rm d}y\,
 {\rm d}z
\end{aligned} 
\ee
and thus
\begin{equation*}\tag*{(2), [13.9]}
\left.
\begin{aligned} 
u_1 &=\frac{dP_1}{dx_1} +  \frac{1}{2\pi} \int \left( Y \cos \gamma - Z \cos \beta \right ) \frac{{\rm d} \omega}{r} +  \frac{1}{2\pi} \iiint \left\{ Y(z-z_1) - Z (y - y_1) \right\}\frac{{\rm d}x\, {\rm d}y\, {\rm d}z }{r^3} \,\,\,\,\\\\
&\text{\qquad\qquad\qquad\qquad and in the same way\,:} \\\\
v_1 &=\frac{dP_1}{dy_1} +  \frac{1}{2\pi} \int ( Z \cos \alpha - X \cos \gamma ) \frac{{\rm d} \omega}{r} +  \frac{1}{2\pi} \iiint \left\{ Z(x-x_1) - X (z - z_1) \right\}\frac{{\rm d}x\, {\rm d}y\, {\rm d}z }{r^3}\,\,\,\, \\\\
w_1 &=\frac{dP_1}{dz_1} +  \frac{1}{2\pi} \int ( X \cos \beta - X \cos \alpha  ) \frac{{\rm d} \omega}{r} +  \frac{1}{2\pi} \iiint \left\{ X(y-y_1) - Y (x - x_1) \right\} \frac{{\rm d}x\, {\rm d}y\, {\rm d}z }{r^3} 
\end{aligned} \right\}
\end{equation*}
where ${\rm d}\omega$ denotes a surface element and $\alpha,\beta, \gamma$ the angle formed by the normal to it with the coordinates axes.\\
\indent\indent
These analytical formulas for the representation of $u_1,\,v_1,\,w_1$ in terms of $X,\,Y,\,Z$  allow for an interesting interpretation leading to a striking analogy
between the effect of  vortex filaments and that of electrical currents.\hspace{0.4cm}Namely, if we indicate the parts of $u_1,\,v_1,\,w_1$ which originate from the elements ${\rm d}x\, { \rm d}y\, {\rm d}z$ of the triple integrals with ${\rm d}u_1 { \rm d}v_1 {\rm d}w_1$ then  one has\,:
\be \tag*{[13.10]}
 \begin{aligned}
{\rm d} u_1 =  \frac{1}{2\pi} \left \{ Y(z - z_1, -Z(y-y_1)\right \} \frac{{\rm d}x\, {\rm d}y\, {\rm d}z }{r^3}\\\\
{\rm d} v_1 =  \frac{1}{2\pi} \left \{ Z (x - x_1, -Z(z-z_1)\right \} \frac{{\rm d}x\, {\rm d}y\, {\rm d}z }{r^3}\\\\
{\rm d} w_1 =  \frac{1}{2\pi} \left \{ X (y - y_1, -Y(x-x_1)\right \} \frac{{\rm d}x\, {\rm d}y\, {\rm d}z }{r^3}
\end{aligned} 
\ee
By known theorems, from these equations one derives the relations\,:
 \begin{align}
  \tag*{[13.11]} &(x - x_1)\, {\rm d} u_1 + (y - y_1) \,{\rm d} v_1 + (z - z_1)\, {\rm d} w_1=0,   \\ \nonumber \\
&X {\rm d}u_1 + Y{\rm d} v_1 + Z {\rm d} w_1 =0,\hspace{6.8cm} \tag*{[13.12]}\\\nonumber\\
&{\rm d}u_1^2 + {\rm d}v_1^2 + {\rm d}w_1^2 = \left(\frac{{\rm d}x {\rm d}y {\rm d}z}{2 \pi r^3}\right)^2 \Big\{ (X^2 + Y^2 + Z^2) [(x-x_1)^2 +  (y-y_1)^2 + (z - z_1)^2] \nonumber\\
&\quad\hspace{5cm} - [X(x-x_1)+ Y(y-y_1) + Z(z-z_{3})]^2 \Big\} \,. \tag*{[13.13]}
\end{align} 
\noindent The first two equations show that the velocity with the components  $ {\rm d}u_1, {\rm d}v_1, {\rm d}w_1$ 
\be \tag*{[13.14]}
{\rm d} U_1 = \sqrt{{\rm d}u_1^2 + {\rm d}v_1^2 + {\rm d}w_1^2 }
\ee
is normal to the plane containing the point  $x_1, y_1, z_1$  and
 the rotation axis of the fluid particles $x,\, y,\,z$.\,\,The previous equation can also be written\,:
\begin{equation*} \tag*{[13.15]}
{\rm d} U_1^2=\Big(\frac{{\rm d}x\, {\rm d}y\, {\rm d}z}{2 \pi r^2} \Delta \Big)^2 \left\{ 1-\Big[\frac{X}{\Delta} \frac{x-x_1}{r } +\frac{Y}{\Delta} \frac{y-y_1}{r } + \frac{Z}{\Delta} \frac{z-z_1}{r }\Big]^2 \right\}
\end{equation*}
where $\Delta$ means the rotational velocity.\hspace{0.4cm}Herefrom one finds\,:
\be \tag*{[13.16]}
{\rm d} U_1= \frac{{\rm d}x\,{\rm d}y\,{\rm d}z}{2 \pi r^2} \Delta \sin \varepsilon
\ee
where $\varepsilon$ indicates the angle between the rotation axis of the particle $x, \,y, \,z$ and the connecting line $r$ between this point and $x_1, y_1, z_1$.\hspace{0.4cm}Each rotating particle  $x,\,y,\,z$  generates, then, into another particle of the fluid $x_1, y_1, z_1$ a velocity which is normal to the plane passing through the rotation axis of the particle $x,\,y,\,z$ and $x_1, y_1, z_1$. This velocity is  directly proportional to the volume ${\rm d}x\, {\rm d}y\, {\rm d}z$ of the particle $x,\,y,\,z$, to its rotational velocity $\Delta$ and to the sinus of the angle $\varepsilon$ between the rotation axis of $x,\, y,\, z$ and the connecting line $r$ of both particles, 
is inversely proportional to the square of the distance $r$
 between both particles.\\
\indent\indent
However, according to  A\ws m\ws p\ws \`e\ws r\ws e's Law, this is the same force that an electrical particle of intensity $\Delta$
located at $x,\,y,\,z$ in a current aligned parallely to the rotational axis,  would exert on a little magnet located at $x_1,\,y_1,\,z_1$.\\
\indent\indent
This highly remarkable analogy, whose discovery is due  
to \mbox{H\ws e\ws l\ws m\ws h\ws o\ws l\ws t\ws z,} is firstly
 of great importance for the theory of the vortex filaments in  liquid fluids --- since equations~(2) are only applicable to such liquid 
 fluids. Indeed, it allows to apply theorems developed in electrodynamics, with minor modifications, to hydrodynamics, making the visualisation significantly simpler. 
Furthermore, this analogy is also of a certain value for electrodynamics, since it allows to analyse the electrodynamical processes not based on the mutual elementary interactions between two particles
 --- as is it is usually done following \mbox{A\ws m\ws p\ws \`e\ws r\ws e's}  procedure --- but by considering infinitely thin closed currents  as a basis for the whole theory,  in analogy with vortex filaments}. \\[1cm]
\centerline{\fett{$\S.\,14.$}}\vspace{0.3cm}
\indent\indent
The principle of virtual velocities for
the motion of liquid 
fluids\:
\begin{equation*}\tag*{(1), [14.1]}
\frac{d \delta x}{dx} + \frac{d \delta y}{dy}  + \frac{d \delta z}{dz} = 0,
\end{equation*}
gives the 
next equation, 
following from (1) of $\S.\,4$\,:
\be \tag*{[14.2]}
\iint \Big[(\frac{d^2 x}{dt^2} -X) \delta x + (\frac{d^2 y}{dt^2} -Y) \delta y + (\frac{d^2 z}{dt^2} -Z) \delta z \Big] {\rm d}x\, {\rm d}y\,{\rm d}z =0 .
\ee
\noindent If one takes, instead of the virtual velocities, the actual velocities, one has  to put, instead of $\delta x, \delta y, \delta z$\,:
\be \tag*{[14.3]}
{\rm d}x = \frac{dx}{dt}{{\rm d}t}, \,\,\, {\rm d}y=\frac{dy}{dt}{{\rm d}t}, \,\,\,{\rm d}z=\frac{dz}{dt}{{\rm d}t},
\ee
and, since (1), turns into the continuity equation, when limiting oneself to liquid 
flow,  one finds, by  the usual assumptions about $X,Y,Z$\,:
\be \tag*{[14.4]}
\iiint \Big(\frac{d^2 x}{dt^2}\frac{d x}{dt} +\frac{d^2 y}{dt^2}\frac{d y}{dt} +\frac{d^2 z}{dt^2}\frac{d z}{dt}\Big ) {\rm d} x\,{\rm d} y\,{\rm d}z = \iiint
 \Big( \frac{dV}{dx} \frac{dx}{dt} + \frac{dV}{dy} \frac{dy}{dt} + \frac{dV}{dz}\frac{dz}{dt}\Big) {\rm d} x\,{\rm d} y\,{\rm d}z 
\ee
From this equation, after integration in time, denoting by

\be \tag*{[14.5]}
K=\frac{1}{2}\iiint \Big\{\Big(\frac{dx}{dt}\Big)^2 + \Big(\frac{dy}{dt}\Big)^2 + \Big(\frac{dz}{dt}\Big)^2 \Big\}{\rm d} x\,{\rm d} y\,{\rm d}z 
\ee
the living force 
\be \tag*{[14.6]}
K = \text{{\em Const.}} + \iiint V {\rm d}x\,{\rm d}y\,{\rm d}z
\ee
wherein the constant is
independent of $t$.\hspace{0.4cm}This equation of the living force 
has been presumably given for the first time  by  \mbox{L\ws e\ws j\ws e\ws u\ws n\ws e-D\ws i\ws r\ws i\ws c\ws h\ws l\ws e\ws t}.\footnote{Untersuchungen \"uber ein Problem der Hydrodynamik. Crelle's Journal,  Bd. 58, S. 202, [\hyperlink{Dirichlet}{1860}].}\\
\indent\indent
Consequently,
the living force will be independent of time, provided that the integral $\iiint V {\rm d}x\,{\rm d}y\,{\rm d}z$ is constant in time, 
in other words, 
 since  $V$ for each point in space is independent of $t$.\hspace{0.4cm}The condition that the fluid always fills the same absolute space may be also expressed by  the fact that the normal velocity components of the particles at the surface is zero.\hspace{0.4cm}One 
may show directly that, in this case, $K$ is constant. 
Namely, one has\,:
\be \tag*{[14.7]}
\frac{d K}{dt }= \iiint \Big(\frac{dV}{dx} u + \frac{dV}{dy}v + \frac{dV}{dz}w\Big) {\rm d } x\,\,{\rm d} y\,\,{\rm d}z 
\ee
and because of
\be \tag*{[14.8]}
 \begin{aligned}
\iiint \frac{dV}{dx} u {\rm d}x\,{\rm d}y\,{\rm d}z &= \iint Vu{\rm d}y\,{\rm d}z  - \iiint V \frac{du}{dx} {\rm d}x\,{\rm d}y\,{\rm d}z  \\\
\iiint \frac{dV}{dy} v {\rm d}x\,{\rm d}y\,{\rm d}z  &= \iint Vv{\rm d}z\,{\rm d}x - \iiint V \frac{dv}{dy} {\rm d}x\,{\rm d}y\,{\rm d}z\\\
\iiint \frac{dV}{dz} w {\rm d}x\,{\rm d}y\,{\rm d}z &= \iint Vw{\rm d}x\hspace{0.01cm}{\rm d}y  - \iiint V \frac{dw}{dx}{\rm d}x\hspace{0.02cm}{\rm d}y\hspace{0.02cm}{\rm d}z
\end{aligned} 
\ee
\noindent one has\,:
\be \tag*{[14.9]}
\frac{d K}{dt} = \iint Vu {\rm d}y\,{\rm d}z+\iint Vv {\rm d}z\,{\rm d}x  + \iint Vw{\rm d}x {\rm d}y   - \iiint V\Big( \frac{du}{dx} + \frac{dv}{dy} +\frac{dw}{dz}\Big) {\rm d}x\, {\rm d}y \,{\rm d}z 
\ee
The density equation cancels out the last integral  and one finds that 
\be \tag*{[14.10]}
\frac{d K}{dt}  =\int V\,U_n\, {\rm d} \omega,
\ee
where $U_n$ denotes the normal velocity of a particle at the external 
 surface, and where
 ${\rm d} \omega$ is the surface element.\hspace{0.4cm}If $U_n =0$ on the whole surface,   we find indeed that $\displaystyle \frac{dK}{dt}=0$.\\\\
\indent\indent
For every stationary motion the normal velocity component at the surface vanishes\,; the same happens when  the fluid is surrounded by motionless rigid walls\,;  
this includes also the case when the liquid mass extends to infinity in all directions, since this amounts to having the flow contained in an infinitely
large sphere.\hspace{0.4cm}In all these cases the living forces are constant in time.\\\\
\indent\indent
When there is no rotational motion in the flow,  according to  $\S.\,12$, one can  set\,:
\be \tag*{[14.11]}
u = \frac{dF}{dx}, \hspace{1cm} v = \frac{dF}{dy}, \hspace{1cm} w = \frac{dF}{dz}
\ee
so that we obtain 
\be \tag*{[14.12]}
K= \frac{1}{2} \iiint \Big\{ \Big(  \frac{dF}{dx}  \Big)^2 +  \Big(  \frac{dF}{dy}  \Big)^2 + \Big(  \frac{dF}{dz}  \Big)^2 \Big\} {\rm d}x\,\,{\rm d}y\,\,{\rm d}z
\ee
\indent\indent
By a well-known theorem, one finds\,:%
\deffootnotemark{\textsuperscript{[T.\thefootnotemark]}}\deffootnote{2em}{1.6em}{[T.\thefootnotemark]\enskip}%
\footnote{The factor $F$ in front of the Laplacian in the right-most term was omitted by Hankel.}
\begin{smallequation} \tag*{[14.13]}
 \iiint \Big\{ \Big(\frac{dF}{dx}\Big)^2 + \Big( \frac{dF}{dy}\Big)^2 + \Big( \frac{dF}{dz}\Big)^2 \Big\}{\rm d}x\,{\rm d}y\,{\rm d}z = -\int F\frac{dF}{dn} {\rm d} \omega -
  \iiint  F \big( \frac{d^2 F}{d x^2} + \frac{d^2 F}{d y^2}+ \frac{d^2 F}{d z^2} \big) {\rm d}x\,{\rm d}y\,{\rm d}z
\end{smallequation}
But, by (3) of $\S.\,12$, the last integral vanishes, and thus one has\,:
\be \tag*{[14.14]}
K=- \int F\frac{dF}{dn} \, {\rm d} \omega
\ee
Furthermore, if the flow is in a stationary motion,  
%
then {\small $\displaystyle\frac{dF}{dn}$}, being the  velocity component normal to the external surface, must vanish, 
and thus  $K=0$ from which follows\,: 
\be \tag*{[14.15]}
 \frac{dF}{dx}=0, \hspace{1cm} \frac{dF}{dy} = 0, \hspace{1cm}  \frac{dF}{dz} = 0
\ee
so that there is no motion at all.\hspace{0.4cm}
Thus, a motion  
driven by a velocity potential can never become stationary and conversely, a stationary motion will never have a velocity potential\;--- an interesting theorem discovered
by \mbox{H\ws e\ws l\ws m\ws h\ws o\ws l\ws t\ws z.}
\vspace{1.2cm}
\\\
\centerline{C\ws o\ws r\ws r\ws e\ws c\ws t\ws i\ws o\ws n\ws s [already implemented] }

\vspace{0.5cm}
\noindent
p.\,20.\,\, l.\,7 must be $\big(\frac{dx}{d \varrho_1} \big)^2$ instead of $\big(\frac{dx}{d \varrho_1} \big)^3$ \\
p.\,30.\,\, l.\,10 in front of $2r \cos \theta \cos \varphi \frac{d \theta}{dt} \frac{d \varphi}{d t}$  there must be a $+$ sign  

\vspace{0.8cm}

\section{Letters exchanged for the prize assignment:
judgements of  Bernhard Riemann  and  Wilhelm Eduard Weber} 
\label{s:letters}

\subsection*{{\bf Judgement on Hankel's manuscript by Bernhard Riemann}}

Decision on  the  received  manuscript on the mathematical prize question  carrying the  motto: \,\,{\em The more signs express relationships  in nature, the more useful they are}:\footnote{This motto is in Latin in the  \textit{Preisschrift}: ``Tanto utiliores sunt notae, quanto magis exprimunt rerum relationes''.}
\indent\indent
The manuscript gives  commendable evidence of the Author's diligence, of his knowledge and  ability in using the methods  of computation recently developed  by  contemporary mathematicians.\hspace{0.4cm}This is  particularly shown in  $\S\,6$ of the manuscript, which contains an elegant method for building the equations of motion for a flow in a fully arbitrary  coordinate system, from a point of view which is  commonly referred to as Lagrangian.
 \hspace{0.1cm}
However, when  further developing the general laws of vortex motion, the Lagrangian approach is unnecessarily left aside, and,
 as a  consequence, the various laws had to be found by other and  
quite independent means.\hspace{0.4cm} Also the relation between the equations of motion and  the investigations of Clebsch, reported in $\S.\,14.\,15 $, is omitted  by  the Author.\hspace{0.4cm}Nonetheless, 
as his derivation  actually begins from the Lagrangian equations, one may consider the requirement  posed by the prize-question as fulfilled  by this part of the  manuscript 
(if one substitutes  the  wrong text of a given proof of a  theorem used by the Author with  the  right one contained  in his note 
and leaves out of consideration the mistakes due to rushing  in  $\S.\,9$).\footnote{The $\S\,9$ refers to the Latin manuscript.} \hspace{0.4cm} In the opinion of the referee, 
the imperfectness just evoked in handling  this  part of the question  did not give any sufficient  reason  to deny the prize to the manuscript. However, the Author  will have to  consolidate this part of the manuscript by a reworking, after which the same would  gain in shortness and uniformity. A more important reason against the assignment of the prize could be the incorrectnesses 
occuring in several places. 
These [incorrectnesses] do not get to the core of the argument,
except [in] two [paragraphs] (\S 3 and \S 8) in which the obscurity of another writer may as well serve as an excuse for the Author.
[These incorrectnesses] may still be passed over,
if our highly esteemed Faculty would like  to assign the prize to this manuscript, 
in view of all manner of good things it contains.\\

\begin{figure}[t]
	\includegraphics[width=\textwidth]{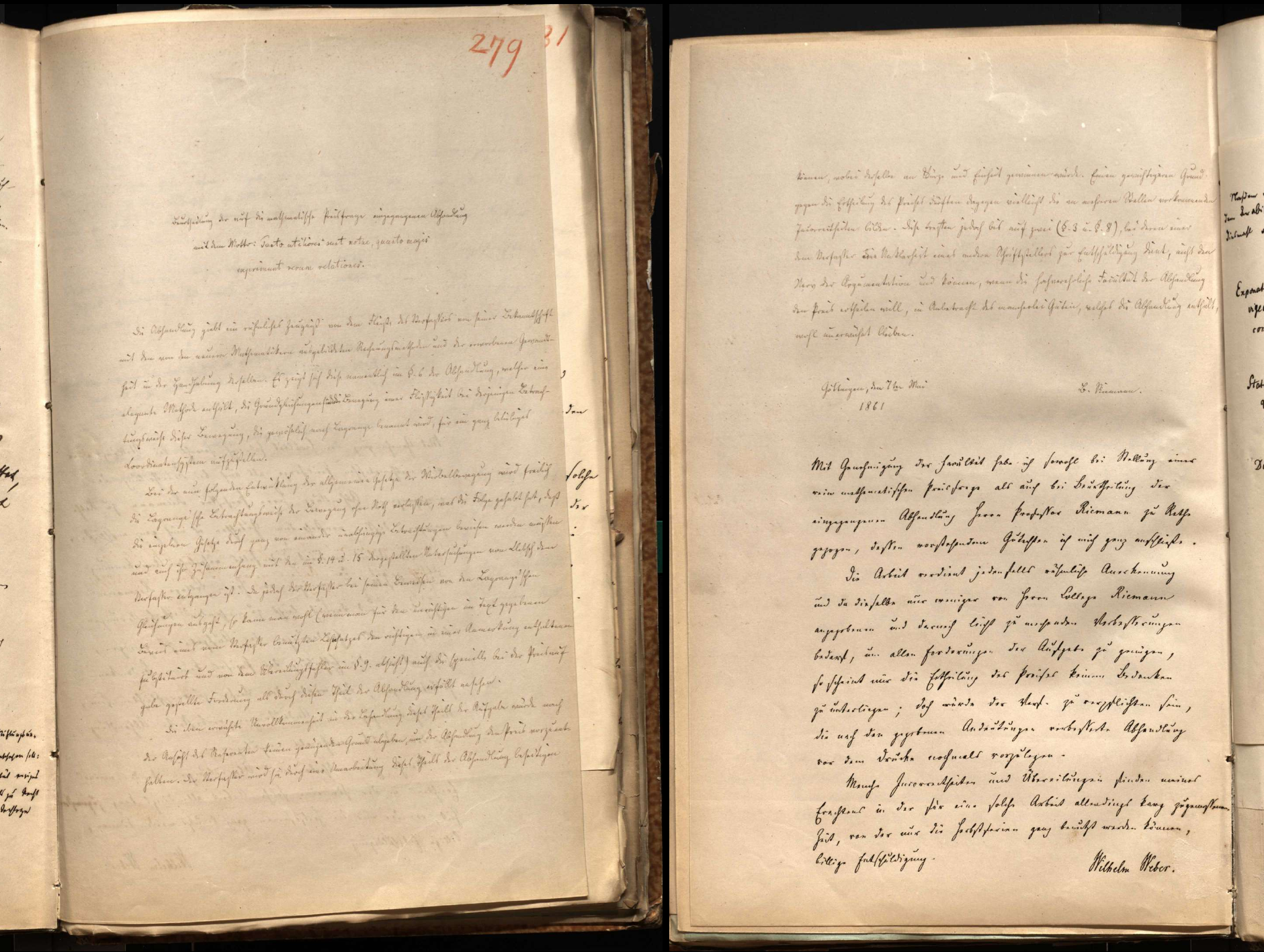}
\caption {Scanned image of the letters from Riemann and Weber, taken from the \textit{G\"ottingen University Archive}. The image is the result of merging two consecutive pages of the full 
text of the exchanged letters amongst the commission members for the prize assignment.
Riemann's letter begins on the left page and ends halfway on the right page. Weber's 
letter follows after Riemann's letter.}
\end{figure}

\vskip-0.3cm
\indent
G\"ottingen, 7th Mai 1861

$\phantom{C}$ \hfill B.\ Riemann \quad \quad
\vspace{0.7cm}

\subsection*{\bf{Judgement on Hankel's manuscript by Wilhelm Eduard Weber}}

With the Faculty approval, I have consulted Prof.\ Riemann, whose judgement  I  
fully agree with, both for the launching
  of a pure mathematical prize   and for the evaluation  of  the  received manuscript. 
In any case, the work deserves
praiseful recognition,
and since it needs  just a few  corrections  indicated by  my colleague Riemann, 
thereby easy to do, in order  to meet the  task requirements, it seems to me that the prize assignment does not give rise to any concern.
Nevertheless, the Author will have to consign again his  revised work before it is  going to print,  according to the given suggestions.
In my consideration, some incorrectnesses and hastinesses find a fair excuse, in the  indeed sparsely proportionated time for such a task,  given that only the autumn  holiday could be used for the scope.\\

$\phantom{C}$ \hfill Wilhelm Weber \quad \quad
\vspace{0.7cm}

\section{Hankel's biographical notes and list of published papers}
\label{s:HHpapers}

\begin{figure}[t]
	\includegraphics[width=0.4\textwidth]{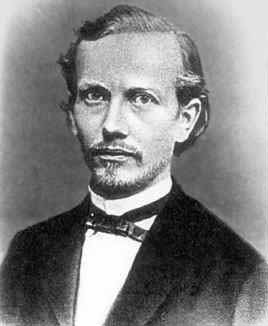}
\caption{Portrait of Hermann Hankel  (from Wikimedia Commons).}
\end{figure}

\subsection*{\bf{Biographical notes about Hermann Hankel}}

\deffootnotemark{\textsuperscript{[T.\thefootnotemark]}}\deffootnote{2em}{1.6em}{[T.\thefootnotemark]\enskip}
Hermann Hankel was born in Halle, near Leipzig, on 14th February 1839. His father was Gottfried Wilhelm Hankel, a renowned physicist.  Hermann Hankel  was  a brilliant student already in high school, with a particular interest in mathematics and its history. From 1857, he studied Mathematics at the Leipzig University under the mathematicians Scheibner, Moebius and Drobisch. Then,  he continued his studies in G\"ottingen, where,   arriving in April 1860,  he could attend, among others, Riemann's lectures. In G\"ottingen he won in 1861 the extraordinary mathematical prize launched
in June 1860  by the Faculty of G\"ottingen with an essay  on the fluid motion theory to be elaborated in a  Lagrangian framework. Also in 1861, he obtained his Doctor degree in Leipzig with the dissertation: ``\"Uber eine besondere Classe der symmetrischen Determinanten''.\footnote{``On a particular class of symmetric determinants''.} Then, in the autumn of the  same year, he went to Berlin, where he could  attend  courses of Weierstrass and Kronecker.
In 1862 he returned to Leipzig and, in 1863, at the same place,
he habilitated as a \textit{Privatdozent} with a thesis
on the Euler integrals with  an unlimited variability of the argument. 
The writing of the habilitation thesis
was probably  firstly induced by the  lectures of Riemann about functions of complex variables. 
In the spring 1867,  he became extraordinary Professor at Leipzig University and, in the same year, ordinary Professor  in Erlangen, then, in T\"ubingen in 1869. He was married to
Marie Dippe,  who much  later became a very important Esperantist. During his life, Hankel was advisor for  doctoral  dissertations in mechanics, real functions and geometry.  He died prematurely on 29th August 1873. Hermann Hankel  is known for his  Hankel functions,  a type  of cylindrical functions,  Hankel transforms, integral transformations whose kernels are Bessel functions of the first kind, and Hankel matrices,  with constant skew diagonals.  Hankel was the first to recognise  the significance of Grassmann's   extension theory (``Ausdehnungslehre''). 
Hankel had a passion for research in  history of mathematics and   published  meaningful writings also in this domain (his inaugural lesson in T\"ubingen was about the development of Mathematics in the last centuries). 
Curiously, his prized work on the fluid-dynamic theory in Lagrangian coordinates written as a student, 
is little known.\footnote{For Hankel's biography  see: Cantor, \hyperlink{Cantor}{1879}; Crowe, \hyperlink{Crowe}{2008} Monna,  \hyperlink{Monna}{1973},  von Zahn, \hyperlink{Zahn}{1874}.}

\vspace{0.7cm}

\subsection*{\bf{List of papers of  Hermann Hankel}}

\begin{enumerate}
 \item[1)] Hankel, Hermann. 1861. \textit {Zur allgemeinen Theorie der Bewegung der Fl\"ussigkeiten.} Eine von der philosophischen Facult\"at der Georgia Augusta am 4.\ Juni 1861 gekr\"onte Preisschrift, G\"ottingen. \HD{http://babel.hathitrust.org/cgi/pt?id=mdp.39015035826760;view=1up;seq=5}{Druck der Dieterichschen Univ.-Buchdruckerei. W.FR.Kaestner, G\"ottingen.}  

 \item[2)] Hankel, Hermann. 1861. \textit {\"Uber eine besondere Classe der symmetrischen Determinanten.} \HD{https://books.google.fr/books/about/Ueber_eine_besondere_Classe_der_symmetri.html?id=vHZaAAAAcAAJ&redir_esc=y}{Inaugural-Dissertation zur Erlangung der philosophischen Doktorw\"urde an der Universit\"at Leipzig von Hermann Hankel.}

 \item[3)]
Hankel, Hermann. 1862. \"Uber die Transformation von Reihen in Kettenbr\"uche. 
\HD{http://gdz.sub.uni-goettingen.de/pdfcache/PPN599415665_0007/PPN599415665_0007___LOG_0026.pdf} {\textit{Zeitschrift f\"ur Mathematik und  Physik},   {\bf 7},  338--343.}
Also  in \textit {Berichte \"uber die Verhandlungen der k\"oniglich s\"achsichen Gesellschaft der Wissenschaften zu Leipzig}, mathematisch-physische Classe {\bf 14}, 17-22, 1862. Verlag der S\"achsischen Akademie der Wissenschaften zu Leipzig.

 \item[4)] Hankel, Hermann,  (signed as Hl.). 1863.  Aufsatz   \"uber  \textit {Ein Beitrag zu den Untersuchungen \"uber die Bewegung eines fl\"ussigen gleichartigen Ellipsoides} by B.Riemann. 
 \HD{https://play.google.com/books/reader?id=zt0EAAAAQAAJ&printsec=frontcover&output=reader&hl=en&pg=GBS.PA50}{\textit{Die Fortschritte der Physik} im Jahre 1861, {\bf 17}, 50--57.}

 \item[5)] Hankel, Hermann, (signed as Hl.). 1863.  Aufsatz \"uber  \textit {Zur allgemeinen Theorie der Bewegung der Fl\"ussigkeiten. Eine von der philosophischen Facult\"at der Georgia Augusta am 4.\ Juni 1861 gekr\"onte Preisschrift, G\"ottingen} by H.\ Hankel. 
\HD{https://play.google.com/books/reader?id=zt0EAAAAQAAJ&printsec=frontcover&output=reader&hl=en&pg=GBS.PA57} {\textit{Die Fortschritte der Physik} im Jahre 1861, {\bf 17}, 57--61.}

\item[6)]
Hankel, Hermann, (signed as Hl.). 1863.  Aufsatz   \"uber  \textit {D\'eveloppements relatifs au \S 3 de recherches de Dirichlet sur un probl\`eme  d'hydrodynamique} by F.\ Brioschi. 
\HD{https://play.google.com/books/reader?id=zt0EAAAAQAAJ&printsec=frontcover&output=reader&hl=en&pg=GBS.PA61}{\textit{Die Fortschritte der Physik} im Jahre 1861, {\bf 17}, 61--62.}

 \item[7)] Hankel, Hermann. 1863. \textit {Die Euler'schen integrale bei unbeschr\"ankter variabilit\"at des Argumentes}: 
\HD{https://ia902205.us.archive.org/28/items/dieeulerschenin00hankgoog/dieeulerschenin00hankgoog.pdf}{zur Habilitation in der Philosophischen Facult\"at der Universit\"at, Leipzig, Voss.}
An extrait is published in \HD{http://gdz.sub.uni-goettingen.de/pdfcache/PPN599415665_0009/PPN599415665_0009___LOG_0006.pdf}{\textit{Zeitschrift f\"ur Mathematik und Physik}, 1864,  {\bf 9}, 1--21.}

 \item[8)]
Hankel, Hermann. 1864. Die Zerlegung algebraischer Functionen in Partialbr\"uche nach den Prinzipien der complexen Functionentheorie. \HD{http://gdz.sub.uni-goettingen.de/dms/load/img/?PID=PPN599415665_0009\%7CLOG_0025} {\textit{Zeitschrift f\"ur Mathematik und  Physik} {\bf 9}, 425--433.}

 \item[9)] Hankel, Hermann. 1864.  Mathematische Bestimmung des Horopters. 
\HD{http://gallica.bnf.fr/ark:/12148/bpt6k15207b/f600.image.r=}{\textit{Annalen der Physik und Chemie}, {\bf 122}, 575--588.}

 \item[10)] Hankel, Hermann. 1864. \textit{\"Uber die Vieldeutigkeit der Quadratur und rectification algebraischer Curven}.
\HD{https://archive.org/stream/bub_gb_M6dWcNSvxJYC\#page/n0/mode/2up}{Eine Gratulationsschrift zur Feier des f\"unfzigjaehrigen Doctorjubilaeums
des Herren August Ferdinand Moebius am 11 Dezember 1864.}

 \item[11)]
Hankel, Hermann. 1867.  Ein Beitrag zur Beurteilung der Naturwissenschaft des griechischen Altertum. \HD{https://babel.hathitrust.org/cgi/pt?id=hvd.hw2918;view=1up;seq=436}{\textit {Deutsche Vierteljahresschrift} {\bf 4}, 120-155.}

 \item[12)]
Hankel, Hermann. 1867.  \textit {Vorlesungen \"uber die complexen Zahlen und ihre Functionen} in \HD{https://books.google.it/books?id=MkttAAAAMAAJ&printsec=frontcover&hl=it&source=gbs_ge_summary_r&cad=0\#v=onepage&q&f=false}{zwei  Theilen, Voss, Leipzig.}

 \item[13)]
Hankel, Hermann. 1867.  Darstellung symmetrischer Functionen durch die Potenzsummen. 
\HD{https://www.digizeitschriften.de/download/PPN243919689_0067/PPN243919689_0067___log8.pdf}{\textit {Journal f\"ur die reine und angewandte Physik}, {\bf 67}, 90--94.}

 \item[14)]
Hankel, Hermann. 1868. Die Astrologie um 1600 mit besonderer R\"ucksicht auf das Verhaeltnis Keppler's und Wallenstein's.   
\HD{http://reader.digitale-sammlungen.de/de/fs1/object/display/bsb10612800_00295.html?contextType=scan&contextSort=score\%2Cdescending&contextRows=10&context=hankel}{\textit {Westermann Monatshefte}, {\bf 25}, 281--294. }

 \item[15)]
Hankel, Hermann. 1869. \textit {Die Entwickelung der Mathematik in den letzten Jahrhunderte}.
 \HD{https://archive.org/details/bub_gb_TE3kAAAAMAAJ}{Ein Vortrag beim Eintritt in den akademischen Senat der Universita\"at T\"ubingen am  29. April 1869 gehalten, Fuessche, T\"ubingen. }

 \item[16)]
Hankel, Hermann. 1869.  Die Entdeckung der Gravitation -- und Pascal - Ein literarisches Bericht von Dr.\ Hermann Hankel. 
\HD{http://gdz.sub.uni-goettingen.de/pdfcache/PPN599415665_0014/PPN599415665_0014___LOG_0014.pdf}{\textit{Zeitschrift  f\"ur Mathematik und Physik}, {\bf 14}, 165--173.}

 \item[17)]
Hankel, Hermann. 1869. Beweis eines Hilfsatzes in der Theorie der bestimmten Integrale.  \HD{http://gdz.sub.uni-goettingen.de/dms/load/img/?PID=PPN599415665_0014\%7CLOG_0030&physid=PHYS_0440}{\textit {Zeitschrift  f\"ur Mathematik und Physik}, {\bf 14},  436--437.}

 \item[18)] Hankel, Hermann. 1869.  Die Cylinderfunctionen erster und zweiter Art. 
\HD{http://gdz.sub.uni-goettingen.de/pdfcache/PPN235181684_0001/PPN235181684_0001___LOG_0041.pdf}{\textit {Mathematische Annalen}, {\bf 1}, 467--501.}

 \item[19)] Hankel, Hermann. 1870. Die Entdeckung der Gravitation durch Newton.  
 \HD{http://reader.digitale-sammlungen.de/de/fs1/object/goToPage/bsb10612803.html?pageNo=496}{\textit {Westermann Monatshefte}, {\bf 27}, 482--493.}

 \item[20)] Hankel, Hermann. 1872. Intorno al volume intitolato: \textit{Geschichte der mathematischen Wissenschaften. 1.\ Theil. Von den \"altesten Zeiten bis Ende des 16.\ Jahrhunderts} of H.Suter.  Relazione del dottor Ermanno Hankel.
\HD{http://gdz.sub.uni-goettingen.de/pdfcache/PPN599471603_0005/PPN599471603_0005___LOG_0022.pdf}{\textit {Bullettino di bibliografia e di storia delle scienze matematiche e fisiche}, {\bf 5}, 297--300.}

 \item[21)] Hankel, Hermann. 1874.  \textit{Zur Geschichte der Mathematik in Althertum und im Mittelalter}, (published \textit{post-mortem}),  
 \HD{http://gallica.bnf.fr/ark:/12148/bpt6k82883t}{Druck und Verlag  von B.G.\ Teubner.}

 \item[22)] Hankel, Hermann. 1875.
\textit{Die Elemente der projectivischen Geometrie in synthetischer Behandlung.} 
  \HD{https://archive.org/stream/dieelementederp00hankgoog\#page/n7/mode/2up}{Vorlesungen von Dr. Hermann Hankel, Teubner, Leipzig.}

 \item[23)] Hankel, Hermann. 1875. Bestimmte Integrale mit Cylinderfunctionen. 
 \HD{http://gdz.sub.uni-goettingen.de/pdfcache/PPN235181684_0008/PPN235181684_0008___LOG_0037.pdf}{\textit {Mathematische Annalen}, { \bf 8}, 453--470. }

 \item[24)] Hankel, Hermann. 1875.
Die Fourier'schen Reihen und Integrale f\"ur Cylinderfunctionen. 
 \HD{http://www.digizeitschriften.de/download/PPN235181684_0008/PPN235181684_0008___log38.pdf}{\textit {Mathematische Annalen}, {\bf 8}, 471--494.}

\item[25)] Hankel, Hermann.  1882 (1870). Untersuchungen \"uber die unendlich oft oszillierenden und unstetigen Funktionen. Abdruck aus dem Gratulationsprogramm der T\"ubinger Universit\"at vom 6. M\"arz 1870. 
 \HD{http://www.digizeitschriften.de/download/PPN235181684_0020/PPN235181684_0020___log14.pdf}{\textit { Mathematische  Annalen}, {\bf 20}, 63--112.}

\item[26)] Hankel, Hermann. 1818--1881, Lagrange Lehrsatz. In  \HD{http://gdz.sub.uni-goettingen.de/pdfcache/PPN358976863/PPN358976863___LOG_0132.pdf} {\textit {Allgemeine Encyclop\"adie der Wissenschaften und K\"unste},  
J.S.\ Ersch, J.G.\ Grube, 353--367.}

\item[27)]
Hankel, Hermann. 1818--1881. Gravitation.  In  \HD{http://gdz.sub.uni-goettingen.de/pdfcache/PPN358787696/PPN358787696___LOG_0301.pdf} {\textit {Allgemeine Encyclop\"adie der Wissenschaften und K\"unste}, 
J.S.\ Ersch, J.G.\ Grube, 313--355.}

\item[28)]
Hankel, Hermann.  1818--1881. Grenze.  In 
\HD{http://gdz.sub.uni-goettingen.de/pdfcache/PPN358976863/PPN358976863___LOG_0132.pdf} {\textit {Allgemeine Encyclop\"adie der Wissenschaften und K\"unste} J.S.\ Ersch, J.G.\ Grube, 185--211.}

\end{enumerate}

\vspace{0.5cm}

\noindent {\it Acknowledgements.} We are grateful to Uriel Frisch and tho the referees for useful remarks. 
We thank the Observatoire de la C\^ote d'Azur and the Laboratoire J.-L.~Lagrange 
for their hospitality.
CR is supported by the DFG through the SFB-Transregio TRR33 ``The Dark Universe''.

{}

\end{document}